\documentclass[11pt,reqno]{article}

\usepackage{amsmath}
\usepackage{amsfonts}
\usepackage{amssymb}
\usepackage{amsfonts} 
\usepackage{latexsym}
\usepackage{graphicx}

\newtheorem{theorem}{\bf Theorem}[section]
\newtheorem{lemma}[theorem]{\bf Lemma}
\newtheorem{cor}[theorem]{\bf Corollary}
\newtheorem{proposition}[theorem]{\bf Proposition}

\newcommand{\intR}{\int^{+\infty}_{- \infty}}
\newcommand{\K}{{\cal K}}
\newcommand{\R}{\mathbb{R}}
\newcommand{\T}{\mathbb{T}}

\newtheorem{definition}[theorem]{\bf Definition}
\newtheorem{remark}[theorem]{\bf Remark}

\newcommand{\re}{\,\mbox{Re}}
\newcommand{\im}{\,\mbox{Im}}

\date{}
\begin{document}

\title{Nonhomogeneous Boundary-Value Problems for
 One-Dimensional Nonlinear Schr\"odinger Equations}

\author{Jerry L. Bona\\
{\small Department of Mathematics, Statistics and Computer Science} \  \\
{\small University of Illinois at Chicago } \\ {\small Chicago, IL
60607, USA} \\ {\small
email: {\tt bona@math.uic.edu}} \\ \\
 Shu-Ming  Sun\footnote{Corresponding author}
\\ {\small  Department of Mathematics} \\ {\small Virginia
Polytechnic Institute and State University}
\\ {\small
Blacksburg, Virginia 24061, USA} \\ {\small email: {\tt sun@math.vt.edu}}\\  \\
Bing-Yu Zhang \\ {\small Department of
Mathematical Sciences} \\ {\small University of Cincinnati}
\\ {\small Cincinnati,
Ohio 45221, USA} \\ {\small  and} \\  {\small Yangtze Center of Mathematics}\\  {\small Sichuan University}\\ {\small Chengdu, China}\\  {\small email: {\tt zhangb@ucmail.uc.edu}}}

\baselineskip=14.3pt

\maketitle



\begin{abstract}
This paper is concerned with initial-boundary-value problems (IBVPs)  for  a
class of nonlinear Schr\"odinger equations posed  either on a half line $\R^+$ or on a
bounded interval $(0, L)$ with nonhomogeneous
boundary conditions. For any $s$ with $0\leq s < 5/2$ and $s \not = 3/2$, it is shown that the relevant IBVPs are
 locally well-posed if the initial data lie in the $L^2$--based Sobolev spaces $H^s(\R^+) $
in the case of the half line and in $H^s (0, L)$ on a bounded interval, provided the
boundary data are selected from $H^{(2s+1)/4}_{loc} ( \R^+)$
and $H^{(s+ 1) /2}_{loc} (\R^+)$, respectively.  (For $s > \frac12$, compatibility
between the initial and boundary conditions is also needed.)   Global
well-posedness is also discussed when $s \ge 1$.  From the point of view of
the well-posedness theory, the results obtained
 reveal a significant difference between the IBVP
posed on $\R^+$ and the IBVP posed on
 $(0,L)$.  The former is reminiscent of the theory for the pure initial-value
problem (IVP) for these Schr\"odinger equations posed on the whole
line $\R$ while the theory on a bounded interval looks more like that of
 the pure IVP posed
on a periodic domain. In particular, the regularity demanded of the
 boundary data for the IBVP on $\R^+$
is consistent with the temporal trace results that obtain
for solutions of the pure  IVP on $\R$, while
the slightly higher regularity of boundary data for the IBVP on $(0, L)$
resembles what is found for temporal traces of spatially periodic solutions.

\end{abstract}

\allowdisplaybreaks

\section{Introduction}

\setcounter{equation}{0}

Studied here are
initial-boundary-value problems for  nonlinear Schr\"odinger
 equations posed either on a half line $\mathbb{R}^+$, {\it viz.}
\begin{equation}\label{half}
\left \{ \begin{array}{l}
           iu_t +u_{xx} + \lambda |u|^{p-2} u =0, \qquad x\in \mathbb{R}^+, \ t\in \mathbb{R
           },\\ \\
           u(x,0) =\phi (x), \qquad
           u(0,t) =h (t),
         \end{array} \right.
         \end{equation}
or on a finite interval $(0,L)$,
 \begin{equation}\label{finite}
    \left \{ \begin{array}{l}
               iu_t+u_{xx}+\lambda |u|^{p-2}u =0, \qquad  \ x\in
               (0,L),\quad t \in \R\, ,
               \\ \\
               u(x,0) =\phi (x),\qquad u(0,t) =h_1 (t), \qquad u(L,t) =h_2 (t).
             \end{array}
             \right.
 \end{equation}
Here, the parameter $\lambda$ is a non-zero real number  and $p\geq 3$.\footnote{Only
the  case  $p \geq 3$ is considered here, but a substantial part of the theory goes through under the weaker hypothesis  $p > 2$.}  Note that, due to the symmetry of
the  equation with respect to the change of variables $x \to -x$, results
established for \eqref{half} carry over {\it mutatis muntandis} to the quarter-plane
problem where $\mathbb{R}^+$  is replaced by $\mathbb{R}^-$.   (The situation
regarding the quarter-plane problems posed on
 $\mathbb{R}^+$ and $\mathbb{R}^-$ for the Korteweg-de Vries equation
are significantly different on the other hand.)
In all cases where \eqref{half} and \eqref{finite} arise in practice,
 the second-order derivative models dispersive effects, which is to say the
tendency of waves to spread out due to the fact that different wavelengths
propagate with different speeds, while the $|u|^{p-2} u$--term accounts for
a variety of nonlinear effects.

Nonlinear Schr\"odinger equations are derived as models for a considerable range of applications.
This includes propagation of light in fiber optics cables, certain types of shallow and
deep surface water waves, Langmuir waves in a hot plasma and in more general forms
in Bose-Einstein condensate theory.  In the case of gravity waves on the surface of
an inviscid liquid, the parameter
$\lambda$ depends upon the undisturbed depth of the water, becoming negative
in water deep with respect to the wavelength of the wavetrain.
A particularly interesting application of nonlinear Schr\"odinger  (NLS henceforth)
equations has been their
use in attempting to explain the somewhat mysterious formation of rogue waves in
the ocean and in optical propagation (see
\cite{bpss},  \cite{BSaut} ,  \cite{CHA2011} and \cite{pere1983}).

In many of the physical applications mentioned above,
the independent variable $x$ is  a coordinate representing
position in the medium of propagation, $t$ is proportional to elapsed
time and $u(x, t)$ is a velocity or an amplitude at the point $x$ at
time $t$.
One configuration that arises naturally
in making predictions of waves in water is to take $x \in  \R^+ = \{ x\,  |\, x \geq 0\} $ and specify
$u(0, t) $ for $t > 0$. This
corresponds to a given wave-train generated at one end by a wave-maker and propagating
into  a region of the medium of propagation (see \cite{hammack} for an example
of this situation).
The domain is natural since solutions of this wavemaker problem for the
 NLS-equation are an approximation
of waves moving in the direction of increasing values of $x $.  The semi-infinite  aspect
of the domain is convenient in that no lateral boundary need be considered downstream of
the wavemaker.

However, real domains are bounded, and in some cases it may be necessary
to impose boundary conditions at both ends of the medium of propagation.
Especially  if one is interested in implementing a numerical scheme to calculate
solutions of the half-line problem or  localized solutions of the
 pure initial-value problem on the whole line, there
arises the need to truncate  the spatial domain.  In such situations, the problem posed
on a finite domain comes to the fore, and one must impose
boundary conditions at both ends  to specify  solutions.  Of course, when
approximating localized solutions of the problem on all of $\R$, it is reasonable to take
$u(0, t) = u(L, t) = 0 $ for $0\le t \le T$ and $ L > 0$ large  enough that essentially
no disturbance reaches the boundary during the time interval $[0,T]$.  However, the wavemaker
problem and its finite domain counterpart demand non-homogeneous boundary
conditions.  Neumann conditions may also be appropriate in some circumstances.

In this paper, the discussion will center around the fundamental questions of existence
and uniqueness of solutions corresponding to specified initial and boundary data.
The issue of the solutions' dependence upon the auxiliary data is also examined, thereby
completing Hadamard's basic idea of well posedness.
The theory developed here will be for initial data in the
$L^2$--based Sobolev spaces $H^s(\R ^+)$ and $H^s(0,L)$.  The spaces from
which the boundary data will be drawn are dictated by these choices of initial data,
as will become apparent presently.    Theory  will be developed wherein the time for
existence depends upon the size of the auxiliary data.  With more restrictive
hypotheses, global well-posedness results will also be provided.  Here and
below, the notation is that which is current in the theory of partial differential
equations.

Theory for the nonlinear Schr\"odinger equation in the form depicted in \eqref{half} and
\eqref{finite} has seen a lot of development in the last four decades, beginning
with the pioneering work of Zakharov and his collaborators \cite{ zm1974, zs1972}.
For the
most part, the mathematical theory for this equation has been concerned with
either the pure initial-value problem posed on the entire real line $\R$ or
the periodic initial-value problem posed on the one-dimensional torus $\T$.
A large body of literature has been concerned with the fundamental questions of
 existence,
uniqueness and continuous dependence of solutions corresponding to initial
data  drawn from Sobolev classes  (again, well posedness {\it a la} Hadamard \cite{hadamard, hadamard1}).
Some highlights of the developments are \cite{bourgain-1,
bourgain-2, cfh2011,
caz-weiss,gini-velo-2,gini-velo,kato-1, kato-2,tsu-1, tsu-4}, for example.  We caution
that this is only a small sample of the extant work on these problems.  The monograph
of Cazenave \cite{cazenave}  provides a good entry into the literature.

The
study of the initial-boundary-value problems (IBVP henceforth)
(\ref{half}) and (\eqref{finite}) with
nonhomogeneous boundary conditions  has not progressed to the same extent (see
\cite{bc,brezis,bu1994,bu2000, bu2001,
bu2005,holmer, holmer-2, kam,strauss,tsu-1,tsu-2,tsu-3} and more recent work on the boundary-value problems of some other dispersive equations \cite{Killip}, and the references
therein).
In this paper, the goal
is to advance the study of the IBVP's (\ref{half}) and (\ref{finite})
to the same level as that obtaining for the relevant pure initial-value problems
posed on all of $\R$.
The local well-posedness theory constructed in the body of the paper
is summarized in the following three theorems.  In all of these results,
we assume that the lowest order compatibility conditions
\begin{equation}   \label{compat}
h(0) = \phi(0) \;\;{\rm for}\; \eqref{half}   \quad {\rm or} \quad
h_1(0) = \phi(0), \; h_2(0) = \phi(L)  \;\; {\rm for} \; \eqref{finite}
\end{equation}
are valid when $s > \frac12$.  These derive simply from the requirement that the solution
be continuous at the corners of the relevant space-time domain, which are
$(0,0)$ for the half-line problem and $(0,0)$ and  $(L,0)$
for the finite interval problem. This point will be elaborated at the end of the next section.
If $ s > 0$ is large, we need to assume that $| u |^{p-2} u$ is differentiable, 
a requirement that imposes a relationship between  $s$ and $p$, {\it viz.}
\begin{align}
&\mbox{if $p$ is even, $s$ is arbitrary;\; if $p$ is odd, $s \leq p-1$; \; otherwise,
 $\lfloor s\rfloor< p-2$ ,}\label{s-assumption-1}
\end{align}
where $\lfloor s\rfloor$ is the largest integer less than $s$. Furthermore, for the convenience of our discussion of the traces of functions in $H^s( \R )$, it is always assumed that
\begin{align}
 s \not = n + \frac12 \qquad \mbox{for}\qquad  n = 0, 1, 2, \cdots. \label{s-assumption-2}
\end{align}
This aspect is not always recalled in the body of the paper.
 
For a given $s\in \R$ and  $\Omega$ being $\R^+$ or a
finite interval $(0,L)$,  the space $H^s (\Omega)$ is defined as the
restriction of the space $H^s (\R)$ to $\Omega $,  {\it viz.}
$$
H^s (\Omega) \, =\,  \{ f|_\Omega \; \big | \; f \in H^s(\R) \}
$$
endowed with the quotient norm 
$$
 \|f\|_{H^s(\Omega)} \, = \,    \inf \Big\{ \|\tilde{f}\|_{H^s(\R)} \,|\,  \tilde{f} \in H^s (\R),  \tilde{f}\big |_{\Omega} =f \Big \}. 
$$
Other equivalent definitions of $H^s (\Omega)$ can be found in Chapter 1 of \cite{LM1972}.
\begin{theorem}\label{th-half}
\begin{itemize}
\item[(i)] Suppose $\frac12 < s < \frac52$ with $3\leq p < \infty $. Then, for $ \phi \in H^s (\mathbb{R}^+)$ and $h\in H^{\frac{2s+1}{4}}_{loc} (\mathbb{R}^+)$, the IBVP (\ref{half}) is locally
well-posed in $H^s (\mathbb{R}^+)$.
\item[(ii)] If $0\leq s< \frac12$ with $3\leq p< \frac{6-4s}{1-2s}$,
the IBVP (\ref{half}) is (conditionally) locally well-posed in $H^s(\R^+)$.
\vspace{.08cm}

In both $(i)$ and $(ii)$, what is meant precisely is
that  for any given $T>0$ and $\gamma >0$, there exists a $T^*$ with $0 <
T^* \leq T$ depending only on $s$, $\gamma$ and $T$ such that if
$\phi \in H^s (\mathbb{R}^+)$ and $h\in H^{\frac{2s+1}{4}}(0,T)$
satisfy
\[\| \phi\|_{H^s (\mathbb{R}^+)}+\|h\| _{H^{\frac{2s+1}{4}}(0,T)} \leq
\gamma, \]
then the IBVP (\ref{half}) admits a  solution $u\in
C([0,T^*]; H^s (\mathbb{R}^+))$.  In case $(i)$, the solution $u\in
C([0,T^*]; H^s (\mathbb{R}^+))$ is unique,
while in $(ii)$, the solution satisfies the auxiliary condition
\begin{equation}
 \|u\|_{L^q(0,T^*;
L^r (\mathbb{R}^+))}< +\infty , \label{condd-1}
\end{equation}
where $(q , r)$ is an admissible pair, and it is the only $C([0,T^*]; H^s (\mathbb{R}^+))$--solution with
this property.   Here, a pair $(q,r)$ is admissible when $\frac1q + \frac{1}{2r} =
\frac14$.   In both cases (i) and (ii), the corresponding solution
map is Lipschitz continuous.
\end{itemize}
\end{theorem}
\begin{theorem}\label{th-finite}
\begin{itemize}
\item[(i)] If $\frac12 < s < \frac52$ with $3\leq p < \infty$,  the IBVP (\ref{finite}) is locally
well-posed in $H^s (0,L)$ for $\phi \in H^s (0,L)$ and $ h_1, \ h_2\in H^{\frac{s+1}{2}}_{loc} (\mathbb{R}^+)$.
\item[(ii)] If  $0\leq s< \frac12$ with $ 3\leq p\leq  4$,
then the IBVP (\ref{finite}) is (conditionally) locally well-posed
in $H^s (0,L)$ for $\phi, h_1, h_2$ in the same spaces.
\vspace{.08cm}

For both cases, this means that for any  $T>0$ and $\gamma >0$, there exists
a $T^*$ with $0< T^* \leq T$ depending only on $s$, $\gamma$ and $T$ such that
if $\phi \in H^s (0,L)$ and $h_1, \ h_2 \in
H^{\frac{s+1}{2}}(0,T)$ satisfy
\[\| \phi\|_{H^s (0,L)}+\|h_1\| _{H^{\frac{s+1}{2}}(0,T)} + \|h_2\| _{H^{\frac{s+1}{2}}(0,T)}\leq
\gamma, \]
the IBVP (\ref{finite}) admits a solution $u\in
C([0,T^*]; H^s (0,L))$.  In case $(i)$ this solution $u\in
C([0,T^*]; H^s (0,L))$ is unique, while in case $(ii)$, the solution
also satisfies
\begin{equation}
 \|u\|_{L^4((0,T^*)\times (0,L))}< +\infty \label{condd-2}
\end{equation}
and is the unique $C([0,T^*]; H^s (0,L))$--solution with this property.
In both cases, the corresponding solution map is Lipschitz continuous.
\end{itemize}
\end{theorem}

  The issue of uniqueness could use some elaboration. In (i) of both Theorems \ref{th-half} and \ref{th-finite}, the uniqueness means that if there are two solutions $u , v \in C ( [ 0, T^*]; H^s)$, then $u \equiv v$. However, for (ii) of both Theorems \ref{th-half} and \ref{th-finite}, the uniqueness means that if there are two solutions $u , v \in C ( [ 0, T^*]; H^s)$ satisfying either \eqref{condd-1} or \eqref{condd-2}, then $ u \equiv v $. Therefore, when $0\leq s< \frac12$, the local well-posedness
results presented in both Theorem \ref{th-half} and Theorem
\ref{th-finite}  are {\it conditional} (see Kato \cite{kato-1} where this distinction
was made in the context of general classes of equations) since (\ref{condd-1}) or
(\ref{condd-2}) is needed to ensure the uniqueness.
It is naturally of interest to know whether these conditions can be removed.
If these auxiliary conditions can be removed, the corresponding results
are called unconditional well-posedness, or simply well-posedness.
In fact, a further argument allows the results  for smaller values of $s$ to be extended,
 so obtaining the following additional
wrinkle appertaining to Theorems \ref{th-half} and
\ref{th-finite}.

\begin{theorem}[unconditional well-posedness]
Let $0\leq s< \frac12$ be given.
Then, both (\ref{condd-1}) and
(\ref{condd-2}) can be removed, so the corresponding
well-posedness is unconditional.
\end{theorem}
As mentioned, the preceding results are all local, which is to say the time interval
$(0,T^*)$  over which  the solution is guaranteed to exist depends on the size of
initial and boundary data.  If $T^*$ can be chosen independently
of the size of the initial and boundary data, then the result is termed  global well-posedness.
The following global well-posedness results
for (\ref{half}) and (\ref{finite}) are proved here.

\begin{theorem}\label{global well-posedness}
\quad
\begin{itemize}
\item[(i)] Assume that either
$p \geq 3$ if $\lambda < 0$ or $3 \leq p \leq 4$ if $  \lambda >
0.$ The IBVP (\ref{half}) is globally well-posed in
$H^s (\R^+)$ for any $1\leq s < \frac52$ with auxiliary date $(\phi,h)$ drawn from
$ H^s (\mathbb{R}^+) \times H^{\frac{s+1}{2}}_{loc} (\mathbb{R}^+)$.
\item[(ii)] Assume that either
$p \geq 3 $ if $ \lambda < 0$ or $3 \leq p \leq \frac{10}{3}$ if $  \lambda >
0.$ The IBVP (\ref{finite}) is globally well-posed in $H^s (0,L)$ for any $1\leq s < \frac52$ with
$\phi \in H^s (0,L)$ and $ h_1 , \ h_2 \in H^{\frac{s+1}{2}}_{loc} (\mathbb{R}^+)$.
\end{itemize}
\end{theorem}

The rest of the paper is organized as follows.   A general overview of the problems together
with an outline of the strategy for analyzing them is
  provided in Section 2. The IBVP (\ref{half}) takes center stage
 in Section 3 which consists
of three subsections. In Subsection 3.1,  explicit
solution formulas are derived for  associated linear problems. In Subsection
3.2, various Strichartz estimates are established using these solution formulas. The local
well-posedness  of the IBVP (\ref{half}) on $\R^+$ is established in Subsection 3.3.
In Section 4,   local well-posedness for the IBVP (\ref{finite}) on the finite interval $(0, L)$ is studied.
The global well-posedness of (\ref{half}) and (\ref{finite}) will be investigated in Section 5.
The paper concludes with an Appendix where a technical lemma  needed in
establishing Proposition \ref{prop3.6}
is proved
and a telling counterexample, which concerns the
optimality of the assumption $h_1 , h_2 \in H^{(s+1)/2} ( 0 , T ) $ in
Theorem \ref{th-finite},   is presented.

\section{Overview}

We begin by reviewing the state of the art for the pure initial-value problems
\begin{equation} \label{ivp}
i u_t + u_{xx} + \lambda|u|^{p-2}u = 0, \quad  u(x,0) = \phi(x), \;\; {\rm for} \;\; x \in \R,
\end{equation}
for the Schr\"odinger equations considered here.   First discussed is the case
of  initial data $\phi$ that is localized on an unbounded domain, which is to
say it evanesces at infinity in at least a weak sense.

\vspace{.15cm}

\noindent {\bf Theorem A} {\it
\begin{itemize}
\item[(i)] For  $s>\frac12 $ with $3\leq p < \infty $ or $0\leq s <\frac12 $ with $3\leq p<
\frac{6-4s}{1-2s}$, the initial-value problem
(\ref{ivp})
 is locally well-posed in  $H^s
(\mathbb{R})$.  That is,  for any $r>0$, there exists a $T>0$ depending
 on $r$ such that if
$\| \phi \|_{H^s(\R)} \leq r ,$
then (\ref{ivp}) admits a unique solution $u\in C([0,T]; H^s
(\mathbb{R}) ) $  and the corresponding solution map is
Lipschitz continuous.
\footnote{For many years since the pioneering work in \cite{kato-2, caz-weiss}, the solution map was only known to be continuous from $H^s (\R) $ to $C([0,T]; H^{s-\epsilon} (\R))$.
It was proved recently by Cazenave, et al. \cite{cfh2011}  to be continuous from $H^s (\R) $ to $C([0,T]; H^{s} (\R))$.}

Moreover, for $0\leq s <\frac12 $, the solution also satisfies
\begin{equation}
\label{cond-1} \| u\|_{L^q_{loc}(0,T; B^{s}_{r,2} (\mathbb{R}))} <
+\infty\, ,
\end{equation}
where
$B^s_{r,2 } (\R)$ is the Besov space and  $1/q  +  1 / (2 r) = 1/4$.  Uniqueness
when $0 \leq s < \frac12$ requires that \eqref{cond-1} holds.

\item[(ii)] If, in addition, $3\leq p< 6$, then the above local
well-posedness results are global, i.e.,  $T$ is  independent of $r$
and can be chosen arbitrarily large.
\end{itemize}
}

Next, the existing results obtained when $\phi$ is periodic are recalled.

\vspace{.2cm}

\noindent {\bf Theorem B} {\it
\begin{itemize}
\item[(i)]
For $s>\frac12 $ with $3\leq p < \infty $ or $0\leq s <\frac12 $ with $3\leq
p<\frac{6-4s}{1-2s}$, the IVP
(\ref{ivp})
 is locally well-posed in  $
H^s(\T)$, i.e., for any $r>0$, there exists a $T>0$ depending only
on $r$ such that if $ \phi \in H^s (\T)$ with
$\| \phi \|_{H^s ( \T)} \leq r ,$
then (\ref{ivp}) admits a unique solution $u\in C([0,T]; H^s
(\mathbb{T} )) $ and the corresponding solution map is
Lipschitz  continuous. Moreover,
for $0\leq s <\frac12 $,  the solution $u$ satisfies
\begin{equation}
\label{cond-2} \| u\|_{  \mathbb{B}^{T}_{s,\frac12}   } < +\infty\, ,
\end{equation}
where
$\mathbb{B}^{T}_{s,\frac12}$ is the restricted Bourgain space
associated to the Schr\"odinger equation (see \cite{bourgain-1}).
As in Theorem A, uniqueness when $0 \leq s < \frac12$  is conditional
and relies upon \eqref{cond-2}.
\item[(ii)] If, in addition, $3\leq p< 6$, then the above local
well-posedness results are global, i.e.,  $T$ is  independent of $r$
and can be chosen arbitrarily large.
\end{itemize}
}

These results may be found in the previously cited references.
We emphasize that at present, the uniqueness part of the
well-posedness results  in the parts (i)
of Theorem A and Theorem
B  requires  the extra conditions (\ref{cond-1}) and (\ref{cond-2}) when $s < \frac12$.
 As mentioned, such well posedness was termed
\emph{conditional} by Kato \cite{kato-1}.   If the auxiliary
conditions can be removed, which is to say  the solution is shown to be
unique only assuming it lies in  $C([0,T]; H^s (\R))$, then  the problem \eqref{ivp} is said
to be unconditionally well-posed.
According to the general discussion presented in \cite{bsz-4}, if $3\leq p< 6$, the
conditional well-posedness results stated in parts (i) of Theorems A and B
 are, in fact, unconditional.\footnote{The reader is referred  to  \cite{win-tsu, h-f} and the references therein for recent progress on the issue of unconditional
 well-posedness of nonlinear Schr\"odinger equations.}
\newcommand{\inxdomain}{\mathbb{R}^+}
\newcommand{\fixdomain}{(0,L)}

\smallskip

The overall goal of the present essay is to bring the well-posedness theory for the
 IBVP's \eqref{half} and \eqref{finite}
into line with what is known for the pure initial-value problem \eqref{ivp}.

The precise terminology used in the paper is now provided
and motivation is developed for the choice of appropriate function
spaces for the initial and boundary conditions.   The  main ideas and methodology
for proving the
results stated in the Introduction are also set forth.

The  notion of well-posedness used for the problems
 (\ref{half}) and (\ref{finite}) is detailed first.
\begin{definition}
Let $s, \ s'\in \R$  and $T>0$ be given.
\begin{itemize}
\item[(i)]
The IBVP (\ref{half}) is
said to be (locally) well-posed in $H^s (\R^+)\times H^{s'
} (0,T)$ if for $\phi \in H^s (\R^+)$ and $h\in H^{s'}(0,T)$
satisfying certain natural compatibility conditions, there exists a
$T'\in (0,T]$ depending only on $\|\phi \|_{H^s(\R^+)}+
\|h\|_{H^{s'} (0,T)}$ such that (\ref{half}) admits a unique
solution $ u\in C([0,T']; H^s (\R^+)).$
Moreover, the solution depends continuously on
$(\phi ,h)$ in the corresponding spaces.
\item[(ii)]
The IBVP (\ref{finite}) is said to
be (locally) well-posed in $H^s (0, L)\times H^{s' }
(0,T)$ if for $\phi \in H^s (0,L)$ and $h_1, \ h_2 \in H^{s'}(0,T)$
satisfying certain natural compatibility conditions, there exists a
$T'\in (0,T]$ depending only on $\|\phi \|_{H^s((0,L))}+
\|h_1\|_{H^{s'} (0,T)}+ \|h_2\|_{H^{s'} (0,T)}$ such that
(\ref{finite}) admits a unique solution
$ u\in C([0,T']; H^s (0,L))$.
Moreover, the solution depends continuously on
$(\phi, h_1, h_2)$ in the corresponding
spaces.

In either case, if $T'$ can be chosen independently of $r$, the
relevant IBVP is said to be globally well posed.
\end{itemize}
\end{definition}
Completing this definition of well-posedness requires making precise what it means
for $u$ to be a {\sl solution} of (\ref{half}) or (\ref{finite}).
The issue  is important for small values of $s$, where the meaning of the
derivatives and nonlinear term has
to be addressed. The usual approach in the
literature is to say that $u$ solves the equation
in the sense of Schwartz distributions.  This, however,
leads to  a further question about how the nonlinear
term $\lambda |u|^{p-2}u$ makes sense as a
distribution, as well as how the solution $u$  takes on the given initial and boundary
values. In this paper,  we
will use the following definitions (see \cite{bsz-4} for a general discussion)
for the solutions of (\ref{half}) and
(\ref{finite}), respectively.
\begin{definition}
Let $s\leq2$, $s'\leq s$  and $T>0$ be given.
\begin{itemize}
\item[(a)]
For  $\phi \in H^s (\R^+)$ and $h\in H^{s'}(0,T)$, we say that
$u\in C([0,T]; H^s (\R^+))$ is a solution of (\ref{half}) if there
exists a sequence $$u_n \in C([0,T]; H^2(\R^+))\cap C^1([0,T]; L^2
(\R^+)), \ n=1,2,3, \cdots  $$ such that
\begin{itemize}
\item[1)]
$u_n$ satisfies the equation in (\ref{half}) in $L^2 (\R^+)$ for
 $0\leq t\leq T$,
\item[2)]$u_n$ converges to $u$ in  $C([0,T]; H^s (\R^+))$ as $n\to
\infty$,
 \item[3)] $\phi_n (x)= u_n (x,0)$ converges to $\phi (x)$ in  $H^s (\R^+)$ as $n\to \infty
 $,
 \item[4)] $h_n (t)= u_n
(0,t)$  is in $H^{s'}(0,T)$ and converges to $ h(t)$ in $H^{s'}(0,T)$  as $n\to
\infty $.
\end{itemize}

\item[(b)] For  $\phi \in H^s (0,L)$ and $h_1, \ h_2 \in H^{s'}(0,T)$, we say
that $u\in C([0,T]; H^s (0,L))$ is a solution of (\ref{finite}) if
there exists a sequence $$u_n \in C([0,T]; H^2(0,L))\cap C^1([0,T];
L^2 (0,L)), \ n=1,2,3, \cdots $$ such that
\begin{itemize}
\item[1)]$u_n$ satisfies the equation of (\ref{finite}) in  $L^2
(0,L)$ for $0\leq t\leq T$,
\item[2)]$u_n$ converges to $u$ in  $C([0,T]; H^s (0,L))$ as $n\to
\infty$,
\item[3)]$\phi_n (x)= u_n (x,0)$  converges to $\phi (x)$ in  $H^s
(0,L)$ as $n\to \infty$,
\item[4)] $h_{1,n} (t)= u_n (0,t),  h_{2,n} (t)= u_n (L,t)$
are in $H^{s'}(0,T)$ and converge to $h_1 (t)$ and $h_2(t)$, respectively, in
$H^{s'}(0,T)$ as $n\to \infty$.
\end{itemize}
\end{itemize}
\end{definition}

Of course, if $s \geq 2$, then a solution in the above sense, sometimes called a mild solution,
is a solution in the ordinary $L^2$--sense.

Attention is now turned to the  relation between $s'$ and $s$ in the
definition of well posedness.
It is well known that the linear Sch\"odinger
equation
\[ iv_t +v_{xx} =0, \quad v(x,0)=\phi (x),   \]
posed on the whole line $\R$
 has the  Kato smoothing property, which is to say $\phi \in H^s (\R)$ implies $v\in L^2_{loc} (\R;
H^{s+\frac12}_{loc} (\R))$. In addition, the Schr\"odinger equation itself entails that
$ \partial _t \sim \partial _{xx}$,  so suggesting that
\begin{equation} \label{optimal}
s'=  \frac12\left( s + \frac12\right) = \frac{2s+1}{4}
\end{equation}
 (see \cite{holmer} for a more detailed discussion and \cite{Audiard-1,Audiard-2} for recent studies of this issue for Schr\"odinger equations).
We are thus led to complete the definition of well posedness
for the IBVP (\ref{half}) with the stipulation \eqref{optimal}.

For the IBVP (\ref{finite}), one might  imagine that the correct value should also be
$s'=\frac{2s+1}{4}.$  However, as  will be seen presently,
 this is not the case.  Instead, the optimal relation
between $s$ and $s'$ for the IBVP (\ref{finite}) is
\begin{equation}
\label{optimal-1} s'= \frac{s+1}{2}.
\end{equation}
Thus, a significant, albeit technical difference, emerges between the IBVP (\ref{half})
 (posed on an unbounded domain)  and the IBVP
(\ref{finite}) (posed on a finite domain).

To put our main theorems into context, we sketch previous work on such
IBVP's.  Carrolle and Bu
in \cite{bc} studied  (\ref{half}) with
$p=4$ and showed that if $\phi \in H^2 (\mathbb{R}^+)$ and
$h\in C^2(\mathbb{R}^+)$ with $\phi (0)=h(0)$, then the
problem admits a  unique global solution
$$
u\in
C^1(\mathbb{R}^+; L^2 (\mathbb{R}^+))\cap C(\mathbb{R}^+; H^2
(\mathbb{R}^+)).
$$
This result was  extended  to the case $p\geq 3$ by Bu in  \cite{bu1994}
for the defocusing case ($\lambda <0$).  In
\cite{strauss}, Strauss and Bu considered the problem
\begin{equation}
\left \{ \begin{array}{l} u_t -\Delta u+ \lambda |u|^{p-2} u =0,
\quad \ \ \,   x\in \Omega , \quad t\in \R \, ,\\ \\
u(x,0)= \phi (x), \quad x\in \Omega, \qquad
u(x,t)=q(x,t),  \quad x\in \partial \Omega \, ,\end{array}\right.
\label{0-1}
\end{equation}
for  the NLS equations
posed on a smooth  (bounded or unbounded) domain $\Omega \subset \mathbb{R}^n$.
Assuming that $\lambda <0$ and $p\geq 3$, they showed  that
for any $\phi \in H^1 (\Omega)$ and $q\in C^3 (\R^n \times
(-\infty, \infty))$ with compact support satisfying the natural
compatibility condition, the IBVP (\ref{0-1}) admits a global
solution
\[ u\in L^{\infty}_{loc} ((-\infty, \infty ); H^1 (\Omega)\cap L^{p}
(\Omega)).\]
Bu, Tsutaya and Zhang \cite{bu2005} extended the above result
to the case of $\lambda >0$ assuming $3\leq p\leq 2+\frac{n}{2}$ and
$n\geq 2$. In all this work, the third leg of Hadamard's conception, namely
continuous dependence of
solutions on the initial and boundary data, was not discussed.
For small $s \geq 0$, Holmer \cite{holmer}
obtained the following result for the half-line problem (\ref{half}).
\vspace{.15cm}

\begin{theorem} \label{holmer} (Holmer)
 Let
$ \frac12 < s< \frac 32 $ with $ 3\leq p< \infty $ or $
0\leq s< \frac12 $ with  $ 3\leq p<\frac{6-4s}{1-2s}$ be given. For any
$r>0$, there exists $T>0$ such  that if $\phi \in H^s
(\mathbb{R}^+)$ and $h\in H^{\frac{2s+1}{4}}_{loc}(\mathbb{R}^+) $
satisfy
\[
 \|\phi \|_{H^s(\R^+)} + \| h\|
_{H^{\frac{2s+1}{4}}(0,T)} \leq r\qquad \quad \Big ( h(0) =\phi (0) \quad\mbox{for}\quad s >\frac12 \Big )\, ,
\]
then the IBVP (\ref{half})
admits a solution $u\in C([0,T]; H^s(\R^+)) $ which depends continuously
upon the auxiliary data in the relevant function classes.  Moreover, for $
\frac12 < s< \frac 32 $, the solution $u$ is unique.
\end{theorem}

\vspace{.15cm}

This result is  very similar to that obtained here for the quarter-plane problem
\eqref{half}.  Our result, which is obtained by a different approach to
be described presently, improves Theorem \ref{holmer} in small ways (the issue
of uniqueness for  $s$ in the range $0 \leq s< \frac 12 $ is clarified and the
range of values of $s$ is extended).  The boundary integral method used in this paper
 and in our earlier work \cite{bsz-5} on the KdV equation, has other points to
recommend it, however.  First, one can read off from our representation of solutions a significant difference
between the IBVPs for the KdV equation and the nonlinear
Schr\"odinger equation.   For the KdV equation, the imposition of a boundary
condition at the left-hand end of $\mathbb{R}^+$ produces a
  strong dissipative smoothing mechanism, whereas  no such dissipative smoothing
appears from solving the same boundary-value problem for
the nonlinear Schr\"odinger equation. This distinction is not so clearly
seen using the earlier methods.  (More detail concerning
 this distinction   will be presented
elsewhere.)   Another point in favor of the boundary-integral method is
that it generalizes immediately to higher space dimensions.  This, also, is
a project for  future investigation.

\medskip

The discussion is now turned in a slightly more technical direction.
The first point we want to make is that at least for relatively small values of
$s$, the case where the boundary data is homogeneous ({\it i.e.} $h \equiv 0$
or $h_1 = h_2 \equiv 0$) can be reduced
to the situation described in Theorem A or Theorem B,
respectively.  (This is  no longer true for larger values of $s$, however.)
Thus, with essentially no effort, the following results obtain.

\begin{theorem}\label{homogeneous-1}
Assume that $h =0$ in (\ref{half}).
\begin{itemize}
\item[(i)] If  $\frac12 < s < \frac52$ with $3\leq p < \infty $ or $0\leq s< \frac12$ with $3\leq p <
\frac{6-4s}{1-2s}$,  the IBVP (\ref{half}) is locally well-posed in
$H^s(\mathbb{R}^+)$.
\item[(ii)] If $3\leq p< 6$, then the IBVP (\ref{half}) is (unconditionally) globally
well-posed in $H^s (\mathbb{R}^+)$ for any $s$ with $0 \leq s < \frac52$.
\end{itemize}
\end{theorem}
\begin{theorem}\label{homogeneous-2} Assume that $h_1 =0$  and $h_2=0$ in
(\ref{finite}).
\begin{itemize}
\item[(i)] If  $\frac12 < s < \frac52$ with $p\geq 3$ or $0\leq s <\frac12 $ with $3\leq
p<\frac{6-4s}{1-2s}$,  the IBVP  (\ref{finite}) is
 locally well-posed in $H^s
(0,L)$.
\item[(ii)] If $3\leq p< 6$, then the IBVP (\ref{finite}) is (unconditionally) globally
well-posed in $H^s (0, L)$ for any $s$ with $0 \leq s < \frac52$.
\end{itemize}
\end{theorem}

For Theorem \ref{homogeneous-1}, the result follows by choosing as initial data
 the odd extension $\tilde{\phi}$ of $\phi$, solving the equation on $\R$ with
$\tilde{\phi}$ as initial data and then restricting the resulting solution to the
half line.  For Theorem \ref{homogeneous-2}, extend $\phi$ to $[-L,L]$ by
taking the odd extension and then extend to all of $\R$ by $2L$--periodicity.
Solve the resulting periodic initial-value problem and then restrict to $[0,L]$.

\vspace{.15cm}

For the nonhomogeneous boundary-value problems that are the focus of
attention here, such simple methods do not appear to give results.
To deal with nonhomogeneous boundary data, a standard
approach is to homogenize the boundary data by a change of the dependent variables.
Define a new dependent variable by subtracting from the original dependent variable a known function that takes
on the given boundary values.  This new variable will satisfy a related equation, but
with zero boundary conditions.  While this works well in some cases, {\it e.g.} BBM-type
equations (see \cite{BCSZ} and the references therein), in the present context it requires that
the boundary data must have stronger regularity than should be needed according to
the heuristic analysis leading to the relation \eqref{optimal} between the 
function classes  of
the initial and the boundary data.  For
instance, this method, applied in a straightforward way for $p = 4$, say, ends up requiring for the
quarter-plane problem \eqref{half} that $h\in H^1 ([0,T])$ to obtain the
well-posedness of the IBVP (\ref{half}) in $L^2 (\R^+)$
rather than $h\in H^{\frac14} ([0,T])$ as advertised in Theorem \ref{th-half}, part $(ii)$.

The initial-boundary-value problem
\begin{equation}
\left \{ \begin{array}{l}u_t+u^ku_x +u_{xxx}=0, \ x\in \R^+, \ t\in (0,T)\, ,\\
\\ u(x,0)=\phi (x), \quad u(0,t) =h(t) , \end{array} \right.
\label{x-1}
\end{equation}
for the generalized Korteweg-de Vries (KdV) equation posed on a half line
$\R^+$\!,
is instructive.  Colliander and Kenig \cite{colliander} introduced a new method to
analyze this problem by solving the pure IVP
\begin{equation}
\left \{ \begin{array}{l}w_t+w^kw_x +w_{xxx}=\delta (x) f(t), \ x\in
\R, \ t\in (0,T)\, ,\\  \\ w(x,0)=\psi (x), \end{array} \right.
\label{x-2}
\end{equation}
of a forced,  generalized KdV equation with an appropriate forcing function $f(t)$.
Here, $\delta (x)$ denotes the Dirac mass at $x=0$ and $\psi (x) $ is an
extension of $\phi (x)$ from $\R^+$ to $\R$. It is demonstrated that
an appropriate  forcing function $f(t)$ can be  chosen  so that the
corresponding solution $w$ of (\ref{x-2}) satisfies
\[ w(0,t) = h(t), \quad {\rm for} \;\;0< t< T .\]
Consequently, the restriction of $w(x,t)$ to the half line $\R^+$ is
a  solution of the IBVP (\ref{x-1}).   The IVP (\ref{x-2}) is
solved using the contraction mapping principle in a
carefully constructed, Bourgain-type  space $X_{s,T}$.  The key step of this
approach is to study the associated linear problem,
\begin{equation}
\left \{ \begin{array}{l}v_t  +v_{xxx}=\delta (x) f(t), \ x\in \R, \
t\in (0,T)\, ,\\ \\ v(x,0)=\psi (x) \, ,\end{array} \right. \label{x-3}
\end{equation}
and show that there exists a real number $s'$ (depending only on
$s$) such that for any given $\psi \in H^s (\R)$,
\begin{itemize}
\item[(i)] if $f\equiv 0$, the solution $v$ of (\ref{x-3}) satisfies
\begin{equation}\label{x-6}
\sup _{x\in \R} \| v(x, \cdot )\|_{H^{s'} (0,T)} \leq C_{s,T} \|
\psi \|_{H^s(\R)}\, ,
 \end{equation}
 \item[(ii)]
if  $h\in H^{s'} (0,T)$, one can find a forcing function $f$ such
that the IVP (\ref{x-3}) admits a solution $v\in X_{s,T}$ and
\begin{equation} \| v\|_{X_{s,T}} \leq C_{s,T} \Big (\| \psi\|_{H^s (\R)}
+ \| h\|_{H^{s'}(0,T)} \Big ).\label{x-4} \end{equation}
\end{itemize}
 It turns out
that for the IBVP (\ref{x-1}),
\begin{equation}
 s'=\frac{s+1}{3} .\label{x-5}
 \end{equation}
 The estimate (\ref{x-6}) is, in fact, the sharp Kato  smoothing property possessed by
 the solutions of the linearized KdV equation. The Riemann-Liouville
 fractional integral is the main tool used to  establish the estimate (\ref{x-4}).

There is another approach to deal with the IBVP (\ref{x-1}) put forward by
the present authors in \cite{bsz-1}.  A major constituent of this latter
approach is the explicit solution formula
\begin{equation}\label{x-8}
q(x,t)= [U_b(t)h](x)+\overline{[U_b(t)h](x)}
\end{equation}
where
\[ [U_b(t)h](x)=\frac{1}{2\pi} \int ^{\infty}_1e^{i\mu ^3t-i\mu
t}e^{-\left (\frac{\sqrt{3\mu^2-4}+i\mu}{2}\right ) x}(3\mu^2-1)\int
^{\infty}_0e^{-i(\mu ^3 +\mu)\xi }h(\xi) \, d\xi \, d\mu,\]
of the linear,
nonhomogeneous boundary-value problem,
\begin{equation}
\left \{ \begin{array}{l} q_t +q_x +q_{xxx}=0, \quad x\in \R^+, \
t\in (0,T), \\ \\ q(x,0)=0, \qquad q(0,t) =h(t).
\end{array}\right.\label{x-7}
 \end{equation}
This explicit formula, which is obtained by formally taking the Laplace transform in time, solving the resulting third-order problem and taking the inverse Laplace transform,  enables one to establish directly
various  estimates needed for proving the well-posedness of the IBVP
(\ref{x-1}).  Moreover, it clearly demonstrates that  the solution $q(x,t)$ of  (\ref{x-7}) becomes infinitely smooth when $x>0$ and $t>0$. It has been further  shown  in \cite{bsz-5} that

\medskip
\noindent {\em
the solution $q(x, t) $   is the restriction  to $\R^+\times \R^+$ of
 a function $w(x, t)$ defined on $\R\times \R$ which is such that
 $$
  \left (\int ^{\infty }_{-\infty }  \int ^{\infty }_{-\infty }
 (1+|\xi |) ^{2s}
  (1+ |\tau - \xi ^3 |)^{2b}  | \hat{w} (\xi , \tau ) |^2 d\xi
 d\tau\right  )^{1/2}\leq C \|h  \| _{H^{\frac{3b+s -1/2}{3}} (\R^+)}
 $$
where $b$ is any value in $[0,  \frac12 -\frac{s}{3}) $ if $-\frac32 \leq  s< \frac32 ,   $  $b$ is any value in $[0,
\frac56-\frac{s}{3} ]$ if  $-\frac12 < s< 1$ and $C$ is a constant depending only on $s$ and $b$. }

\medskip
It then follows that the IBVP  (\ref{x-7}) possesses the following strong dissipative smoothing property:
\[ h\in H^{(s+1)/3}_{loc} (\R^+) \implies  q\in L^2 (0,T; H^
{s+\frac{3}{2}}(\R^+)).\]

In \cite{holmer}, Holmer applied the Colliander--Kenig approach to
study the IBVP (\ref{half}) and obtained  the results
described in Theorem \ref{holmer}. However, as we will
show in this paper, this approach may fail for the IBVP
(\ref{finite}). More precisely, we show that   for the solution $u$ of the IBVP
\begin{equation}
\left \{ \begin{array}{l} iu_t + u_{xx}=0, \quad
x\in (0,L), \ t\in (0,T), \\ \\ u(x,0)=0, \quad u(0,t)=h_1 (t), \quad u(L,t) =h_2
(t),\end{array} \right.  \label{x-10}
\end{equation}
for the linear Schr\"odinger equation posed on $(0,L)$,
the estimate
\begin{equation} \label{e-1}
\| u\|_{L^2 ((0,L)\times (0,T))} \leq C\left (\| h_1\|_{H^{\alpha} (0,T)}
 +\|h_2 \|_{H^{\alpha}(0,T)} \right )
 \end{equation}
holds if $\alpha \geq \frac12$,  but fails if $\alpha < \frac12$.
(Indeed, Example A2 in the Appendix shows the optimality of the assumption $h_1 , h_2 \in H^{1/2}(0 , T)$ for this estimate to hold).
 By contrast,  for solutions of the IBVP
\begin{equation}
\left \{ \begin{array}{l} iv_t + v_{xx}=0, \quad
x\in \R^+, \ t\in (0,T), \\ \\ v(x,0)=0, \quad  v(0,t)=h (t),
\end{array} \right. \label{x-11}
\end{equation}
for the linear Schr\"odinger
 equation posed on $\R^+$\!,
it is indeed the case that
 \begin{equation} \label{e-2}
\| v\|_{L^2( \R^+\times (0,T))} \leq C\| h\|_{H^{\frac14} (0,T)}.
\end{equation}
And,   solutions of the pure IVP
\begin{equation}
iw_t +w_{xx}=0, \quad w(x,0)=\psi (x), \ x\in \R, \quad t\in (0,T)
\label{x-12}
\end{equation}
for the linear Schr\"odinger
equation posed on $\R$,
comply with the inequality
\begin{equation} \label{e-3}
\sup _{x\in \R} \| w(x, \cdot)\|_{H^{\frac14}(0,T)} \leq C \| \psi
\| _{L^2 (\R)} .
\end{equation}
Thus, while it is possible to solve the nonhomogeneous IBVP
(\ref{half}) by solving a forced IVP of the form
\[
iu_t +\lambda |u|^{p-2} u +u_{xx}=\delta (x) f(t), \ u(x,0) =\psi
(x), \quad x\in \R, \ t\in (0,T) ,
\]
 with
an appropriate forcing function $f(t)$, it may not be feasible to apply the same
approach to the two-point IBVP
(\ref{finite}).

In this paper, the  approach
developed earlier in \cite{bsz-1} for studying  nonhomogeneous boundary-value
 problems of the KdV equation will be used  to establish  local well-posedness results for
(\ref{half}) and (\ref{finite}).  Analogous to  the solution formula \eqref{x-8} in the KdV case,
 the  nonhomogeneous, linear IBVP (\ref{x-11}) has the
 explicit solution
\begin{equation}\label{e-4}
v(x,t)= \frac{1}{\pi} \int ^{\infty}_0 e^{-i\beta^2 t}e^{i\beta x}
\beta \int ^{\infty}_0 e^{i\beta ^2 \tau} h(\tau) d\tau d\beta +
\frac{1}{\pi}\int ^{\infty}_{0} e^{i\beta ^2t}e^{-\beta x}
 \beta \int ^{\infty}_0 e^{-i\beta ^2 \tau} h(\tau) d\tau d\beta.
\end{equation}
Similarly, the solution formula for the nonhomogeneous, linear IBVP (\ref{x-10}) is
 \begin{equation}
\label{e-6}
 u(x,t)=\sum ^{\infty }_{n=1} 2 in\pi e^{-i(n\pi )^2 t} \int
^t_0 e^{i(n \pi)^2 \tau} \Big(h_1  (\tau)-(-1)^nh_2 (\tau )\Big )d\tau \sin
n\pi x.
\end{equation}
As in the case of the KdV equation, these formulas are derived by taking the Laplace transform of $u$
in the temporal variable, solving the resulting, second-order, ordinary
differential equation and taking the inverse Laplace transform of the result.
The inequalities needed to advance the local well-posedness theory
 obtain directly from these explicit
solution formulas.  Moreover, from these  formulas, one ascertains
 that, unlike the KdV equation, the imposition of boundary
conditions brings no smoothing effect.  For example,  consider 
the IBVP (\ref{x-11}).  The second term
\[ B(x,t)= \frac{1}{\pi}\int ^{\infty}_{0} e^{i\beta ^2t}e^{-\beta x}
 \beta \int ^{\infty}_0 e^{-i\beta ^2 \tau} h(\tau) d\tau d\beta\]
 on the right-hand side of  the solution formula (\ref{e-4})
becomes infinitely smooth as soon as $x>0$ and $t>0$.  
On the other hand, 
  the first term  on the right-hand side of   (\ref{e-4}) 
can be written as
\[ A (x,t)=  \frac{1}{\pi} \int ^{\infty}_{-\infty} e^{-i\beta^2 t}e^{i\beta x} \hat{\psi} (\beta)
 d\beta,\]
 where $\psi $ is the function whose  Fourier transform is 
 \[ \hat{\psi} (\beta)= \left \{ \begin{array}{ll} \beta \int ^{\infty}_0 e^{i\beta ^2 \tau} h(\tau) d\tau & \ if \ \beta >0, \\  \\ 0 & \ if \ \beta < 0. \end{array} \right.\]
Thus, $A(x,t)$ solves the pure initial-value problem for the 
linear Schr\"odinger equation, 
 posed on the whole line $\R$, with the initial value $\psi (x)$.
 It follows that  $\psi \in H^s(\R) $ if and only if $h\in H^{\frac{2s+1}{4}}_0 (\R^+)$. Consequently, in contrast to the KdV equation, there is no boundary smoothing 
for the Schr\"odinger equation.

This section is concluded with remarks on higher-order
regularity and global well-posedness.
  The theory outlined above, and which is developed in detail
in the remainder of the essay, has upper limits on the regularity of the
auxiliary data.  As we will see momentarily, these restrictions are
necessary.   They can be relaxed only by asking for additional properties
of the auxiliary data.

When equations like the Schr\"odinger equation are derived to describe
physical phenomena, they often come as a simplification of a more
complete model.   Justifying the simpler model as an approximation of
a more elaborate model typically requires smoothness of the solutions
of both the full and the approximate models (see \cite{BCL2005, Crai1985, SW2000} for justification of
the KdV equation as an approximation of the full water-wave problem,
 for instance).  Without smoothness, the comparisons
are not in fact valid.  Thus, it is not only of academic interest to understand
higher regularity solutions.

An example will illustrate the problem that arises when smoother solutions
are in question.  Take the classical case $p = 4$ so that the nonlinearity is cubic and
smooth.   Suppose that the quarter-plane problem \eqref{half} is
locally well posed in $H^3(\R^+)$, say.
Then, there is a $T>0$ and a solution
$u \in C(0,T;H^3(\R^+))$.  Because $u_{xx} \in C(0,T;H^1(\R^+))$ and
$u$ satisfies the equation, it must be the case that $u \in C^1(0,T;H^1(\R^+))$.    It
follows that each term in the evolution equation is a continuous function of
both space and time in $\R^+ \times [0,T]$.       Evaluating the equation
at the point $(0,0)$ and using the initial and boundary conditions then
yields
\begin{equation} \label{higher}
ih'(0) \,+\, \phi ''(0) \,+\, \lambda |\phi(0)|^{2} \phi(0) \, = \, 0.
\end{equation}
Thus, the auxiliary data necessarily satisfies a higher-order compatibility
condition in addition to the lower-order condition \eqref{compat} that
has been assumed throughout the discussion.   It is straightforward to
calculate yet  higher-order conditions on the auxiliary data that must
obtain for well-posedness to hold in smaller Sobolev spaces.   This issue
also arises for the KdV equation posed on the half-line or on a
bounded interval.   In that case, higher-order regularity theory has
been developed in the presence of higher-order compatibility conditions
(see \cite{bsz-6, BW1983}).

When the nonlinearity is smooth, {\it e.g.} when $p = 4,6,8,\cdots$, local
well posedness in the presence of higher regularity and the associated
compatibility conditions can be established by the methods put
forward here.  However, we eschew this task in the present script.

Finally, we come to the issue of global well-posedness.  As is standard in the
theory of evolution equations, local well-posedness
 coupled with suitable {\it a priori} bounds on solutions is the path to
global well posedness.   For the pure initial-value problem \eqref{ivp},
  the  bounds  provided by the  conserved quantities
\begin{equation}\label{conserved}
I(t):=\int ^{\infty}_{-\infty}|u(x,t)|^2 dx \quad {\rm and} \quad
II(t):= \int ^{\infty}_{-\infty}\left ( |u_x(x,t)|^2
-\frac{2\lambda}{p} |u(x,t)|^p \right )d x
\end{equation}
suffice for the global results mentioned earlier.
However, corresponding to the quarter-plane problem (\ref{half}), one has (cf.  \cite{bc})
\begin{equation}\label{c-3-1}
I'(t):= \frac{d}{dt} \int ^{\infty}_{0}|u(x,t)|^2 dx =-2\im
\Big (u_x(0,t)\overline{h}(t)\Big )
\end{equation}
and
\begin{equation}\label{c-4}
II'(t):=\frac{d}{dt} \int ^{\infty}_{0}\left ( |u_x(x,t)|^2
-\frac{2\lambda}{p} |u(x,t)|^p \right )d x =-2\re \Big (u_x
(0,t)\overline{h}\, '(t)\Big )
\end{equation}
while the two-point IBVP (\ref{finite}) has
\begin{equation}\label{c-5}
I'(t):= \frac{d}{dt} \int ^{L}_{0}|u(x,t)|^2 dx =2\im \Big (u_x
(L,t)\overline{h}_2 (t)- u_x(0,t)\overline{h}_1(t)\Big )
\end{equation}
and
\begin{equation}\label{c-6}
II'(t):= \frac{d}{dt} \int ^{L}_{0}\left ( |u_x(x,t)|^2
-\frac{2\lambda}{p} |u(x,t)|^p \right )d x =2\re \Big (u_x
(L,t)\overline{ h}\, '_2(t)-u_x (0,t)\overline{h}\, '_1(t)\Big ),
\end{equation}
for all $t \in \R$ for which the solutions exist.
In case  the  boundary conditions are homogeneous, {\it viz.}  $h\equiv 0$  or
$h_1=h_2 \equiv 0$,  both $I(t)$ and $II(t)$ are formally conserved just as
in the case of the pure initial-value problem \eqref{ivp}.
At least for small values of the Sobolev index $s$,   global
well-posedness results for the homogeneous IBVP's (\ref{half}) and (\ref{finite})
then follow readily.   For the nonhomogeneous cases, both $I(t)$
and $II(t)$ are no longer conserved and the task of obtaining
global \emph{a priori} estimates becomes interesting (see Section 5).

We turn now to the body of the paper where
 detailed analysis is given leading to the conclusions
 recounted in the Introduction.  The explicit solution formulas
\eqref{e-4} and \eqref{e-6} will play a central role in our development.

\section{The Schr\"odinger equation posed on the half line $\R ^+$}
\setcounter{equation}{0}

Considered first is the IBVP (\ref{half})
\begin{align}\label{n-1-homo}
\left \{ \begin{array}{l} iu_t +u_{xx} + \lambda u |u |^{p-2} =0,
 \quad x\in \mathbb{R}^+, \ t\in
\mathbb{R},
\\ \\   u(x,0) =\phi (x) ,\quad u(0,t) = 0\, .
\end{array}\right.
\end{align}
 with a homogeneous boundary
condition.
It transpires that  this can be reduced to the pure
IVP
\begin{align}\label{n-1-R}
\left \{ \begin{array}{l} iw_t +w_{xx} + \lambda w |w |^{p-2} =0,
\qquad x\in \mathbb{R} , \ t\in \mathbb{R},
\\ \\ w(x,0) =\psi (x)
\end{array}\right.
\end{align}
of the NLS equation posed on the whole line $\mathbb{R}$.
Indeed, observe that if $w=w(x,t)$ is a solution of
(\ref{n-1-R}) which is an odd function with respect to $x$, then its restriction
\[ u(x,t):=w(x,t), \quad x\in \R^+,\]
to the half-line is a solution of  (\ref{n-1-homo}) with
$\phi (x) =\psi (x), x\in \R^+ .$  On the other hand,  the
IVP (\ref{n-1-R}) possesses the following invariance property.
\begin{lemma}
If $\psi $ is an odd and smooth function, then for any $t \in \R$, the corresponding solution $w$
of (\ref{n-1-R}) is odd with respect to $x$.
\end{lemma}
\noindent {\bf Proof:} Consider first the associated linear problem
\begin{align}\label{n-1-R-1}
\left \{ \begin{array}{l} iw_t +w_{xx}  =f, \qquad x\in \mathbb{R} ,
\ t\in \mathbb{R},
\\ \\ w(x,0) =\psi (x)
.\end{array}\right.
\end{align}
Using the Fourier transform, its solution $w(x,t)$ is
\[
w(x,t)=\int _{\R} e^{i\xi^2 t} e^{i\xi x}\hat{\psi} (\xi) d\xi +
\int ^t_0 \int _{\R}e^{i\xi ^2 (t-\tau)} e^{i\xi x} \hat{f} (\xi,
\tau) dx\xi d\tau \, .
\]
It then follows
directly that the solution $w(x,t)$ of (\ref{n-1-R-1}) is odd with
respect to $x$ if $\psi $ and $f$ are odd in $x$. For  the IBVP (\ref{n-1-R}), suppose
$\psi$ is odd and consider the map
$\Gamma: v\mapsto w$,
where $v=v(x,t)$ is an odd function in $x$  and $ w$ is the
solution of
\begin{align}
\left \{ \begin{array}{l} iw_t +w_{xx}   =-\lambda |v|^{p-2} v, \qquad
x\in \mathbb{R} , \ t\in \mathbb{R},
\\ \\ w(x,0) =\psi (x)
.\end{array}\right.
\end{align}
It follows from the previous remark about \eqref{n-1-R-1} that $\Gamma (v)$ is odd in $x$  if $v$ is odd in $x$.
The classical contraction mapping principle provides the
solution $w$ of the nonlinear IVP (\ref{n-1-R}).  This solution is necessarily odd
 as a function of $x$ if
its initial value $\psi $ is odd, as one determines by iterating $\Gamma$ starting at $v = 0$.
 $\Box$

\vspace{.07cm}

Thus, the following well-posedness result for the IBVP
(\ref{n-1-homo}) follows from the well-posedness of the IVP
(\ref{n-1-R}).
\begin{theorem}
For any $s$ satisfying either $ \frac12 < s < \frac52$ for $3 \leq p < \infty$, 
$\lfloor s \rfloor < p-2 $ if $s\not = 1, 2$, or $0\leq s< \frac12$ for $3\leq p <
\frac{6-4s}{1- 2s}$, the IBVP (\ref{n-1-homo})
is locally well-posed in $H^s(\R ^+)$ (for $\frac12 < s < \frac52$, it is required that
$\phi (0) = 0$).
\end{theorem}

Now, (\ref{half}) is considered with nonhomogeneous boundary data. The analysis of this problem is
carried out in several subsections.

\subsection{Solution formulas for linear problems}

 Consideration is first given to  the linear, nonhomogeneous, boundary-value problem
\begin{equation}\label{2.1}
\left \{ \begin{array}{l} iu_t +u_{xx} =0, \qquad x\in \R^+, \ t\in
\R^+, \\ \\ u(x,0) =0, \qquad u(0,t) = h(t) .\end{array}\right.
\end{equation}
By taking the  Laplace transform with respect to \emph{t} of both
sides of (\ref{2.1}), the IBVP is converted to a one-parameter
family of second-order boundary-value problems, {\it viz.}
\begin{equation}
\left \{ \begin{array}{l} i\lambda
\tilde{u}(x,\lambda)  +\tilde{u}_{xx}(x,\lambda )=0,  \\ \\
\tilde{u}(0,\lambda)=\tilde{h}(\lambda ), \qquad \tilde{u}(+\infty ,
\lambda)=0, \end{array}\right. \label{2.2}
\end{equation}
where
$\tilde{u}=\tilde{u}(x, \lambda)$ is the Laplace transform of
$u=u(x,t)$ with respect to
 $t$ and $\re \, \lambda >0$ is the dual variable. The solution of
 (\ref{2.2}) is given by
 \[ \tilde{u}(x, \lambda )= e^{r (\lambda)x} \tilde{h}(\lambda)\]
 where $r (\lambda)$ is the solution of the quadratic equation
 \[ i\lambda +r ^2 =0 \]
 for which $\mbox{Re}\, r <0$.   In consequence, the solution of (\ref{2.2}) is given
 formally by
 \[ u(x,t)=\frac{1}{2\pi i}\int ^{+\infty i +\gamma}_{-\infty i
 +\gamma}e^{\lambda t}e^{r(\lambda) x}\tilde{h}(\lambda )d\lambda \]
 for $x, \ t>0$, where $\gamma >0$ is fixed. Letting $\gamma \to 0$,
 one arrives at
\begin{eqnarray*}
 u(x,t) &=&  \frac{1}{2\pi } \int ^{\infty}_{-\infty} e^{i\beta t}
 e^{r (i\beta )x} \tilde{h} (i\beta ) d\beta  \qquad \Big( -\beta
 +r ^2 =0, \qquad \re\,r \leq 0 \Big) \\
 &=& \frac{1}{2\pi} \int^0_{-\infty } e^{i\beta t}
 e^{i\sqrt{-\beta}x} \tilde{h} (i\beta ) d\beta +  \frac{1}{2\pi} \int^{\infty }_0 e^{i\beta t}
 e^{-\sqrt{\beta}x} \tilde{h} (i\beta ) d\beta \\ \\
 &=&
 \frac{1}{2\pi} \int^{\infty }_0 e^{-i\beta t+
  i\sqrt{\beta}x} \tilde{h} (-i\beta ) d\beta +   \frac{1}{2\pi} \int ^{\infty }_0 e^{i\beta t-
  \sqrt{\beta}x}\tilde{h} (i\beta ) d\beta \\ \\
  &=&
 \frac{1}{\pi} \int^{\infty }_0 e^{-i\beta ^2 t+
  i \beta x} \beta \tilde{h} (-i\beta ^2) d\beta   + \frac{1}{\pi} \int ^{\infty }_0 e^{i\beta ^2 t-
  \beta x} \beta \tilde{h} (i\beta ^2) d\beta \\ \\
  &=& I(x,t)+II(x,t).
  \end{eqnarray*}

For $I(x,t)$, define \begin{equation} \nu _1 (\beta)= \left \{
\begin{array}{ll}\frac{1}{\pi}\beta \tilde{h}(-i\beta ^2)& \ \mbox{
for $\beta \geq 0$}, \\ \\ 0 & \ \mbox{ for $\beta < 0$}
\end{array} \right.\label{2.3-1} \end{equation}  and
\begin{equation}
\phi _h =\phi _h (x) \label{2.3-2}
\end{equation}
to be the inverse Fourier transform of $\nu _1 $, so that the Fourier
transform $\hat{\phi_h}$ of $\phi_h$ is
 \[ \hat{\phi}_h
(\beta) =\nu_1 (\beta) , \qquad \beta \in \R.\]
Then,
$I(x,t)$ can be rewritten as
\[ I(x,t) = \int ^{\infty}_{-\infty} e^{-i\beta ^2 t+i\beta x}
\hat{\phi}_h (\beta ) d\beta \, ,\] which is exactly the
solution formula of the Cauchy problem for the linear Schr\"odinger
equation on $\R$.
Similarly, for $II(x,t)$, define \begin{equation} \nu _2 (\beta)=
\left \{
\begin{array}{ll}\frac{1}{\pi}\beta \tilde{h}(i\beta ^2)& \ \mbox{
for $\beta \geq 0$}, \\ \\ 0 & \ \mbox{ for $\beta < 0$}
\end{array} \right.\label{2.4-1} \end{equation} and
\begin{equation}
\psi _h =\psi _h (x) \label{2.4-2} \end{equation} to be the inverse
Fourier transform of $\nu _2 $, i.e., \[ \hat{\psi}_h (\beta)
=\nu_2(\beta) , \qquad \beta \in \R.\]
Thus, $II(x,t)$ can be written as
\[ II(x,t)=   \int ^{\infty}_{-\infty} e^{i\beta ^2 t-\beta
x}\hat{\psi}_h (\beta)d\beta \]
for $x>0$.

\begin{proposition} The solution of (\ref{2.1}) may be written as
\[ u(x,t)= [W_b (t)h](x) := [W_{b,1} (t)h](x) + [W_{b,2}(t)h](x) \]
where for $x, \ t>0$,
\[ \big [W_{b,1} (t)h \big](x)= \int ^{\infty}_{-\infty} e^{-i\beta ^2 t+i\beta x}
\hat{\phi}_h (\beta ) d\beta , \]
\[ \big [W_{b,2} (t)h\big ](x)= \int ^{\infty}_{-\infty} e^{i\beta ^2 t-\beta
x}\hat{\psi}_h (\beta)d\beta \] and $\phi _h, \psi _h$ are defined
by (\ref{2.3-1})-(\ref{2.3-2}) and (\ref{2.4-1})-(\ref{2.4-2}),
respectively.
\end{proposition}
\begin{remark}\label{rem2.4}

\quad
\begin{itemize}
\item[(i)]
It follows from their definitions that   for any $s\in \mathbb{R}, \phi _h  $ and
$\psi _h $ belong to the space $H^s (\mathbb{R})$ if and only if
$h\in H^{\frac{2s+1}{4}}_0 (\mathbb{R}^+)$.
\item[(ii)] The function $v(x,t)= [W_{b,1} (t)h](x)$ is, in fact, defined for $x,
\ t\in \mathbb{R}$ and solves the  IVP
\[ iv_t + v_{xx} =0, \quad v(x,0)= \phi _h (x), \qquad x,\ t\in \mathbb{R}\]
 for
the linear Schr\"odinger equation
posed on $\mathbb{R}$.
As for $[W_{b,2} (t)h](x)$, it is defined only for $x>0$. However, it may be
 extended for $x< 0$  by setting
\begin{equation}  \big [W_{b,2} (t)h\big ](x)= \int ^{\infty}_{-\infty}
e^{i\beta ^2 t-\beta |x|}\hat{\psi}_h (\beta)d\beta .\label{2.5}
\end{equation}
Note that this extension is not necessarily differentiable at $x =0$. Therefore,
this small trick is not applicable when $ s > 3/2$.
\end{itemize}
\end{remark}

Next, consider the same linear equation
\begin{equation}
\label{initial} \left \{ \begin{array}{l} iu_t +u_{xx} =0, \qquad
x\in \mathbb{R}^+, \ t\in \mathbb{R}^+, \\ \\ u(x,0) =\phi(x) ,
\qquad u(0,t) = 0\end{array}\right.
\end{equation}
with zero boundary condition,  but non-trivial initial data.
By semigroup theory, its solution $u$ may be obtained in the
form
$$
u(t) = W_0(t) \phi
$$
where the spatial variable is suppressed and  $W_0(t)$ is the
$C_0$-group in $L^2(\mathbb{R}^+)$ generated by the
operator $A$ defined by \[ A v=iv''\] with domain \[{\cal D}(A)=
\{ v\in H^2(R^+)\ |  \ v(0)=0 \} .\] By Duhamel's principle, one may
use the semi-group $W_0(t)$ to formally write the solution of the
forced linear problem
\begin{equation}
\label{2.9-1} \left \{ \begin{array}{l} iv_t +v_{xx} =f, \qquad x\in
\mathbb{R}^+, \ t\in \mathbb{R}^+,
\\ \\ v(x,0) =0 , \qquad v(0,t) =  0\end{array}\right.
\end{equation}
in the form
\[ v(x,t)=-i \int ^t_0 W_0 (t-\tau ) f(\cdot, \tau ) d\tau .\]



Let a function $\phi $ be defined on the half line $\mathbb{R}^+$
and let $\phi ^*$ be an extension to the whole line $\R$. The
mapping $\phi \mapsto \phi ^*$ can be organized so that it defines a
bounded linear operator $B$ from $H^s (\mathbb{R}^+)$ to $H^s
(\mathbb{R})$.
Henceforth $\phi ^*= B\phi $ will refer to such an
extension operator applied to $\phi \in H^s (\mathbb{R}^+)$. Assume
that $v=v(x,t) = W_\R (t) \phi^*$ is the solution of
\[ iv_t + v_{xx} =0, \quad v(x,0)=\phi^*(x),\]
for $ x, \ t\in \mathbb{R}$. If $g(t) =v(0,t),$ then $v_g =v_g
(x,t)= W_b (t)g$ is the corresponding solution of the nonhomogeneous
boundary-value problem (\ref{2.1}) with boundary condition
$h(t)=g(t)$, for $t\geq 0$.
Similarly, the function
\[ w\equiv w(x,t)=\int ^t_0
W_{\R} (t-\tau ) f^* (\tau ) d\tau
\]
with $f^* (x,t)= Bf(x,t)$
solves
\[ iw_t +w_{xx}=f^*(x,t), \qquad w(x,0)=0
\, ,
\]
for $x, \ t \in
\R$. If $p(t)=w(0,t)$, then $w_p\equiv w_p (x,t)=W_b (t)p = W_{bdr} (t) p $
is the corresponding solution of the non-homogeneous boundary-value
problem (\ref{2.1}) with boundary condition $h(t)=p(t)$, for $t\geq
0$.  The following integral representation thus obtains for
solutions of the fully non-homogeneous linear initial-boundary-value
problem
\begin{equation}
\label{iforcing} \left \{ \begin{array}{l} iu_t +u_{xx} =f, \quad
x, \ t\in \R^+ , \\ \\
u(x,0) =\phi (x), \qquad u(0,t)
=h(t)\, .\end{array} \right.
\end{equation}

\begin{proposition}  The solution $u$ of (\ref{iforcing}) can be
expressed as
\begin{equation} \label{integral}
 u(t) = W_{\R} (t)\phi ^* +\int ^t_0 W_{\R} (t-\tau ) f^* (\tau)
 d\tau + W_{bdr}(t)\big (h(t)-g(t)-p(t)\big )
 \end{equation}
 where
 \[ \phi^ * (x,t) =\big [B\phi\big ] (x,t), \qquad f^* (x,t)= \big [Bf\big ] (x,t) \]
 and
 \[ g(t)= \left. W_{\R} (t) \phi ^* \right |_{x=0}, \qquad p (t)=\left. \int ^t_0 W_{\R}
 (t-\tau ) f^*(\tau ) d\tau \right |_{x=0} .\]
\end{proposition}

\subsection{Linear estimates}

\quad
As before, for any $q\geq 2 $ and $r\geq 2$, the pair $(q,r)$ is called {\it admissible} if
\begin{equation}  \label{admissible}  \frac{1}{q}+\frac{1}{2r}=\frac14\, .\end{equation}
For any  $q$ with $ 1 \leq q \leq \infty$, $q'$ will denote the 
Lebesgue index conjugate to $q$, which is to say,
$\frac{1}{q} + \frac{1}{q'}=1$.

 \smallskip
 The following estimates for solutions of the
linear Schr\"odinger equation posed on the whole line $\R$ are well known
in the subject and will find use here.
\begin{proposition} \label{prop2.9}
Let $s\in \mathbb{R}$ and $T>0$ be given.
For any $\phi \in H^s (\mathbb{R})$, let $u=W_\mathbb{R} (t) \phi .$
Then, there exists a constant $C$ depending only on $s$ such that
\[ \sup _{t\in (0,T)} \|u(\cdot, t)\| _{H^s (\mathbb{R})} \leq C\|
\phi \| _{H^s (\R)}, \]
\[
\sup _{x\in \mathbb{R}} \|  u(x, \cdot
)\|_{H^{\frac{2s+1}{4}}(0,T)}\leq C \| \phi \| _{H^s (\mathbb{R})} \quad \]
and \[ \|u\| _{L^q (0,T;W^{s,r}(\R))}\leq C\|\phi
\|_{H^s(\mathbb{R})}\, ,\] for any given admissible pair $(q,r)$.
\end{proposition}
This proposition is same as Lemma 4.1 in \cite{holmer}.
\begin{proposition} \label{prop2.10}
Let $(q,r)$ be admissible and $T>0$ be given. Suppose  $f\in L^{q'}(0,T;
W^{s, r'} (\R))$  and define
\[ u=\int ^t_0 W_{\R} (t-\tau ) f(\tau) d\tau .\]
\begin{itemize}
\item[(i)] For any $s\in \R$, there exists a constant $C>0$ depending
only on $s$ such that
\begin{equation}\label{f-1}  \| u\| _{C([0,T]; H^s (\R))} + \| u\| _{L^q(0,T; W^{s,r}
(\R))} \leq C \| f\|_{L^{q'} (0,T; W^{s, r'} (\R))} . \end{equation}
\item[(ii)] For any $s\in (-\frac32, \frac12)$, there exists a constant $C>0$ depending
only on $s$ such that
\begin{equation}\label{f-2}
\sup _{x\in \R} \| u(x, \cdot )\| _{H^{\frac{2s+1}{4}} (0,T)} \leq
C(1+T)^{\frac12}\| f\| _{L^{q'} (0,T; W^{s, r'} (\R))}  .
\end{equation}
\item[(iii)] For any $s\in \R$, there exists a constant $C>0$  such that
\begin{equation}\label{f-3}
\sup _{x\in \R} \|  u(x, \cdot )\| _{H^{\frac{2s+1}{4}} (0,T)} \leq
C\| f\| _{L^{1} (0,T; H^{s} (\R))} .
\end{equation}
\end{itemize}
\end{proposition}
\noindent {\bf Proof:} The proof of (\ref{f-1}) can be found in
\cite{cazenave}.  A proof of (\ref{f-2}) is provided in \cite{holmer}.  For (\ref{f-3}),
note that \begin{eqnarray*}
 \sup _{x\in \R} \|  u(x, \cdot )\| _{H^{\frac{2s+1}{4}} (0,T)}&\leq& \int ^T_0\sup _{x\in \R} \| W_{\R}
 (t-\tau ) f(\tau ) \|_{H^{\frac{2s+1}{4}}_t (0,T)} d\tau \\
 &\leq & C\int ^T_0 \| f(\cdot, \tau)\|_{H^s (\R)} d \tau \, ,\end{eqnarray*}
thereby completing the analysis. $\Box$

\medskip
Next, consider  the  boundary integral operator $W_{bdr} (t)$.

\begin{proposition}\label{prop2.11} Let $0\leq s\leq 1 $ and $T>0$ be given and
suppose  $(q,r)$ is an admissible pair. There exists a constant $C>0$ such that
\begin{equation}
\| W_{bdr} (\cdot )h\| _{L^q (0,T; W^{s,r}(\R ))} \leq C \| h\|
_{H^{\frac{2s+1}{4}} (\mathbb{R}^+)} ,\label{2.7}
\end{equation}
\begin{equation}
\sup _{0<t<T} \| W_{bdr} (\cdot )h\| _{H^s(\R)} \leq C \| h\|
_{H^{\frac{2s+1}{4}} (\mathbb{R}^+)} \label{2.8}
\end{equation}
and \begin{equation} \sup _{x\in \R} \| W_{bdr} (\cdot )h\|
_{H^{\frac{2s+1}{4}}_t(0,T)} \leq C \| h\| _{H^{\frac{2s+1}{4}} (\mathbb{R}^+)}\, ,
\label{2.8-1}
\end{equation}
for  any $h\in H^{\frac{2s+1}{4}}_0 (\R^+)$.
\end{proposition}
{\bf Proof:}  It is sufficient  to prove that
\begin{equation}
\| W_{b,2} (\cdot )h\| _{L^q (0,T; W^{s,r}(\R ))} \leq C \| h\|
_{H^{\frac{2s+1}{4}} (\mathbb{R}^+)} \label{2.9}
\end{equation}  since
\begin{equation}
\| W_{b,1} (\cdot )h\| _{L^q (0,T; W^{s,r}(\R ))} \leq C \| h\|
_{H^{\frac{2s+1}{4}} (\mathbb{R}^+)}  \label{2.10}
\end{equation}
can be obtained from the result for the whole real line  given in
\cite{cazenave} and Remark \ref{rem2.4}. To show (\ref{2.9}), note
that
\begin{eqnarray*}
\big [W_{b,2} (t)h\big ](x)&=& \frac{1}{\pi} \int ^{\infty}_0 e^{i\beta ^2
t-\beta
|x|} \beta\hat{h}(i\beta ^2) d\beta
  \\ \\
&=& \frac{1}{\pi} \int ^{\infty }_0 e^{i\beta ^2t-\beta |x|} \int
^{\infty }_{-\infty } e^{-iy\beta}\psi_h (y)dy d\beta
= \frac{1}{\pi}\int ^{\infty }_{-\infty } \psi_h  (y)    \int
^{\infty }_0 e^{i\beta ^2t-\beta |x|-iy\beta } d\beta dy \\ \\
&:=& \int ^{\infty }_{-\infty } \psi_h  (y)K_t(x,y)dy
\end{eqnarray*}
where \[ K_t (x,y)= \frac{1}{\pi} \int ^{\infty }_0 e^{i\beta ^2
t-\beta |x|-iy\beta} d\beta.\]

\medskip
{\bf Claim:} {\em There exists a constant $C>0$ independent of $t, x, y$ such that for any
$t\ne 0$, $x, \ y\in \R \, ,$
\begin{equation} \left |K_t (x,y)\right |
\leq \frac{C}{\sqrt{|t|}} \, .\label{claim}
\end{equation}}

\medskip
\noindent {\bf Proof of the Claim:}   Note that although the Van Der Corput
lemma (Corollary 1.1 in \cite{LP2009}) can be used to shorten the proof of the claim, 
we present a  self-contained
argument in favor of \eqref{claim} here. Our approach is the following:
\begin{eqnarray*}
K_t (x,y)&=& \frac{1}{\pi} \int ^{\infty }_0 e^{i\beta ^2 t-\beta
|x|-iy\beta} d\beta = \frac{1}{\pi \sqrt{t}} \int ^{\infty }_0
e^{i\beta ^2 -\beta |x| t^{-\frac12} - iy \beta t^{-\frac12}} d\beta
\\ \\
&=& \frac{1}{\pi \sqrt{t}} \int ^{\infty }_0 e^{i \left ( \beta
-\frac{y}{2\sqrt{t}}\right )^2 -\frac{\beta | x|}{\sqrt{t}}} d\beta
e^{-i\frac{y^2}{4t}}= \frac{1}{\pi \sqrt{t}} \int ^{\infty
}_{-\frac{y}{2 \sqrt{t}}} e^{i\beta ^2 -\frac{|x|}{\sqrt{t}}\left
(\beta
+\frac{y}{2\sqrt{t}}\right )}d\beta  e^{-i\frac{y^2}{4t}} \\ \\
&=& \frac{1}{\pi \sqrt{t}} e^{-\frac{|x|y}{2t} -i\frac{y^2}{4t}} \int ^{\infty
}_{-\frac{y}{2 \sqrt{t}} }e^{i\beta ^2 - \frac{|x|}{\sqrt{t}}\beta }
d\beta \, .
\end{eqnarray*}
If $y\leq 0$,
\begin{eqnarray*}
\left | K_t (x,y)\right | &=& \frac{1}{\pi \sqrt{t}}
e^{-\frac{|x|y}{2t}} \left | \int ^{\infty}_{
- \frac{y}{2\sqrt{t}}} e^{i\beta
^2 -\frac{|x|\beta}{\sqrt{t}}} d\beta \right | =
\frac{1}{2\pi\sqrt{t}} e^{-\frac{|x|y}{2t}} \left | \int
^{\infty}_{\frac{y^2}{4t}} \frac{ e^{i\beta -\frac{|x|\sqrt{\beta
}}{\sqrt{t}}}}{\sqrt{\beta }} d\beta \right | .
\end{eqnarray*}
But $e^{-\frac{|x|\sqrt{\beta}}{t}}/\sqrt{\beta }$ is monotone
decreasing as $\beta \to \infty $. Standard results about oscillatory integrals,
then imply that
\[ \left | K_t (x,y)\right | \leq \frac{1}{\pi \sqrt{t}}
e^{\frac{-|x|y}{2t}} e^{-\frac{|x| |y|}{2t}}\left
(\frac{|y|}{2\sqrt{t}}\right )^{-1} \leq Ct^{-1/2}\] if
$\frac{|y|}{2\sqrt{t}}\geq 1$. For $0\leq \frac{|y|}{2\sqrt{t}}\leq
1$,
\begin{eqnarray*}
\left | K_t (x,y)\right | &\leq & \frac{1}{2\pi \sqrt{t}}
e^{\frac{-|x|y}{2t}} \left ( \left | \int ^{\infty }_1 e^{i\beta
-\frac{|x|}{\sqrt{t}}\sqrt{\beta}}\frac{1}{\sqrt{\beta }} d \beta
\right | +\left |  \int ^{1 }_{\frac{y^2}{4t}}
 e^{i\beta
-\frac{|x|}{\sqrt{t}}\sqrt{\beta}}\frac{1}{\sqrt{\beta }}d\beta
\right |
\right ) \\ \\
&\leq & \frac{1}{2\pi\sqrt{t}}e^{-\frac{|x|y}{2t}}\left (
e^{\frac{-|x|}{\sqrt{t}}}+\int
^1_{\frac{y^2}{4t}}\frac{e^{  \frac{-|x|}{\sqrt{t}} \cdot \frac{|y|}{2\sqrt{t}} } }{\sqrt{\beta}}
d\beta \right ) \leq  \frac{1}{2\pi\sqrt{t}}\left (1+\int
^1_{\frac{y^2}{4t}}\frac{d\beta}{\sqrt{\beta}} \right )\leq
\frac{C}{\sqrt{t}}.
\end{eqnarray*}
Hence, if $y\leq 0$,
\[ \left | K_t (x,y) \right | \leq \frac{C}{\sqrt{t}} .\]
On the other hand, if $y>0$,
\begin{eqnarray*}
 \left |K_t (x,y)\right | & \leq & \frac{1}{\pi
\sqrt{t}} e^{-\frac{|x|y}{2t}} \left ( \left | \int
^0_{-\frac{y}{2\sqrt{t}}} e^{i \beta ^2 -\frac{|x|\beta}{\sqrt{t}}}
d\beta \right | +\left | \int ^{\infty }_0 e^{i \beta ^2
-\frac{|x|\beta}{\sqrt{t}}} d\beta \right | \right )=\frac{1}{\pi \sqrt{t}} e^{-\frac{|x|y}{2t}} (I_1 +I_2) ,
\end{eqnarray*}
where
\begin{eqnarray*}
I_2 &=&\left | \int ^{\infty }_0 e^{i \beta ^2
-\frac{|x|\beta}{\sqrt{t}}} d\beta \right | \leq C, \\ \\
I_1 &=& \frac12 \left | \int ^0_{\frac{y^2}{4t}}
\frac{e^{i\beta +\frac{|x|\sqrt{\beta}}{\sqrt{t}}}}{\sqrt{\beta}}d \beta \right |
 \leq  \frac12 \left | \int ^{\frac{y^2}{4t}}_0 \frac{\cos \beta
 e^{\frac{|x|\sqrt{\beta}}{\sqrt{t}}}}{\sqrt{\beta}}d\beta \right | +
\frac12 \left | \int ^{\frac{y^2}{4t}}_0 \frac{\sin \beta
e^{\frac{|x|\sqrt{\beta}}{\sqrt{t}}}}{\sqrt{\beta}}d\beta \right |.
\end{eqnarray*}
If $\frac{y^2}{4t}\leq 2\pi$, then
$ |I_1|\leq Ce^{\frac{|x|y}{2\sqrt{t}}}.$
 If $\frac{y^2}{4t}> 2\pi$, let $k_0=\left \lfloor \frac{y^2}{8\pi t}
 \right \rfloor $ and obtain
 \begin{eqnarray*}
 \left |I_1\right | &\leq & \frac12 \left | \sum _{k=0}^{k=k_0-1}
 \int ^{2(k+1)\pi}_{2k\pi} \frac{\cos \beta
 e^{\frac{|x|\sqrt{\beta}}{\sqrt{t}}}}{\sqrt{\beta}} d\beta \right | +
 \frac12 \left | \sum _{k=0}^{k=k_0-1}
 \int ^{2(k+1)\pi}_{2k\pi} \frac{\sin  \beta
 e^{\frac{|x|\sqrt{\beta}}{\sqrt{t}}}}{\sqrt{\beta}} d\beta \right |
 \\ \\ & & \qquad \qquad + \frac12 \left |
 \int ^{\frac{y^2}{4t}}_{2k_0\pi} \frac{\cos \beta
 e^{\frac{|x|\sqrt{\beta}}{\sqrt{t}}}}{\sqrt{\beta}} d\beta \right | +
 \frac12 \left |
 \int ^{\frac{y^2}{4t}}_{2k_0\pi} \frac{\sin  \beta
 e^{\frac{|x|\sqrt{\beta}}{\sqrt{t}}}}{\sqrt{\beta}} d\beta \right |
 \\ \\
 &=& II_1 +II_2 +II_3 +II_4 .
 \end{eqnarray*}
 It is clear  that $
 |II_3 |+ | II_4| \leq Ce^{\frac{|x|y}{2\sqrt{t}}} .$
  The integral $II_2$ is now analyzed; $II_1$ can be treated similarly. First, notice that
 \begin{eqnarray*}
 |II_2| &=& \frac12 \left | \sum ^{k_0-1}_{k=0} \int
 ^{2(k+1)\pi}_{2k\pi} \frac{\sin \beta }{\sqrt{\beta}}
 e^{\frac{|x|\sqrt{\beta}}{\sqrt{t}}} d\beta \right | \\ \\
 &=&\frac12 \left | \sum ^{k_0-1}_{k=0} \int ^{2\pi}_0 \frac{\sin
 \beta}{\sqrt{2k\pi +\beta}}e ^{\frac{|x|}{\sqrt{t}}\sqrt{2k\pi
 +\beta}} d\beta \right |\\ \\
 &=& \frac12 \left |\sum ^{k_0-1}_{k=0} \left ( \int ^{\pi}_0 \frac{\sin
 \beta}{\sqrt{2k\pi +\beta}}e ^{\frac{|x|}{\sqrt{t}}\sqrt{2k\pi
 +\beta}} d\beta -\int ^{\pi}_0 \frac{\sin
 \beta}{\sqrt{2k\pi +\pi +\beta}}e ^{\frac{|x|}{\sqrt{t}}\sqrt{2k\pi
 +\pi
 +\beta}} d\beta \right )\right |.
 \end{eqnarray*}
 Since
 \begin{eqnarray*}
 \frac{\partial}{\partial u} \left ( \frac{1}{\sqrt{u}}e^{\frac{|x|\sqrt{u}}{\sqrt{t}}}\right )
 &=& \frac{x}{2u\sqrt{t}}e^{\frac{|x|\sqrt{u}}{\sqrt{t}}}-
 \frac{1}{2\sqrt{u^3}}e^{\frac{|x|\sqrt{u}}{\sqrt{t}}}=
 \frac{e^{\frac{|x|\sqrt{u}}{\sqrt{t}}}(\frac{|x|\sqrt{u}}{\sqrt{t}}-1)}{2u^{\frac32}}\\
 \\
 &\;\;& \left \{ \begin{array}{ll} >0,&  \quad \mbox{if }\
 \frac{|x|\sqrt{u}}{\sqrt{t}}>1, \\ <0\, , & \quad \mbox{if }\
 \frac{|x|\sqrt{u}}{\sqrt{t}}<1\, ,\end{array} \right.
 \end{eqnarray*}
if $a_k$ is defined by
\[
 a_k =\int ^{\pi}_0 \frac{\sin
 \beta}{\sqrt{k\pi +\beta}}e ^{\frac{|x|}{\sqrt{t}}\sqrt{k\pi
 +\beta}} d\beta ,
\]
 then there is an $N\geq 0$ such that $a_k$ is increasing in $k$ if
 $k>N$ and decreasing if $k\leq N$. In consequence, it transpires that
 \begin{eqnarray*}
 |II_2| &\leq & \frac12 \left | \sum ^N_{k=0} a_k (-1)^k \right | +
 \frac12 \left | \sum ^{2k_0-1}_{k=N+1} (-1)^k a_k \right | \\ \\
 &\leq & \frac12 \left (|a_0|+|a_N|\right ) + \frac12
 (|a_N|+|a_{2k_0-1}| ) \\ \\
 &\leq & C \Bigg ( \left | \int ^{\pi}_0 \frac{\sin
 \beta}{\sqrt{\beta}}e ^{\frac{|x|}{\sqrt{t}}\sqrt{\beta}} d\beta
 \right | + \left | \int ^{\pi}_0 \frac{\sin
 \beta}{\sqrt{N\pi +\beta}}e ^{\frac{|x|}{\sqrt{t}}\sqrt{N\pi
 +\beta}} d\beta \right |     \\ \\
 &&\qquad + \left | \int ^{\pi}_0 \frac{\sin
 \beta}{\sqrt{(2k_0-1)\pi +\beta}}e ^{\frac{|x|}{\sqrt{t}}\sqrt{(2k_0-1)\pi
 +\beta}} d\beta \right | \Bigg )\\ \\
 &\leq & Ce^{\frac{|x|}{\sqrt{t}}\sqrt{2k_0\pi}}\leq  C e^{\frac{|x|y}{2t}}.
 \end{eqnarray*}
The integral $II_1$ has a similar bound, whence
\[ \left |K_t (x,y)\right | \leq \frac{C}{\pi
\sqrt{t}}e^{-\frac{|x|y}{2t}}\left ( e^{\frac{|x|y}{2t}}+1\right ) \leq
\frac{C}{\sqrt{t}} \, \] for any $x, \  y\in \R$ and $t>0$.
Similar remarks apply to  $K_{-t} (x,y)$ so that for all $t\in \R \backslash \{0\}$,
$ \left |K_t (x,y)\right |  \leq
\frac{C}{\sqrt{|t|}}.$ This completes the proof of the Claim. \quad $\Box$

\medskip

To prove   inequality (\ref{2.9}), let $\K ( t ) \psi_h  = \intR
\psi_h  (y) K_t ( x, y ) dy $. The result of the Claim yields
$$
\| \K (t)\psi_h  \|_{L^{\infty }(\R)} \leq C |t|^{-\frac{1}{2}} \|
\psi_h \|_{L^{1} (\R)}.
$$
Also, Proposition 2.2.3 in \cite{cazenave} provides the inequality
$$
\| \K (t)\psi_h  \|_{L^2(\R)} \leq C  \| \psi_h \|_{L^{2} (\R)}.
$$
The Riesz-Thorin interpolation theorem then implies that
$$
\| \K (t)\psi_h  \|_{L^p(\R)} \leq C |t|^{-(\frac{1}{2} - \frac{1}{p})}
\| \psi_h \|_{L^{p'} (\R)}\, ,
$$
where $p'$ is the index conjugate to $p\,$ as before.
From this, there follows the inequality
\begin{align*}
\left \| \int^t_0 \K ( t - \tau ) f(\tau ) d\tau\right\|_{L^r(\R)}
& \leq C \int^T _0 |t- \tau |^{- (\frac{1}{2} - \frac{1}{r})}\|
f(\tau ) \|_{L^{r'} ( \R )}d\tau \\
& \leq C \int^T _0 |t - \tau |^{- \frac{2}{q} }\| f(\tau )
\|_{L^{r'} ( \R )}d\tau,
\end{align*}
valid for any $f ( \cdot, t) = f(t) \in L^{q'}(
(0, T) , L^{r'} (\R) ) $.
The Riesz potential inequalities (see
\cite{stein}, Theorem 1, p. 119)  then imply that
\begin{equation}
 \left \| \int^t_0 \K ( t - \tau ) f(\tau )  d\tau \right
\|_{L^q((0, T), L^r(\mathbb{R}))} \leq C \| f \|_{L^{q'} ((0, T),
L^{r'} ( \mathbb{R}))}. \label{lpq}
\end{equation}
A similar estimate holds for $\int^T_0 \K ( t - \tau ) f(\tau )
d\tau$.

Now, compute the $L^2(\mathbb{R})$--norm of the function
$$
\K_1 (y)  = \int^T_0 \int^{+\infty}_{-\infty}\overline {K_t(x,y)} f( x, t)
dx dt\, ,
$$
{\it viz.}
\begin{align*}
\| \K_1 (y) \|_{L^2(\mathbb{R})} ^2 =\intR \left ( \int^T_0
\int^T_0\int^{+\infty}_{-\infty}\int^{+\infty}_{-\infty} \overline {K_s(x,y)} f( x,
s) K_\sigma (w, y ) \overline { f ( w, \sigma)} \, dx\,dw\,ds\,
d\sigma\right ) dy.
\end{align*}
Note that
\begin{eqnarray*}
&&\intR \int^\infty_{-\infty} f (x, s ) \overline {K_{s} ( x, y )} dx
\int^\infty_{-\infty}
\overline{f (w, \sigma )} K_{\sigma} ( w, y ) \, dw  \, dy \\ \\
&&\quad= \int^\infty_{-\infty} f( x , s ) \int^\infty_{-\infty} \overline{ f ( w,
\sigma )} \intR \overline{K _ { s} ( x, y ) }  {K _ {\sigma } ( w
, y)}  \, dy   \, dw  \, dx\\ \\
&&\quad= \int^\infty_{-\infty}\int^\infty_{-\infty} f( x , s )  \overline{ f ( w,
\sigma )} K _ { s,\sigma } ( x, w)  \, dw  \, dx.
\end{eqnarray*}
The inequality (\ref{claim}) implies that
\begin{eqnarray*}
&&K _ { s,\sigma } (x, w)  =  \intR \overline{K _ { s}
( x, y ) } K _ {\sigma } ( w , y) dy \\ \\
& = &\frac{1}{\pi^2 } \int^\infty_0\intR \int^\infty_0 {e ^{ -i
\tilde \beta^2  s - \tilde \beta |x| + i y \tilde \beta  } e ^{ i
\beta^2  \sigma  - \beta |w|- i
y  \beta } }   \, d \beta  \, dy  \, d\tilde \beta\\ \\
& = &\frac{2}{\pi} \int^\infty_0  e ^{ - i \beta^2 (s-\sigma)  -
\beta ( |x|+ |w| ) } d\beta \leq \frac{C}{\sqrt{|s-\sigma|}}
\end{eqnarray*}
for $s \not = \sigma$, where the  constant $C$ is independent of $x, w \in \R$.                                                                                                                                                                                                                                                                                                                                                                                         . Rewrite $\| \K_1 (y) \|_{L^2(\R)}
^2$ as
\begin{align*}
\| \K_1 (y) \|_{L^2(\R)} ^2 =& \int^T_0 \int^T_0
\int_{-\infty}^{+\infty}\int_{-\infty}^{+\infty} f( x , s ) \overline{ f ( w, \sigma
)} K _ {
s,\sigma } ( x, w)  \, dw  \, dx  \, d\sigma  \, ds \\ \\
&\quad  = \int^T_0 \int_{-\infty}^{+\infty} f( x , s ) \int^T_0
\int_{-\infty}^{+\infty} \overline{ f ( w, \sigma ) }K _ {  s,\sigma } ( x,
w)  \, dw   \, d \sigma  \, dx  \, ds\, .
\end{align*}
Then,
using the procedure described for proving (\ref{lpq}), it is inferred that
$$
\left \| \int^T_0 \int_{-\infty}^{+\infty}  \overline{ f ( w, \sigma ) }K _
{ s,\sigma } ( x, w) dw  d \sigma \right \|_{L^q((0, T), L^r(\R))}
\leq C \| f \|_{L^{q'} ((0, T), L^{r'} ( \R))}\, ,
$$
which in turns gives
\begin{align*}
\| \K_1 (y) \|_{L^2(\mathbb{R})} ^2 & \leq C \|
f\|_{L^{q'}((0,T),L^{r'}(\R))} \Bigg (\left \| \int^T_0 \intR
\overline{ f ( w, \sigma ) }K _ { s,\sigma } ( x, w) dw  d \sigma
\right \|_{L^q((0, T), L^r(\mathbb{R}))}\Bigg
) \\ \\
& \qquad \leq C \| f \|^2_{L^{q'} ((0, T), L^{r'} ( \mathbb{R}))}.
\end{align*}

Finally, consider the integral
$$
\intR \Big ( \K (t)\psi_h  , \phi (\cdot, t)\Big )_{L^2} dt =  \intR \left (
\int_{-\infty}^{+\infty} \intR K_t ( x, y) \psi_h  ( y ) dy \overline{\phi ( x,
t) } dx \right ) dt\, ,
$$
where $\phi ( x, t) \in C_c ( [ 0, T], {\cal D} (\mathbb{R})),
\psi_h  \in L^2 ( \mathbb{R})$. Applying the just obtained estimates yields
\begin{align*}
\intR ( \K (t)\psi_h  , \phi (\cdot, t))_{L^2(\R)} dt & = \intR \psi_h
(y ) \overline { \intR \int_{-\infty}^{+\infty} \overline{K_t ( x, y ) }
\phi ( x, t) dx dt}
dy \\
& \leq \| \psi_h  \|_{L^2} \left \| \int^T_0  \int_{-\infty}^{+\infty}
\overline{K_t ( x,
y ) } \phi ( x, t) dx dt\right \| _{L_y ^2 ( \mathbb{R}) }\\
& \leq C \| \psi_h  \|_{L^2(\R)} \| \phi \|_{L^{q'} ((0, T), L^{r'} (
\mathbb{R}))}\, .
\end{align*}
By duality,
$
\| \K (t)\psi_h  \|_{L^{q} ((0, T), L^{r} ( \mathbb{R}))} \leq C \|
\psi_h  \|_{L^2(\mathbb{R})}\, ,
$
which gives (\ref{2.7}) with $s=0$.   Since
\[\partial _x  [W_{b,2} (t)h](x)= \frac{1}{\pi} \int ^{\infty}_0 e^{i\beta ^2
t-\beta
|x|} \  sign (x) \ \frac{\beta ^2} {b-a\beta}\hat{h}(i\beta ^2) d\beta, \]
the same argument suffices to show  that (\ref{2.7}) holds for $s=1$. When $0<s< 1$, the
relevant estimate follows by interpolation. The inequality (\ref{2.8}) is a special case of (\ref{2.7}) and (\ref{2.8-1}) is straightforwardly
obtained using a classical trace argument and the Fourier transform. $\Box$

\medskip
The following estimates of the temporal regularity of $W_{bdr}$
will also be helpful.
\begin{proposition} \label{prop2.12}
Let $(q,r)$ be a given admissible pair, $T>0$ and $ s \geq 0$.
  For any $h\in H^{\frac14+s} (\R^+)$, the correspondence $t\mapsto
\frac{\partial ^s W_{bdr} (t)}{\partial  t^s}h $ belongs to the space
\[ L^q (0,T; L^r (\R^+))\cap C([0,T], L^2 (\R^+)) \]
and there exists a constant $C $ such that
\begin{equation}
\left \| \frac{\partial^s W_{bdr} (\cdot )}{\partial t^s} h\right \|
_{L^q (0,T; L^r(\R^+))} \leq C \| h\| _{H^{\frac14 + {s}} (\R^+)}\, .
\end{equation} In particular, for $h\in H^{\frac14+ {s}}(\R^+)$,
\begin{equation}
\sup _{0<t<T} \left \| \frac{\partial^s W_{bdr} (\cdot )}{\partial
t^s} h\right \| _{L^2(\R^+)} \leq C \| h\| _{H^{\frac14+ {s}}
(\R^+)}
\end{equation}
and
\begin{equation}
\sup _{0<x<\infty} \left \| \frac{\partial^s W_{bdr} (\cdot )}{\partial
t^s} h\right \| _{H^{\frac14} _t(0,T)} \leq C \| h\| _{H^{\frac14+
{s}} (\R^+)}\, .
\end{equation}

\end{proposition}

\noindent{\it Proof:} As above, we only have to study $W_{b,2} h$. It
is straightforward to calculate that
\begin{align*}
\frac{\partial W_{b,2} (t )h}{\partial t}   &= \frac{i}{\pi} \int
^{\infty}_0 e^{i\beta ^2 t-\beta x}{\beta^3} \hat{h}(i\beta ^2) d\beta
= \frac{i}{\pi} \int ^{\infty }_0 e^{i\beta ^2t-\beta x} \hat{\psi_1
} (\beta ) d\beta = i\int ^{\infty }_{-\infty } \psi  _1 (y)K_t(x,y)dy.
\end{align*}
It follows immediately that
$$
\left \| \frac{\partial W_{b,2} (t )}{\partial t} h\right \| _{L^q
(0,T;;L^r(\R^+))} \leq C \| \psi_1\|_{L^2(\R)} \leq  C \|
h\| _{H^{\frac14 + 1} (\mathbb{R}^+)}.
$$
A similar proof holds for all integers $s\geq 0$. The general case then
follows by interpolation. Since there are no boundadry conditions
involved in the argument, we do not run into trouble when the interpolation index
is equal to $\frac12$.  In particular, the Sobolev space
 $H^{\frac12} (\mathbb{R}^+) $ is  the mid-point  interpolation space between $L^2 (\mathbb{R}^+)$ and $H^1 (\mathbb{R}^+)$ in this case.
$\Box$

\medskip
Note that from the equation $ iu_t +u_{xx}=0,$ one
$t$-derivative of $u$ is equivalent to two $x$-derivatives
of $u$. The following proposition holds as a corollary of this observation.
\begin{proposition}  \label{prop2.13}
Let $(q,r)$ be a given admissible pair, $T>0$  and $ s\geq 0$. There
exists a constant $C>0$ such that
  for any $h\in H^{\frac14+s} (\R^+)$, $u=W_{bdr} (t) h$ satisfies
   \[ \|u\| _{L^q_t (0,T; W^{s,r}_x(\R^+))}+  \sup _{0<t< T} \|u(\cdot, t)\| _{H^s(\R^+)}
   + \sup_{0<x<\infty} \|\partial ^j_xu(x, \cdot)\|_{H^{\frac{2s+1-2j}{4}}(0,T)}\leq
   C\|h \|_{H^{\frac{2s+1}{4}}(\R^+)} \]
   \end{proposition}
   for $j=0,1$ and $2s+1-2j \geq 0$.

\medskip
Finally,  we consider the IBVP (\ref{iforcing}).  The next proposition follows
readily from Propositions \ref{prop2.9}-\ref{prop2.13}.
\begin{proposition}\label{linear-half} Let $T>0$ and $0\leq s < \frac52$ be given. Assume $f\in L^1 (0,T; H^s (\R^+))$,  $\phi \in H^s
(\mathbb{R}^+)$, $h\in H^{\frac{2s+1}{4} }(0,T)$ and $\phi (0) =h(0)$ if $\ \frac12 < s< \frac52 $. Then there exists a constant $C>0$ depending only on $s$ such that the solution $u$ of the IBVP (\ref{iforcing})
respects the inequality
\begin{eqnarray*} & & \sup_{0\leq t\leq T} \| u(\cdot, t)\|_{H^s(\R^+)} +\sup _{0<x< \infty}
\| u(x, \cdot)\|_{H^{\frac{2s+1}{4}} (0,T)} + \|u\| _{L^q(0,T;W^{s,r}(\R^+))}   \\
& & \qquad \qquad \leq C\left (\|\phi \|_{H^s(\R^+)} +\|h\|_{H^{\frac{2s+1}{4}} (0,T)} +\|f\|_{L^1 (0,T; H^s (\R^+))}  \right )\, ,
\end{eqnarray*}
 where  $(q,r)$ is  any admissible pair.
\end{proposition}

\subsection{Local well-posedness}
In this subsection, the  local well-posedness of the full nonlinear problem
\begin{align}\label{n-1}
\left \{ \begin{array}{l} iu_t +u_{xx} + \lambda u |u |^{p-2} =0,
\qquad x\in \mathbb{R}^+, \ t\in (0,T),
\\ \\ u(x,0) =\phi (x) , \qquad  u(0,t) = h(t)\, ,
\end{array}\right.
\end{align}
is the topic of conversation. Let ${\phi}^*  = B \phi$ be an extension of $\phi$ from $\R^+$ to $\R$  as before, with
\[
\|{\phi}^* ||_{H^s (\R)} \leq C_s \| \phi \|_{H^s (\R^+)} .
\]
Suppose $ 0\leq s < 5/2$ and let the operator  ${\cal W}_{bdr}$  be as  introduced in Section 2.  Rewrite (\ref{n-1})
as an  integral equation on the domain $(x, t) \in \R\times \R^+$, {\it viz.}
\begin{equation}\label{int}
u(t)=W (t) {\phi}^* +{\cal W}_{bdr} (t)\big ( h(t)-
g_{\phi}(t)\big ) -i\lambda \int ^t_0 W(t-\tau )|u|^{p-2}
 u (\tau) d\tau - {\cal W}_{bdr} (t) f_u(t)\,
\end{equation}
where $W(t) = W_{\mathbb{R}}(t)$ and $g_{\phi}(t), f_u(t)$ are the trace of $W(t)\phi^*$
and $-i\lambda \int ^t_0  W(t-\tau )|u|^{p-2} u (\tau) d\tau$  at
$x=0$.   That  is to say,
\[ g_{\phi}(t)=  \left. W(t) \tilde{\phi} \right  |_{x=0},  \quad \quad   f_u(t)= \left.  -i\lambda
\int ^t_0  W(t-\tau )|u|^{p-2} u (\tau) d\tau \right  |_{x=0}.\]

\begin{proposition}\label{prop3.1}  Assume $$0\leq s<\frac12\quad\mbox{and}\quad 3\leq {p}
<  \frac{6-4s}{1-2s}.$$  Let $(\gamma, \rho)$ be the admissible pair
defined by \[ \rho= \frac{p}{1+s{(p-2)}}, \quad \gamma =
\frac{4p}{{(p-2)}(1-2s)} .\]
For any given $\phi \in
H^s(\R^+)$ and $h\in H^{\frac{2s+1}{4}}(0, T)$, there exists a $T_{max}$ with $0 < T_{max} \leq T$
such that the integral equation (\ref{int}) admits a unique solution
$u\in C([0,T_{\max}); H^s (\R))$
 satisfying
\begin{equation}\label{cond-11}
u\in L^{\gamma} _{loc}([0, T_{\max}); W^{s, \rho}(\R)).
\end{equation}
Moreover, this solution possesses the following additional properties:
\begin{itemize}
\item[(i)] The solution $u\in  L^{q} _{loc}([0, T_{\max}); W^{s, r}(\R))$ for
every admissible pair $(q,r)$.
\item[(ii)] The solution $u$ depends   continuously on $\phi $ and $h$ in the sense that
if $\phi _n \to \phi$  in $H^s (\R^+)$ and $h_n \to h$ in
$H^{\frac{2s+1}{4}}(\R^+)$, then,  for any $T$ with  $0<T<T_{max}$,  the
corresponding solutions  $u_n $  tend to $u$  in $ C([0,T]; H^s (\R
))  $  as $n\to \infty $.
\item[(iii)] If $3\leq {p} < \frac{6-4s}{1-2s}$ and  $T_{max} < +\infty $, then
\[ \lim _{t\to T_{max}} \|u(\cdot, t)\|_{H^s (\R)} = +\infty.\]
\end{itemize}
\end{proposition}
\begin{proposition}\label{prop3.2} Let $\frac12 < s< \frac52 $  and $[s]\leq {p-2} $
 be given. For any $\phi \in H^s (\R^+) $ and $h\in H^{\frac{2s+1}{4}}(\R ^+)$ satisfying the
 compatibility condition
\[ \phi (0)=h(0) ,\]
there exists a $T_{max} >0$ such that the integral equation
(\ref{int}) admits a unique solution $$u\in C([0, T_{max}); H^s
(\R^+)).$$ Moreover, the solution $u$ possesses the following
properties:
\begin{itemize}
\item[(i)] The solution $u$  belongs to the space  $ L^{\infty}_x (\R^+;
H^{\frac{2s+1}{4}}_t (\R))$.
\item[(ii)] The solution $u$ depends on $\phi $ and $h$  continuously in the sense that
if $\phi _n \to \phi$  in $H^s (\R^+)$ and $h_n \to h$ in
$H^{\frac{2s+1}{4}}(\R^+)$, then,  for any $T$ with  $0<T<T_{max}$,  the
corresponding solutions  $u_n $  tends to $u$  in $ C([0,T]; H^s (\R^+
))\cap  L^{\infty}_x (\R^+;  H^{\frac{2s+1}{4}}_t (\R^+))$  as $n\to \infty $.
\item[(iii)]
If $T_{max} < +\infty $, then
\[ \lim _{t\to T_{max}} \|u(\cdot, t)\|_{H^s (\R^+)} = +\infty.\]
\end{itemize}
\end{proposition}

The proofs of Propositions \ref{prop3.1} and \ref{prop3.2}
follow just as does the local existence theory  laid out in Holmer \cite{holmer}.
The chain rule and product rule for fractional derivatives and the propositions in the
last subsection provide the necessary estimates for applying the contraction mapping theorem
to the right-hand side of (\ref{int}). The details are omitted.

\begin{remark}\label{rem3.3.1} Proposition \ref{prop3.2} also holds for $  s \in (\frac52, \frac92) $ if the following compatibility conditions are satisfied;
$$
\phi ( 0 ) = h (0) \, ,\quad h_t( 0) = i \phi _{xx}(0) + i \lambda |\phi (0) |^{p-2} \phi (0) \, .
$$
\end{remark}
The only difference from the proof of the local existence in Holmer \cite{holmer} is to use the function space
$$
C( (0, T) ; H^s_x (\R^+ ))\cap C( \R^+_x, H^{\frac{2s+1}4} (0, T) ) \cap C^1_t  ( (0, T) ; H^{s-2}_x (\R^+ ))\, .
$$
Note again that one $t$-derivative of $u$ corresponds to two $x$-derivatives of $u$.

\begin{remark}\label{rem3.3} In case of $s=1$ or $s=2$,  the assumption $\lfloor s \rfloor
\leq {p-2} < +\infty$ is not needed.   The result of
 Proposition \ref{prop3.2} holds for any $p$ with $p > 2$ in these situations.
\end{remark}

\begin{remark}\label{rem3.4}  According to Proposition \ref{prop3.1}, for $\phi \in H^s (\R^+)$,
$h\in H^{\frac{2s+1}4}_{loc}(\R^+)$  with $0\leq s< \frac12 $, there exists a $T_{max}$ depending only
on $s$ such that
 the corresponding solution $u\in C([0, T_{max}); H^s (\R))$ blows up at $T_{max}$, i.e.,
\[ \lim _{t\to T_{max}}\| u(\cdot, t)\|_{H^s(\R)} =+\infty\]
if  $T_{max} <\infty $. However, if $(\phi, h)$ also belongs to the space
$H^2 (\R^+)\times H^{\frac54}_{loc} (\R^+)$,  then  by Proposition \ref{prop3.2},
 there exists a $T_{max}^*>0$
such that  $ u \in C([0, T^*_{max});  H^2 (\R)) $ and
 \[ \lim _{t\to T_{max}} \|u(\cdot, t)\|_{H^2 (\R)} = +\infty \]
if $T^*_{max} <\infty $.  It is obviously the case that  $T^*_{max} \leq T_{max}$.
Is it true that $T^*_{max}= T_{max}?$
This is a well-known regularity issue (see \cite{cazenave}).   For the pure Cauchy problem (\ref{n-1-R}),
 the answer is positive.  The  same proof can be applied  to the IBVP considered here  to show that  $T^*_{max}=T_{max}$.
\end{remark}

A solution of the integral equation (\ref{int})  on $\R$ as given in
Propositions \ref{prop3.1} and \ref{prop3.2}, when restricted to
$\R^+$,  is a distributional solution of the IBVP (\ref{n-1}) with
strong traces. However, as the IBVP (\ref{n-1}) can be converted to other
 integral  equations similar to (\ref{int})  on $\R$, whose solutions,
when restricted to $\R^+$,  yield distributional solutions to the
IBVP (\ref{n-1}), the following  question arises naturally.

\smallskip
{\it Are solutions of the various integral equations on $\R$  equal to each other when restricted  to $\R^+$?}
In other words, Propositions \ref{prop3.1} and \ref{prop3.2}  lead to the existence of  distributional
solutions with strong traces for the IBVP (\ref{n-1}). As for its uniqueness,
in the case of $s>\frac12$, since the space $H^s (\R^+)$ is continuously  imbedded into the space
$L^{\infty} (\R^+)$, it is straightforward to ascertain that the IBVP (\ref{n-1}) admits at most one distributional
solution with strong traces
in the space $C([0,T]; H^s (\R^+))$.   The following well-posedness theory for
 the IBVP (\ref{n-1}) results as a corollary of Proposition \ref{prop3.2}.

\begin{cor}\label{cor3.3} Let $\frac12 < s< \frac52 $  and $3\leq {p} < +\infty $
 be given. For any $\phi \in H^s (\R^+) $ and $h\in H^{\frac{2s+1}{4}}(\R ^+)$ satisfying the
 compatibility condition
\[ \phi (0)=h(0) ,\]
there exists a $T_{max} >0$ such that (\ref{int}) admits a unique solution $u\in C([0, T_{max}); H^s (\R ^+))$.
Additionally,  the solution $u\in L^{\infty}_x (\R;  H^{\frac{2s+1}{4}}_t (\R ^+))$ and
if $T_{max} < +\infty $, then
\[ \lim _{t\to T_{max}} \|u(\cdot, t)\|_{H^s (\R ^+)} = +\infty.\]
Moreover,  the solution $u$ depends continuously on $\phi $ and $h$   in the sense that
if $\phi _n \to \phi$  in $H^s (\R^+)$ and $h_n \to h$ in $H^{\frac{2s+1}{4}}(\R^+)$,
then,  for any $0<T<T_{max}$,  the corresponding solutions  $u_n $  tend to $u$  in
 $ C([0,T]; H^s (\R ^+ )) $  as $n\to \infty $.
\end{cor}

The  uniqueness of the IBVP (\ref{n-1}) in the space $C([0,T]; H^s (\R^+)) $
remains open in case            $0\leq s< \frac12.$  To resolve this issue,
 we first show that the solution given in
Proposition \ref{prop3.1}, when restricted on $\R^+$ is a mild solution of the IBVP (\ref{n-1}).

\begin{proposition}\label{prop3.3} Let $0\leq s<\frac12$ be given and assume that $ 3\leq {p}
< \frac{6-4s}{1-2s}$.  For any given  $\phi \in H^s(\R^+)$ and
$h\in H^{\frac{2s+1}{4}}_{loc}(\R^+)$,
there exists a $T_{max}>0$ such that the IBVP (\ref{n-1}) admits a mild solution $u\in C([0,T_{\max}); H^s (\R))$.
\end{proposition}

\noindent
{\bf Proof:} It suffices to show that for $0\leq s < 1/2$ the solution $u\in C([0,T_{max}); H^s (\R))$
of (\ref{int}) given by Proposition \ref{prop3.1}, when restricted to $\R^+\!$, is a mild
solution of the IBVP (\ref{n-1}).   To this end, let $(\phi _n, h_n)\in H^2
(\R^+)\times H^{\frac54}_{loc} (\R^+)$ with $\phi _n (0)=h_n (0)$ and
\[
\lim _{n\to \infty} \| (\phi _n , h_n)- (\phi , h)\|_{H^s (\R^+)\times H^{\frac{2s+1}4}_{loc} (\R^+)} =0 .
\]
Then by Proposition \ref{prop3.1}, there exists $u_n \in C([0,T_{max}); H^2 (\R^+))$ solving
the integral equation (\ref{int}) with $(\phi,  h)$ replaced by $(\phi _n , h_n)$. Moreover,  $u_n$
tends to $u$ in the space $C([0,T]; H^s (\R^+))$ as $n\to \infty $ for any $T<T_{max}$.
 According to Remarks \ref{rem3.3} and \ref{rem3.4}, when restricted to $\R^+$,
$u_n$ lies in $ C([0,T_{max}); H^2 (\R^+))$ and it  solves the IBVP (\ref{n-1}).
In particular,  $u_n$ tends to $u$ in the space $C([0,T]; H^s (\R^+))$ as $n\to \infty$
 for any $T<T_{max}$.  Thus, the solution $u$
of (\ref{int}) when restricted to $\R^+$ is a mild solution of the IBVP (\ref{n-1}).
 $\Box$

\medskip
Next, we show that  the IBVP (\ref{n-1}) admits at most one mild solution.
\begin{proposition}\label{prop3.7} Assume that $s$ and $p$ are such that
\[ 0\leq s < \frac12 \quad {\rm and} \quad  3\leq {p} < \frac{6-4s}{1-2s} .\]
For any $\phi \in H^s(\R^+)$, $h\in H^{\frac{2s+1}4}_{loc} (\R^+)$,
the IBVP (\ref{n-1})
admits at most one mild solution.
\end{proposition}

\medskip
\noindent {\bf Proof}:   Suppose that
for a given $\phi\in H^s (\R^+)$ and $h\in H^{\frac{2s+1}{4}}_{loc}
(\R^+)$, the IBVP (\ref{n-1}) admits two mild solutions $u$ and $v$
which lie in
in the space $C([0,T']; H^s (\R^+))$ for some $T'>0$. By
definition, there exist two sequences $\{ u_n\} $ and $\{v_n\}$ in
the space $C([0,T']; H^2 (\R^+))$ such that both $ u_n  $ and $v_n $
solve the equation in (\ref{n-1}) for $n=1,2,\cdots , $ and if
\[
u_n(x,0)=\phi _n (x), \ v_n(x,0)= \psi _n (x), \ u_n(0,t)= h_n
(t), \ v_n (0,t)= g_n (t),
\]
then as $n\to \infty$,
\[ u_n \to u , \quad v_n\to v \quad \mbox{in} \quad C([0,T']; H^s(\R^+)),\quad
\phi _n \to \phi, \quad \psi _n \to \phi \quad \mbox{in} \quad  H^s(\R^+)\]
and
\[ g_n \to h , \quad h_n \to h \quad \mbox{in} \quad  H^{\frac{2s+1}{4}}(0,T) .\]
Let $u_n^*$,  $v^*_n$ and $w$ be the solutions of the integral equation  (\ref{int})
corresponding to $(\phi _n, h_n)$, $(\psi _n, g_n)$ and $(\phi , h)$,
respectively, given by Proposition \ref{prop3.1} (restricted to $\R^+$). It follows that
$u_n^* $,  $v_n^*$   and $w$ lie in $C([0,T;
H^s (\R^+))\cap L^q(0,T; W^{s,r} (\R^+))$ for some $T>0$.  Then, by Proposition \ref{prop3.2}
and  Remarks \ref{rem3.3} and \ref{rem3.4},
$u_n^* $ and $v_n^*$ are in $C([0,T]; H^2(\R^+)) $.  Note that the time
interval over which $u_n^*$ and $v_n^*$ exist in the space $H^2 (\R^+)$ is $(0,T)$ for any $n$,
as guaranteed by  Remark \ref{rem3.4}. By
the uniqueness result in Corollary \ref{cor3.3}, it must be the case that
\[ u_n\equiv u_n^* ,\qquad v_n\equiv v_n^*, \ n=1,2,\cdots  \, . \]
Since $(\phi _n, h_n)$ and $(\psi_n , g_n)$ are both convergent to
$(\phi, h)$ in $H^s(\R^+)\times H^{\frac{2s+1}{4}}(0,T)$,
it follows from Proposition \ref{prop3.1}  that both $u_n$ and $v_n$
converge to $w$ in $C([0,T]; H^s (\R^+))$. Consequently,
$  u\equiv v.$
The proof is complete. $\Box$
\medskip

The last result of the section summarizes the previous ruminations.
\begin{theorem}
Assume either
\[  3\leq {p} < \frac{6-4s}{1-2s}, \quad  0\leq s< \frac12 , \quad s=1,2 \]or

\[    \frac12 <s<\frac52 , \qquad  \lfloor s\rfloor < {p-2} < \infty.\]
For any $\phi \in H^s (\R^+)$ and $H^{\frac{2s+1}4}_{loc} (\R^+)$ satisfying $\phi (0)= h(0)$ if
$s> \frac12$, there exists a $T_{\max} >0$ such that
the IBVP (\ref{n-1}) admits a unique  mild solution $u\in C([0,T_{max}); H^s (\R^+))$.
Moreover, the solution $u$ has the following properties:
\begin{itemize}
\item[(i)]  The solution $u\in L^{\infty}_x (\R^+; H^{\frac{2s+1}{4}}_{loc} (\R^+))$.
\item[(ii)] The solution $u$ depends on $\phi $ and $h$ continuously in the  sense that
 if $\phi _n \to \phi$  in $H^s (\R^+)$ and $h_n \to h$ in $H^{\frac{2s+1}{4}-{loc}}(\R^+)$,
then,  for any $T$ with $0<T<T_{max}$,  the corresponding solutions  $u_n $  tend to $u$  in
 $ C([0,T]; H^s (\R ^+ ))$  as $n\to \infty $.

\end{itemize}
\end{theorem}

\section{The Schr\"odinger equation posed on a finite interval}
\setcounter{equation}{0}

In this section, consideration is given to  the well-posedness in $H^s (0,L)$ of the IBVP
\begin{align}\label{4.0-1}
\left \{ \begin{array}{l} iu_t + u_{xx} + \lambda |u|^{p-2} u =0,
\qquad x\in (0,L), \ t\in \R,
\\ \\ u(x,0) =\phi(x) , \qquad u(0,t) = h_1 (t), \qquad u(L, t) =h_2
(t)  ,\end{array}\right.
\end{align}
for the NLS
equation posed on a finite interval $(0,L)$.  Without loss of
generality, take $L=1$.

First, the homogeneous boundary-value problem
\begin{align}\label{0-1-1}
\left \{ \begin{array}{l} iu_t +u_{xx} + \lambda u |u |^{p-2} =0,
\qquad x\in (0,1), \ t\in \mathbb{R},
\\ \\ u(x,0) =\psi (x) , \qquad u(0,t) = 0, \quad u (1, t) = 0 ,
 \end{array}\right.
\end{align}
is discussed. The well-posedness of (\ref{0-1-1}) in $H^s (0,1)$  can be reduced to a special case  of the IVP
\begin{align}\label{0-2}
\left \{ \begin{array}{l} iu_t +u_{xx} + \lambda u |u|^{p-2} =0,
\qquad -1< x< 1, \ t\in \mathbb{R},
\\ \\ u(x,0) =\psi (x), \qquad
u(-1,t)=u(1,t),\qquad u_x (-1,t)=u_x (1,t), \end{array}\right.
\end{align}
of the NLS equation posed on the interval $(-1,1)$ with periodic boundary
conditions. Observe that solutions of the IVP (\ref{0-2}) are even (odd) in $x$ if
 $\psi $ is even (odd).
On the other hand, if $u$ is an odd function with respect to $x$ and
solves the IVP (\ref{0-2}), then its restriction to the interval
$(0,1)$  solves  the IBVP (\ref{0-1-1}) since the boundary conditions $u (0, t) = u(1, t) =0$
are automatically satisfied. Thus, the following
well-posedness result follows immediately from the known results for
(\ref{0-2}).
\begin{theorem} Assume that $3\leq p < \infty $ if $\lambda <0$ and
$3\leq p < 6$ if $\lambda >0$. Then, for any $s \in [0,\frac52)$ ($s$ not equal to $ \frac12$ or $\frac32$, see \eqref{s-assumption-2}),
 the IBVP (\ref{0-1-1}) is
unconditionally  locally well-posed in $H^s (0,1)$  under the conditions that
$\lfloor s\rfloor < p-2$ if $p$ is not an integer and $\psi (0) = \psi (1) = 0$ if
$ \frac12 < s < \frac52$.
\end{theorem}

Now, consider (\ref{4.0-1}) with nonhomogeneous boundary data.
This is analyzed in several stages.

\subsection{Linear problem}
First, consider the IBVP
\begin{align}\label{4.1-1-1}
\left \{ \begin{array}{l} iu_t +u_{xx} =0, \qquad x\in (0,1), \ t\in
\R,
\\ \\ u(x,0) =\phi(x) , \qquad u(0,t) = u(1, t) =
0,\end{array}\right.
\end{align}
for the linear Schr\"odinger
equation posed  on the finite interval $(0,1)$. According to standard semigroup theory, for any $\phi \in L^2
(0,1)$, the IBVP  admits a unique solution $u\in C(\R^+; L^2 (0,1))$
given by
$$
u(t) = W_0(t) \phi
$$
where $W_0(t)$ is the $C_0$-group in $L^2(0,1)$ generated by
the operator $ Av=iv''$ with domain ${\cal D}(A)= H^2(0,1)\cap
H^1_0 (0,1)$. Moreover, the solution of the following nonhomogeneous
problem
\begin{align}\label{4.1-1}
\left \{ \begin{array}{l} iu_t +u_{xx} =f, \qquad x\in (0,1), \ t\in
\R,
\\ \\ u(x,0) =0 , \qquad u(0,t) = u(1, t) =
0\end{array}\right.
\end{align}
can be expressed, via Duhamel's principle, as
\[ u(t)=- i \int ^t_0 W_0 (t-\tau )  f(\cdot, \tau) d\tau .\]
\begin{proposition}\label{pro3.1}
Let $0\leq s\leq 2 $  and $T>0$ be given.  Let \[ u(t)=W_0 (t) \phi,
\quad v(t)=\int ^t_0 W_0 (t-\tau)f(\cdot, \tau)dt\tau\] and \[ w(t)=\int
^t_0 W_0 (t-\tau ) g (\cdot ,\tau ) d\tau, \]
 with
 $\phi \in H^s (0,1)$, $f\in L^1 (0,T;H^s (0,1))$ and $g\in
 W^{\frac{s}{2},1}(0,T; L^2 (0,1))$
satisfying
\[\phi(0)=\phi (1)=0, \quad f(0,t)=f(1,t) \equiv 0 \] if $s>\frac12$.
Then, $u, \ v, \ w\in C([0,T]; H^s (0,1))$ and
\[ \|u\|_{C([0,T]; H^s (0,1))}\leq C_{T,s} \|\phi \|_{H^s (0,1)},  \]
\[ \|v\|_{C([0,T]; H^s (0,1))}\leq C_{T,s} \|f\|_{L^1(0,T; H^s (0,1))}  \]
and
\[ \|w\|_{C([0,T]; H^s (0,1))}\leq C_{T,s} \|g \|_{W^{\frac{s}{2},1} (0,T;L^2 (0,1))}, \]
 where  the constant $C_{T,s}$ depends only on $s$ and $T$.
\end{proposition}
{\bf Proof}: The cases $s=0$ and $s=2$ follow from  standard
semigroup theory.    When $0< s< 2$,  these inequalities
are follow from standard interpolation theory. $\Box$

 \medskip
In terms of Fourier sine series, the solution $u$ is given explicitly
by
$$
u (x, t) =\big [W_0(t) \phi\big ] (x) =  \sum_{n=1}^{+\infty} c_ne^{-i(n\pi)
^2 t } \sin (n\pi x ) \quad \mbox{where}\quad c_n =2
\int ^{1}_0\phi ( x ) \sin (n\pi x )  dx.
$$
This can be written in the complex form
\[ u(x,t)= \sum ^{\infty}_{n=-\infty}\tilde{c}_n e^{-i(n\pi)^2t+ i n\pi x}\]
where
\[ \tilde{c}_n = \left \{ \begin{array}{ll} c_n & \ \mbox{ if } \ n\geq
1, \\ 0& \ \mbox{ if } \ n=0, \\ -c_n & \ \mbox{ if } \ n\leq -1. \end{array} \right.
\]
In this form, it is clear that $u$ may be viewed as  the solution $u(x,t)$ of the  Cauchy
problem
\begin{equation}\label{4.3-1}
\left \{ \begin{array}{l} iu_t +u_{xx}=0, \ u(x,0)= \phi ^* (x), \
x\in (-1,1), \\ \\ u(-1,t)=u(1,t), \ u_x (-1,t)=u_x (1,t)\, ,
\end{array} \right.
\end{equation}
where $\phi ^* $ is the odd
extension of  $\phi $ from $(0,1)$ to $(-1,1)$. On the other hand,
if $u$ is a solution of (\ref{4.3-1}) and is also an odd function,
then its restriction to $(0,1)$ solves (\ref{4.1-1}). Thus
\[ \big [W_{\T} (t) \phi ^*\big ](x) = \big [W_0 (t) \phi \big ] (x), \quad x\in (0,1).\]
Here, $W_{\T} (t)$ is the $C_0$-group in $L^2 (\T)$
generated by the operator $A_{\T}$ in $L^2 (\T)$ with
domain ${\cal D}(A_{\T}) = H^2(\T)$. Consequently, the following
proposition follows from the theory developed in \cite{bourgain-1}.
\begin{proposition} \label{prop3.3-1}
Let $0\leq s< \frac12 $ and $T>0$ be given and
let $\Omega _T= (0,1)\times (0,T)$.
For any $\phi \in H^s (0,1)$, $u=W_0 (t) \phi \in L^4(\Omega _T)\cap
C([0,T]; H^s (0,1))$ has
\[ \| u\|_{L^4(\Omega _T)\cap
C([0,T]; H^s (0,1))} \leq C\| \phi \|_{H^s(0,1)} \, ,\] where $C>0$
depends only on $s$ and $T$.
\end{proposition}

Next is discussed the IBVP of the associated linear problem with
nonhomogeneous Dirichlet boundary data, namely,
\begin{align}\label{boundary}
\left \{ \begin{array}{l} iu_t +u_{xx} =0, \qquad x\in (0,1), \ t\in
\R,
\\ \\ u(x,0) =0, \qquad u(0,t) = h_1(t), \quad u (1, t) = h_2 (t)
 ,\end{array}\right.
\end{align}
with the compatibility conditions $h_1 (0) = h_2(0) =0$ if necessary.
\begin{proposition}
The solution of (\ref{boundary}) can be expressed as
\begin{align}\label{explicit}
u( x, t) & = \sum^{+\infty}_{n=1} 2 i n\pi e^{-i (n \pi)^2t } \int^t _ 0
e^{i (n \pi)^2\tau } \Big (h_1 ( \tau ) - (-1)^n h_2 (\tau ) \Big ) d \tau
\sin n\pi x \nonumber \\
&=  W_{h} h_1 +  (W_{h} h_2)\Big|_{x\rightarrow 1-x}\, .
\end{align}
\end{proposition}
\noindent {\bf Proof:} Consider first the special case where $h_2\equiv
0$ and $h_1(0)=0$. Define $v$ by
\[ u(x,t)=v(x,t)+(1-x)h_1(t).\]
Then $v(x,t)$ solves
\begin{align*}
\left \{ \begin{array}{l} iv_t +v_{xx} =-i(1-x)h'_1 (t), \qquad x\in
(0,1), \ t\in \R,
\\ \\ v(x,0) =0, \qquad v(0,t) = 0, \quad v (1, t) = 0\, ,
 \end{array}\right.
\end{align*}
if $u$ solves (\ref{boundary}). As above, write $v(x,t)$ as
\[ v(x,t)=\sum ^{\infty}_{k=1} \alpha_k (t)\sin k\pi x .\]
Then, for $k=1,2,\cdots ,$
\[
\frac{d}{dt}\alpha _k  (t) + i (k\pi )^2 \alpha _k (t) = \beta_k
h'_1 (t), \qquad \alpha _k (0)=0,
\]
where
\[ \beta _k =- 2 \int ^1_0
(1-x)\sin k \pi x dx= -\frac{2}{k\pi}.\]
 It follows that
\[
\alpha _k (t)=\beta_k \int ^t_0 e^{-i (k\pi )^2 (t-\tau)} h'_1
(\tau ) d\tau = \beta _k h_1 (t) - i\beta _k (\pi k)^2 \int
^t_0e^{-i (k\pi )^2 (t-\tau)} h_1 (\tau ) d\tau \, .
\]
Substituting the latter into the original Fourier series representation yields
\[
v(x,t) =-(1-x) h_1(t) - \sum ^{\infty}_{k=1} i \beta _k (\pi k)^2 \int ^t_0e^{-i
(k\pi )^2 (t-\tau)} h_1 (\tau ) d\tau \, ,
\]
which in turn implies that
\[
u(x,t)= \sum ^{\infty}_{k=1} 2 i\pi k \int ^t_0e^{-i
(k\pi )^2 (t-\tau)} h_1 (\tau ) d\tau \sin k\pi x.
\]
Next,
consider the case of $h_1\equiv 0 $ and $h_2 (0)=0$. If we let
$x'=1-x$, this situation can be reduced to the case just studied.
Thus, if $h_1\equiv 0$ and $h_2 (0) =0$,
\[ u(x,t)=\sum ^{\infty}_{k=1} (-1)^{k+1} 2 i\pi k \int ^t_0e^{-i
(k\pi )^2 (t-\tau)} h_2 (\tau ) d\tau \sin k\pi x.
\]
The full representation \eqref{explicit} now follows. $\Box$

\begin{remark} One may view   the solution
$u$ in (\ref{explicit}) of (\ref{boundary}) as being written in
the form
\begin{equation} \label{special}
u(x,t)=\int ^t_0 W_0 (t-\tau )q(\cdot, \tau) d \tau
\end{equation}
where \[ q(x,t) =\Big (h_1
(t)-(-1)^nh_2(t)\Big ) \sum ^{\infty}_{n=1} 2 in\pi \sin n\pi x \, .
\]
Of course, $q$ belongs to the space $H^s \left (0,T; H^{-(3/2)-\epsilon} (0,1)\right )$  for any $\epsilon >0$ if $h_1, h_2 \in H^s (0,T)$. By semigroup theory, if $h_1, \ h_2 \in W^{1,1} (0,T)$, then
$u\in C\left ([0,T]; H^{(1/2)-\epsilon}(0,1)\right )$.
\end{remark}

Attention is now turned to the
boundary
integral
\begin{eqnarray}u_h =W_{h} h &=& \sum ^{\infty}_{n=1} 2i n\pi e^{-i(n\pi )^2
t}\int ^t_0 e^{i(n\pi)^2 \tau }h(\tau) d \tau \sin n\pi x \nonumber \\
&=& \sum ^{\infty}_{n=-\infty}  n\pi e^{-i(n\pi )^2
t}\int ^t_0 e^{i(n\pi)^2 \tau} h(\tau) d \tau e^{i n\pi x}. \label{4.9.1}
\end{eqnarray}
In the following, we will use the Lions-Magenes space  $ H_{00}^{1/2} (0 , T ) $ \cite{LM1972},
which is the interpolation space $[ H_0^1(0, T) , L^2(0, T) ]_\theta$ with $\theta = 1/2$.
\begin{proposition} \label{prop3.6}
For a given $T>0$, let $\Omega _T= (0,1)\times (0,T)$. If $h\in H^{\frac{1}{2}}_{00} (0,T)$, then
$$u_h=  W_h (\cdot )h \in L^4(\Omega _T)\cap C \left ([0,T]; L^2
(0,1)\right )$$ and there is a constant $C_T$ depending only on $T$ such that
\begin{equation}\label{estimate}
\| u_h\| _{L^4 (\Omega _T)} \leq C _T\| h\| _{H^{\frac{1}{2}}_{00} (0,
T)}
\end{equation}
and
\begin{equation}\label{estimate-}
\sup _{0\leq t\leq T} \| u_h(\cdot, t)\| _{L^2 (0,1)} \leq C_T \|
h\| _{H^{\frac{1}{2}}_{00} (0, T)}.
\end{equation}
\end{proposition}

\noindent{\bf Proof:} These results follow from analysis provided in Bourgain's paper \cite{bourgain-1}. In more detail, let $h(\tau) = \int^{\infty}_{-\infty}
e^{-\pi ^2 i \lambda \tau}\hat{h}(\lambda )d\lambda $. Write $u_h$ as follows:
\begin{eqnarray*}
u_h &=& \sum ^{\infty}_{n=-\infty} e^{-i(n\pi)^2t} e^{in\pi x}n\pi
\int ^{\infty}_{-\infty} \hat{h} (\lambda )\int ^t_0
e^{i(n\pi)^2\tau-\pi ^2 i\lambda \tau} d\tau d\lambda \\ \\
&=& \sum ^{\infty}_{n=-\infty} e^{-i(n\pi)^2t} e^{in\pi x}n\pi \int
^{\infty}_{-\infty} \hat{h} (\lambda )\frac{e^{i(n\pi)^2t-\pi ^2
i\lambda t}-1}{(n^2-\lambda )\pi ^2 i} d\lambda \\ \\
&=& \sum ^{\infty}_{n=-\infty} e^{-i(n\pi)^2t} e^{in\pi x}n\pi \left (
\int ^{0}_{ -\infty} +\int ^{\infty}_0 \right ) \hat{h} (\lambda
)\frac{e^{i(n\pi)^2t-\pi ^2 i\lambda t}-1}{(n^2-\lambda )\pi ^2 i}
d\lambda \\ \\
&=& I^-(x,t) +I^+ (x,t).
\end{eqnarray*}
Note that $I^+(x,t)$ also takes the form
\[
I^+(x,t)= \sum ^{\infty}_{n=1} 2 e^{-i(n\pi)^2t} n\pi\sin n\pi x
 \int ^{\infty}_0 \hat{h} (\lambda )\frac{e^{i(n\pi)^2t-\pi
^2 i\lambda t}-1}{(n^2-\lambda )\pi ^2 } d\lambda \, .
\]  

The quantity 
$I^+(x,t)$ is studied first.  Write
\begin{eqnarray*}
I^+(x,t)&=& \sum ^{\infty}_{n=-\infty} e^{-i(n\pi)^2t} e^{in\pi
x}n\pi \int ^{\infty}_{0} \hat{h} (\lambda )\psi (n^2-\lambda ) \sum
^{\infty}_{k=1} \frac{\left ((n^2-\lambda )t\pi ^2 i\right
)^k}{k! (n^2-\lambda )\pi ^2 i} d\lambda \\ \\
& & \qquad + \sum ^{\infty}_{n=-\infty} e^{in\pi x}n\pi \int
^{\infty}_{0} \hat{h} (\lambda )\Big (1-\psi (n^2-\lambda
)\Big )\frac{e^{-\lambda \pi ^2 i t}}{(n^2-\lambda )\pi ^2 i} d\lambda \\ \\
& & \qquad - \sum ^{\infty}_{n=-\infty} e^{-i(n\pi)^2t} e^{in\pi
x}n\pi \int ^{\infty}_{0} \hat{h} (\lambda )\Big (1-\psi (n^2-\lambda
)\Big )\frac{1}{(n^2-\lambda )\pi^2 i} d\lambda\\ \\
&=& I^+_1 +I^+_2 +I^+_3,
\end{eqnarray*}
where $\psi$ is a suitable $C^\infty$ cut-off function (see \cite{bourgain-1}).
For $I^+_1$, consider the individual summand
\[ I_{1,k}^+ = \sum ^{\infty}_{n=-\infty} e^{-i(n\pi)^2t} e^{in\pi
x}n\pi \int ^{\infty}_{0} \hat{h} (\lambda )\psi (n^2-\lambda )
 ((n^2-\lambda  )^k  d\lambda \, , \]
for $k =1,2, \cdots$.
 By Proposition 2.1 in  \cite{bourgain-1},
 \begin{eqnarray*}
 \left \| I^+_{1,k} \right \| ^2_{L^4 (\Omega _T)\cap L^{\infty}(0,T; L^2(0,1))}&\leq & C \left ( \sum
 ^{\infty}_{n=-\infty} n^2 \left | \int ^{\infty}_0
 \hat{h}(\lambda)\psi (n^2-\lambda ) \Big (n^2-\lambda \Big )^k d\lambda \right
 |^2 \right )\\ \\ &\leq & CB^k \left ( \sum
 ^{\infty}_{n=-\infty} n^2 \left | \int _{|\lambda -n^2|\leq B}
 \hat{h}(\lambda ) d\lambda \right |^2\right ) \\ \\
 &\leq & CB^{k+1} \left ( \sum
 ^{\infty}_{n=-\infty} n^2 \left | \int _{|\lambda -n^2|\leq B}
| \hat{h}(\lambda )|^2  d\lambda \right |\right )
\\ \\
 &\leq & CB^{k+1} \left ( \sum
 ^{\infty}_{n=-\infty}  \left | \int _{|\lambda -n^2|\leq B}
|\lambda | | \hat{h}(\lambda )|^2  d\lambda \right |\right ) \\ \\
& \leq & C B^{k+1} \int ^{\infty }_0 |\lambda | | \hat{h}(\lambda
)|^2 d\lambda \leq  C B^{k+1} \| h\| ^2_{H^{\frac12} (\R)} .
\end{eqnarray*}
Bounds on $I_1^+$ follow. Rewrite $I^+_2$ as
\begin{eqnarray*}
I^+_{2}(x,t)&=& \sum ^{\infty}_{n=0} 2 \sin n\pi x  \int ^{\infty }_0
 \hat{h} (\lambda )n\pi \Big (1-\psi (n^2-\lambda
)\Big )\frac{e^{-\lambda \pi ^2 i t}}{(n^2-\lambda )\pi ^2 } d\lambda \\
\\ &=& \sum ^{\infty }_{n=1} \frac{1}{\pi } \int ^{\infty }_{0} \hat{h}(\lambda )
 e^{-\lambda \pi ^2 i t}\Big (1-\psi (n^2-\lambda
)\Big ) \left (\frac{1}{n-\sqrt{\lambda}} +\frac{1}{n+\sqrt{\lambda }}
\right )\sin n \pi x   \,  d\lambda \\ \\ & =& \sum ^{\infty }_{n=1}
\frac{1}{\pi } \int ^{\infty }_{0} 2\mu\hat{h}(\mu ^2 ) e^{-\mu^2
\pi ^2 i t}\Big (1-\psi (n^2-\mu ^2 )\Big ) \left (\frac{1}{n-\mu}
+\frac{1}{n+\mu} \right )\sin n \pi x   \, d\mu .
\end{eqnarray*}
Applying Lemma A-1 in the Appendix leads to
 \begin{eqnarray*}
 \sup_{0\leq t\leq T}\| I^+_2(\cdot ,t)\|_{L^2 (0,1)} &\leq & C \sum ^{\infty }_{n=1}
 \left |\int ^{\infty}_0
  \mu  \hat{h}(\mu ^2) \Big (1-\psi (n^2-\mu^2)\Big ) \left ( \frac{1}{n-\mu} +\frac{1}{n+\mu}
\right )  d\mu \right |^2 \\ \\ & \leq & C\left \| \left ( 1 + |\mu|\right ) ^{3/2} \hat{h}(\mu
^2)\right \|_{L^2(\R)} \leq   C\|h\|^2_{H^{\frac12} (\R^+)}.
\end{eqnarray*}
To estimate the $L^4(\Omega_T)$-norm,  rewrite $I^+_{2}(x,t)$  as
\begin{eqnarray*}
I^+_{2}(x,t)&=& \sum ^{\infty}_{n=-\infty} e^{in\pi x}n\pi \left (
\int ^{\frac{n^2}{2}}_{0} + \int ^{\infty }_{\frac{n^2}{2}} \right )
 \hat{h} (\lambda )\Big (1-\psi (n^2-\lambda
)\Big )\frac{e^{-\lambda \pi ^2 i t}}{(n^2-\lambda )\pi ^2 i} d\lambda \\
\\ &:=& I^+_{2,1} +I^+_{2,2}\, .
\end{eqnarray*}
Proposition 2.6 in  \cite{bourgain-1} implies
\begin{align*}
\| I^+_{2,2} \|^2_{L^4 (\Omega _T)} \leq & \, C  \left ( \sum
^{\infty}_{n=-\infty } \int ^{\infty}_0 \frac{n^2\pi ^2
|\hat{h}(\lambda )|^2}{\big (|\lambda -n^2|+1\big )^2} \Big ( |\lambda
-n^2|+1\Big  )^{\frac34} \chi _{[\frac{n^2}{2} , \infty )} (\lambda
) \Big (1-\psi (n^2-\lambda )\Big )^2 d\lambda \right ) \\ \\
\leq & \, C \int ^{\infty}_0 |\lambda | |\hat{h} (\lambda )|^2 \sum
^{\infty }_{n=-\infty } \frac{1}{(1+|\lambda -n^2|)^{\frac54}}
d\lambda \leq  C \int ^{\infty}_0 |\lambda | |\hat{h} (\lambda )|^2
d\lambda \leq  C \| h\| ^2_{H^{\frac12} (R)}.
\end{align*}
Rewrite $I^+_{2,1}$ as
\begin{eqnarray*}
\left | I^+_{2,1} \right | &=& \left | 2 \int ^{\infty }_0 \left ( \sum ^{\infty }_{n=1}
 \chi _{[0, \frac{n^2}{2}]} (\lambda ) (1- \psi
(n^2-\lambda )) \frac{n\pi
\sin n\pi x }{(n^2-\lambda )\pi ^2}\right ) e^{- \lambda
\pi ^2 it} \hat{h} (\lambda )\,  d\lambda \right |  \\ \\ &=& \left | \frac{1}{\pi }
\int ^{\infty }_0 e^{- \lambda \pi ^2 it} \hat{h} (\lambda ) \left (
\sum ^{\infty }_{n=[\sqrt{2\lambda}]} \sin n\pi x \chi _{[0,
\frac{n^2}{2}]} (\lambda ) \left (\frac{1}{n-\sqrt{\lambda}}
+\frac{1}{n+\sqrt{\lambda }} \right ) \right ) d\lambda \right | \\ \\
&\leq& \frac{1}{\pi } \int ^{\infty }_0 | \hat{h} (\lambda )| \left
| \sum ^{\infty }_{n=[\sqrt{2\lambda}]} \left
(\frac{1}{n-\sqrt{\lambda}} +\frac{1}{n+\sqrt{\lambda }} \right
)\sin n \pi x  \right | d\lambda \, .
\end{eqnarray*}
To estimate the last sum, let
\[ S_k= \sum ^k_{n=1} \sin n\pi x
= \frac {\sin ((k+1)\pi x /2)\sin (k\pi x /2)}{\sin (\pi x/2)} \qquad (S_0 =0)\, .\]
For any $\alpha \in [0,1]$ and $ 0  < x \leq 1$, $
|S_k | \leq  C k^\alpha /|x|^{1-\alpha}$.
Consequently,
\begin{eqnarray*}
& & \sum ^{k }_{n=\lfloor\sqrt{2\lambda }\rfloor}
\frac{1}{n-\sqrt{\lambda}}(S_n-S_{n-1})= \sum ^{k
}_{n=\lfloor\sqrt{2\lambda }\rfloor}\frac{1}{n-\sqrt{\lambda}}S_n-\sum ^{k
}_{n=\lfloor\sqrt{2\lambda }\rfloor}\frac{1}{n-\sqrt{\lambda}} S_{n-1} \\ \\
& & \qquad \quad = \sum ^{k -1}_{n=\lfloor\sqrt{2\lambda }\rfloor} \left (
\frac{1}{n-\sqrt{\lambda }}-\frac{1}{n+1-\sqrt{\lambda}} \right )
S_n + \frac{1}{k-\sqrt{\lambda }} S_k
-\frac{1}{\lfloor\sqrt{2\lambda}\rfloor-\sqrt{\lambda}}S_{\lfloor\sqrt{2\lambda }\rfloor-1}\,
.
\end{eqnarray*}
Choose $ 3/4 < \alpha < 1$ and let $k\to \infty$ to come to the inequality
\begin{eqnarray*}
\left |\sum ^{\infty
}_{n=\lfloor \sqrt{2\lambda}\rfloor}\frac{1}{n-\sqrt{\lambda}} \sin n\pi x\right
|&\leq & C|x|^{\alpha -1}  \left ( \left (\sum ^{\infty
}_{n=\lfloor\sqrt{2\lambda}\rfloor}\frac{n^\alpha }{(n-\sqrt{\lambda})^2}\right
)+\frac{\lambda^{\alpha /2} }{\sqrt{\lambda}+1}\right ) \\ \\
&\leq & C |x|^{\alpha -1} \left (\frac{\lambda^{\alpha /2}}{\sqrt{\lambda}+1}+\sum ^{\infty}_{n=1}
\frac{1}{(n+\sqrt{\lambda})^{2-\alpha} } \right )\\ \\
&\leq & \frac{C}{|x| ^{1-\alpha} (1+\sqrt{\lambda})^{1 -\alpha} } \, .
\end{eqnarray*}
Using a similar argument for the other term gives
\begin{eqnarray*}
\left |I^+_{2,1} (x,t)\right | &\leq &C |x|^{\alpha -1} \int ^{\infty}_{0}
\frac{|\hat{h}(\lambda)|}{(1+\sqrt{\lambda})^{1-\alpha}} d\lambda \\ \\
&\leq & C |x|^{\alpha -1} \int ^{\infty}_0 (1+|\lambda |)^{\tilde \alpha } |\hat{h}(\lambda
)|\frac{d\lambda}{(1+\sqrt{\lambda})^{1-\alpha} (1+|\lambda|)^{\tilde \alpha}} \qquad
\Big(1/2 \geq \tilde \alpha >\alpha /2 \Big) \\ \\ &\leq &C |x|^{\alpha -1} \left (\int ^{\infty}_0 ( 1+|\lambda
|)^{2\tilde \alpha} |\hat{h}(\lambda )|^2 d\lambda \right )^{\frac12}
\left (\int ^{\infty}_0
\frac{d\lambda}{(1+\sqrt{\lambda})^{2-2\alpha} (1+|\lambda|)^{2\tilde \alpha}}\right
)^{\frac12} \\ \\&\leq &C |x|^{\alpha -1} \| h\|_{H^{\tilde \alpha}(\R ^+)} .
\end{eqnarray*}
Combining the foregoing result leads to the desired bound,
\[ \| I^+_2 \| ^2_{L^4 } \leq C \|h\|^2_{H^{\frac12}(\R ^+)} .\]
To study $I^+_3 (x,t)$, use again  Proposition 2.1 in
\cite{bourgain-1} to write
\begin{eqnarray*}
\|I^+_3 \| _{L^4 (\Omega _T)\cap L^{\infty} (0,T; L^2 (0,1))}  &\leq
& C \left (\sum ^{\infty}_{n=1 } n^2 \left | \int ^{\infty}_0
\hat{h} (\lambda )\frac{1-\psi (n^2-\lambda
)}{\lambda -n^2} d \lambda \right |^2 \right ) \\ \\
&=& C \sum ^{\infty}_{n=1} \left | \int ^{\infty }_0 \hat{h}
(\lambda )\left ( \frac{1}{\sqrt{\lambda }-n}-\frac{1}{\sqrt{\lambda
}+n}\right ) \left ( 1-\psi (n^2-\lambda )) d\lambda \right |^2 \right ) \\
\\
&\leq & C \left ( \sum ^{\infty }_{n=1} \left | \int ^{\infty}_0 \mu
\hat{h}(\mu ^2) \frac{1}{\mu -n}(1-\psi (n^2-\mu ^2))d\mu \right |^2
\right. \\ \\
 & + & \left. \sum ^{\infty }_{n=1} \left | \int ^{\infty}_0 \mu
\hat{h}(\mu ^2)
\frac{1}{\mu +n}(1-\psi (n^2-\mu ^2))d\mu \right |^2 \right ) \\ \\
&\leq & C \| \mu \hat{h}(\mu ^2)\|^2_{L^2 }+ \int ^{\infty }_0 \left
( \int ^{\infty }_0 \frac{|\mu \hat{h} (\mu ^2)|}{\mu +y}d\mu \right
)^2 dy \\ \\
&\leq & C \| \mu \hat{h}(\mu ^2) \|^2_{L^2} \leq C \|
h\|^2_{H^{\frac14}(\R ^+)} .
\end{eqnarray*}
In summary, it appears that \[\|
I^+\|^2_{L^4} \leq C \| h\| ^2_{H^{\frac12} (\R ^+)} .\]

Now consider $I^-(x,t)$ and  express it in the form
 \begin{eqnarray*}
 I^-(x,t)&=& \sum ^{\infty }_{n=-\infty } e^{-i(n\pi )^2 t} e^{in
 \pi x} n \pi \int ^{\infty }_0 \frac{e^{i(n\pi )^2 t+i\lambda \pi
 ^2 t}-1}{(\lambda +n^2)\pi^2 i}\hat{h}(-\lambda ) d\lambda  \\ \\
 &=& \sum ^{\infty }_{n=-\infty } e^{in\pi x} n\pi \int ^{\infty }_0
 \frac{e^{i\lambda \pi ^2 t}-e^{-i(n\pi )^2 t}}{(\lambda +n^2)\pi ^2
 i}\hat{h}(-\lambda ) d\lambda  \\ \\
 &:=& I^-_1 - I^-_2 .
 \end{eqnarray*}
For $I^-_2$, Proposition 2.1 of \cite{bourgain-1} implies
\begin{eqnarray*}
&&\| I^-_2\|^2_{L^4(\Omega _T )} + \| I^-_2\|^2_{ L^{\infty}(0,T;
L^2 (0,1))} \leq C \sum ^{\infty }_{n=-\infty} n^2 \pi ^2 \left |
\int ^{\infty }_0 \frac{\hat{h}(-\lambda )}{\lambda +n^2} d
\lambda \right |^2 \\ \\
 &\leq & C \sum ^{\infty}_{n=1} \left | \int ^{\infty }_0
\frac{\hat{h}(-\lambda )n}{\lambda +n^2} d \lambda \right |^2 \leq C
\int ^{\infty }_{1} \left |\int ^{\infty }_0 \frac{\hat{h}(-\lambda
)y}{\lambda +y^2} d \lambda \right |^2 dy
  \\ \\
&\leq & C \int ^{\infty }_1 \left ( \int ^{\infty }_0 (1+ |\lambda |
^{\frac12} )^2 |\hat{h} (-\lambda )|^2 \, d\lambda \int ^{\infty }_0
\frac{y^2}{(y^2+\lambda )^2 (1+|\lambda |^{\frac12})^2} \, d\lambda \right ) dy
\leq  C \| h\|^2_{H^{\frac12}(\R ^+)} .
\end{eqnarray*}
The formula
\[ \sum ^{\infty }_{n=1} \frac{n\sin nx}{n^2 + a^2}= \frac{\pi}{2}
\frac{\sinh a (\pi -x)}{\sinh a \pi}, \qquad \mbox{for}\quad 0< x< 2\pi, \]
 which holds for all $a$, allows us to write
\begin{eqnarray*}
I^-_1 &=& \sum ^{\infty }_{n=1} 2n\pi \sin n\pi x\int ^{\infty }_0
\hat{h} (-\lambda )e^{i\lambda \pi ^2 t} \frac{1}{(n^2+\lambda )\pi
^2 } d\lambda \\ \\
&=& \int ^{\infty }_0 \frac{2\hat{h} (-\lambda ) e^{i\lambda \pi
^2t}}{\pi } \sum ^{\infty }_{n=1} \frac{n\sin \pi n x}{n^2+ \lambda
} d\lambda = \int ^{\infty}_0 \hat{h}(-\lambda ) e^{i\lambda \pi ^2
t} \frac{\sinh \sqrt{\lambda } (\pi -x)}{\sinh \sqrt{\lambda} \pi }
d\lambda\, .
\end{eqnarray*}
Consequently, it is seen that
\begin{eqnarray*}
 \left | I^-_1 (x,t) \right | &\leq & C
\int ^{\infty}_0 |\hat{h} (-\lambda )| e^{-\sqrt{\lambda} \pi x}
d\lambda \\ \\
&\leq & \left | \int ^{\infty }_0 |\hat{h}(-\lambda )|^2 (1+|\lambda
|)d \lambda \right |^{\frac12} \left | \int ^{\infty }_0
\frac{e^{-2\sqrt{\lambda}\pi x}}{1+|\lambda |}d \lambda \right
|^{\frac12}\!\!, \end{eqnarray*} which implies
\begin{align*}
 \sup_{0\leq t\leq T}\int ^1_0\left | I^-_1 (x,t) \right |^2 dx
& \leq
 C \int ^1_0 \left ( \int ^{\infty }_0 |\hat{h}(-\lambda )|^2 (1+\lambda )\, d \lambda
\int ^{\infty }_0 \frac{e^{-2\sqrt{\lambda}\pi x}}{1+\lambda }\, d
\lambda \right ) dx \\ \\ &  \leq C \| h\|^2_{H^{\frac12} (\R ^+)}
\int ^{\infty }_0  \int ^1_0 \frac{e^{-2\sqrt{\lambda}\pi
x}}{(1+\lambda )} dx
 d\lambda  \leq  C\|h\|^2_{H^{\frac12} (\R ^+)}
\end{align*}
and
\begin{eqnarray*}
\int ^T_0 \int ^1_0 |I^-_1 (x,t) |^4 dx dt &\leq & C \int ^1_0 \left
| \int ^{\infty }_0 |\hat{h} (-\lambda)|e^{-\sqrt{\lambda} \pi x}
d\lambda \right | ^4 dx
\\ \\
&\leq & C \int ^1_0 \left | \int ^{\infty }_0 |\hat{h}(-\lambda )|^2
(1+|\lambda |)d \lambda \right |^2 \left | \int ^{\infty }_0
\frac{e^{-2\sqrt{\lambda}\pi x}}{1+|\lambda |}d \lambda \right |^2
dx
\\ \\
&=& C \| h\|^4_{H^{\frac12} (\R ^+)} \left ( \int ^{\infty }_0 \left
( \int ^1_0 \frac{e^{-4\sqrt{\lambda}\pi x}}{(1+|\lambda |)^2} dx
\right )^{\frac12} d\lambda \right )^2 \\ \\
&\leq & C\|h\|^4_{H^{\frac12} (\R ^+)} \left ( \int ^{\infty }_0
\frac{1}{1+|\lambda |}\frac{d\lambda }{\lambda ^{\frac14}} \right
)^2 \leq  C\| h\|^4_{H^{\frac12} (\R ^+)}\, .
\end{eqnarray*}
Hence, we arrive at
$ \|u_h \|_{L^4 (\Omega _T)\cap L^{\infty} (0,T; L^2 (0,1))} \leq C\|h\| _{H^{\frac12}(\R ^+)}$
and the proof is complete. $\Box$

\medskip
If the regularity of $h(t)$ is higher, $W_h (t)
h$ is smoother.
\begin{proposition} \label{prop3.7-1}
Let $s\geq 0$ be given. For any $h\in H^{\frac{1+s }{2}}_0
(0,T)$ (here for $ s$ an even integer, $h$ should be in $ H^{\frac{1+s }{2}}_{00}
(0,T)$), let $u= W_h h$. Then, $\partial ^s_x u$  belongs to $L^4((0,1)\times [0, T])\cap C([0,T]; L^2(0,1))$ and satisfies
\begin{equation*}
\left \|  \partial_x^su  \right \| _{L^4 (\Omega _T)} \leq C \| h\|
_{H^{\frac{1+s }{2}} (0, T)}
\end{equation*}
and
\begin{equation*}
\sup _{0\leq t\leq T}\left \| \ \partial ^s_x u \right \| _{L^2
((0,1)} \leq C \| h\| _{H^{\frac{1+s }{2}} (0, T)}
\end{equation*}
where $C>0$ is a constant independent of $h$.
\end{proposition}

\noindent{\it Proof:} We only need to prove it for $s=2$. The cases
where $s \in (0, 2)$ can then be obtained by interpolation, where we note that $H^s_0 ( 0, T )$ is an interpolation space for $ s \not = \mbox{integer } + 1/2$ while for $ s  = \mbox{integer } + 1/2$, the corresponding interpolation space is the Lions-Magenes space $H^s_{00} ( 0, T )$ \cite{LM1972}. The proof
for $s> 2$ is same as for $s=2$.

Notice that the $t$-derivative of $  W_h (\cdot ) h$
 satisfies the system (\ref{boundary}) with boundary
condition $h'(t)$ and zero initial condition. Hence, by Proposition \ref{prop3.6},
there obtains
\begin{equation*}
\left \| \frac{\partial  W_h (\cdot ) h }{\partial t} \right \|
_{L^4 ((0,1)\times [0, T])} \leq C \| h' (t) \| _{H^{\frac{1 }{2}} (0,
T)}\leq C \| h (t) \| _{H^{\frac{3 }{2}} (0,
T)}\, .
\end{equation*}
But, bounds on
one $t$-derivative of $W_h (t) h$ give bounds on two $x$-derivatives of $W_h
(t) h$. Thus, the case for $s=2$ is established. $\Box$

\medskip
\begin{remark}
Notice that
\[ \| W_0 (t) \phi \|_{L^4 ((0,1)\times (0,T))} \leq C \| \phi \|
_{L^2 (0,1)}\] for any $\phi \in L^2 (0,1)$ and, in addition, for
the linear Schr\"odinger equation posed on the  half-line,
\[ \| W_b (\cdot) h\|_{L^q (\R^+; L^r (\R^+))} \leq C\|h\|
_{H^{\frac14} (\R^+)} \] for any $h\in H^{\frac14} (\R^+)$, where $ (q, r) $ is an admissible pair satisfying $ \frac1{q} + \frac1{2r} = \frac14$.  One
thus wonders whether the estimate  (\ref{estimate}) or (\ref{estimate-}) can be improved.
Example A-2 in the Appendix  shows that if $\| W_b (\cdot )h \|_{L^2([0,1]\times [0, T])} \leq C \| h\|_{H^s ([0, T])}$ for all $h(t) \in H^s([0, T])$, then it must be the case that $s\geq \frac12$. Thus, the estimates in (\ref{estimate}) and (\ref{estimate-}) are optimal.
 \end{remark}

\subsection{The nonlinear problem}

In this subsection,  the full nonlinear IBVP
\begin{align}\label{5.1}
\left \{ \begin{array}{l} iu_t +u_{xx} + \lambda u |u |^{p-2} =0,
\qquad x\in (0,1), \ t\in \R ^+,
\\ \\ u(x,0) =\phi (x) , \qquad u(0,t) = h_1(t), \quad u (1, t) = h_2
(t) \end{array}\right.
\end{align}
 with $\phi \in
H^s(0,1)$ and $h_1, \ h_2 \in H^{\frac{s+1}{2}}_{loc} (\R^+)$ is studied. 
A local well-posedness theorem
is formulated and proved.
\begin{theorem}
Let $3\leq p< \infty$,  $\frac12 < s < \frac52$ 
and $\lfloor s\rfloor < p-2$,  $T>0$  and $r>0$ be given.
There exists a $T^*
>0$ such that if $(\phi, h_1 , h_2 )\in {\cal
X}_{s,T}:=H^s(0,1)\times H^{\frac{s+1}{2}} (0,T)\times H^{\frac{s+1}{2}} (0,T)$ satisfies $h_1 (0) = \phi (0), h_2 (0) = \phi(1)$ and $ \|
(\phi , h_1, h_2 )\| _{{\cal X}_{s,T}} \leq r, $ the IBVP
(\ref{5.1}) admits a unique solution $ u\in C([0, T^*];
H^s(0,1)).$ Moreover, the solution $u$ depends on $(\phi, h_1,h_2)$
continuously in the corresponding spaces.
\end{theorem}
\noindent {\bf Proof:} We only consider the cases where $\frac12 < s\leq
2$. In addition, without loss of generality, we assume that $\phi (0)=h_1 (0)=0$ and $\phi
(1)=h_2 (0)=0$. For if not, we can homogenize the   boundary conditions by
writing $u = v + h_1 (0) ( 1- x) + h_2 (0) x = v + \gamma(x)$.   Then $v$ satisfies homogeneous
compatibility conditions and the equation
$$
           iv_t +v_{xx} + \lambda |v+ \gamma|^{p-2}( v+ \gamma) =0.
$$
 As $\gamma$ is smooth and the direct estimates made of the nonlinear term, {\it e.g.}
\eqref{algebra}, are very simple, theory for either $u$ or $v$ follows exactly the
same lines.

For $s > \frac12$, $H^s (0,1)$ is a Banach
algebra. It follows that  there is a constant $C= C_s$ such that
\begin{equation} \label{algebra}
 \big\| v|v|^{p-2}\big\|_{H^s (0,1)} \leq C \big\| v\big\|^{p-1}_{H^s
(0,1)},
\end{equation}
when $ s =1, 2$.   Indeed, for any $s$ with $\lfloor s\rfloor < p -2$,  the chain rule for
fractional derivatives implies the same  result.

For any $\theta $ with $0< \theta \leq T$ and $v\in C([0, \theta ]; H^s
(0,1))$,
Propositions \ref{pro3.1} and \ref{prop3.7-1} imply that the linear IBVP
\begin{align}\label{5.1-1}
\left \{ \begin{array}{l} iu_t +u_{xx} + \lambda v |v |^{p-2} =0,
\qquad x\in (0,1), \ t\in \mathbb{R},
\\ \\ u(x,0) =\phi (x) , \qquad u(0,t) = h_1(t), \quad u (1, t) = h_2
(t) \, ,\end{array}\right.
\end{align}
admits a unique solution $u\in C([0, \theta ]; H^s  (0,1))$.
Moreover, there exists a constant $C>0$ independent of $\theta $
such that
\[ \| u\|_{C([0, \theta ]; H^s  (0,1))} \leq C\|(\phi , h_1,
h_2)\|_{{\cal X}_{s,T}} + C\theta \| v\|^{p-1}_{C([0, \theta
]; H^s (0,1))} .\]
Thus, for any given $(\phi, h_1, h_2 )\in {\cal
X}_{s,T}$, the IBVP (\ref{5.1-1}) defines a nonlinear map $\Gamma$
from  $ Y_{s, \theta}:= \{ w\in C([0,
\theta ]; H^s (0,1)) \} $ to $Y_{s, \theta}$. A well understood argument, similar to the contraction mapping argument in Section 7 of \cite{holmer} using the chain rule, now reveals that if $\theta
>0$ is chosen small enough, there exists  an $M>0$ such that
\[ \| \Gamma (v_0) \| _{C([0, \theta ]; H^s  (0,1))} \leq M \]
and
\[ \| \Gamma (v_1)- \Gamma (v_2)\|_{C([0, \theta ]; H^s  (0,1))}
\leq \frac12 \| v_1-v_2\|_{C([0, \theta ]; H^s  (0,1))} \] for any
$v_0, \ v_1, \ v_2 \in C([0, \theta ]; H^s (0,1))$ with
\[ \| v_j\|_{C([0, \theta ]; H^s (0,1))}\leq M,\qquad j=0,1,2.\]
Hence, the map $\Gamma$ is a contraction  whose unique fixed point
is the desired solution $u$ of (\ref{5.1}). The proof is complete.
$\Box$

\medskip
Next, we aim to show the well-posedness of the IBVP (\ref{5.1}) in $H^s (0,1)$ for $0\leq s< \frac 12$. To this end, consider
 the integral equation
\begin{equation}\label{5.2}
u( \cdot , t) = W_0 ( t) \phi + W_{h} h_1 +(W_{h} h_2 )\Big |
_{x\rightarrow 1-x}+ i \lambda \int^t_0 W_0 ( t-\tau ) \Big ( u (\cdot ,
\tau ) | u ( \cdot , \tau ) |^{p-2} \Big ) d \tau \, ,
\end{equation}
associated with the IBVP (\ref{5.1}).
\begin{proposition}\label{pro3.2}
Let $0\leq s< \frac12$ and $T>0$. Suppose $r>0$ to be given and
$ 3\leq p \leq 4$.
There exists a $T^* = T^*(r) >0$
such that for any
\[ (\phi, h_1, h_2)\in {\cal X}_{s,T}\] with $\| (\phi , h_1,
h_2 )\|_{{\cal X}_{s,T}}\leq r$, (\ref{5.2}) admits a unique
solution \[ u\in {\cal Y}_{s, T^*}:= L^4((0,1)\times (0,T^*))\cap
C([0,T^*]; H^s (0,1))\] which depends continuously on $(\phi, h_1, h_2)$
 in the corresponding spaces.
\end{proposition}
{\bf Proof:} Solving (\ref{5.2}) can be  viewed as
a problem of finding a fixed point of a nonlinear operator.
Consequently,  the proposition follows
using the argument  that appears already in  \cite{bourgain-1} along with
our Proposition \ref{prop3.6} for the boundary integrals.
$\Box$

\medskip
The solution $u$ of (\ref{5.2}) given by Proposition \ref{pro3.2} is
a mild solution of the IBVP (\ref{5.1}).  By the same arguments as put forward already in the
proofs of Propositions \ref{prop3.3} and \ref{prop3.7}, it is deduced that the IBVP
(\ref{5.1}) admits at most one mild solution, thereby settling the validity of  the following
theorem.
\begin{theorem} Under the conditions in Proposition \ref{pro3.2}, the IBVP (\ref{5.1}) is
unconditionally  locally well-posed in $H^s (0,1)$ for
$0\leq s < \frac12$.
\end{theorem}

\section{Global Well-Posedness}
\setcounter{equation}{0}

In this section, consideration is given to the issue of global
well-posedness for both the problems,
\begin{align}\label{6.1}
\left \{ \begin{array}{l} iu_t +u_{xx} + \lambda u |u |^{p-2} =0,
\qquad x\in\mathbb{ R}^+, \ t\in \mathbb{R},
\\ \\ u(x,0) =\phi (x) , \qquad u(0,t) = h(t) \end{array}\right.
\end{align}
and
\begin{align}\label{6.2}
\left \{ \begin{array}{l} iu_t +u_{xx} + \lambda u |u |^{p-2} =0,
\qquad x\in (0,1), \ t\in \mathbb{R},
\\ \\ u(x,0) =\phi (x) , \qquad u(0,t) = h_1(t), \quad u (1, t) = h_2
(t) \end{array}\right.
\end{align} in  $H^s (\mathbb{R}^+)$ and $H^s(0,1)$,
respectively. Since the local well-posedness of both problems
has been established, global well-posedness will follow from
suitable \emph{a-priori} estimates.


\medskip
First, recall that  if $u (x,t)$ is a smooth solution of the
NLS equation
\[ iu_t +u_{xx} + \lambda u |u
|^{p-2} =0, \] then the following identities
\begin{align}
&\frac{\partial }{\partial  t}  ( | u|^2 )  = -2 \im (u_x ( x , t) \bar u ( x, t)
)_x  \, ,\label{4.1}\\ \nonumber \\ &\frac{\partial  }{\partial  t}\left (  |u_x |^2
- \frac{2 \lambda }{ p} | u |^p
 \right ) = 2 \re \left ( u_x ( x, t) \bar u_t ( x, t)
\right )_x  \, ,\label{4.2}\\ \nonumber \\ &\left ( |u_x ( x, t) |^2
+ \frac{ 2 \lambda  } { p } |u ( x, t)|^p \right ) _x = - i \left (
\frac{\partial  }{\partial  t}  ( u \bar u_x ) - ( u ( x , t) \bar u _t ( x, t))_x
\right )\, \label{4.3}
\end{align}
were obtained in \cite{bc}.
Multiply both sides of (\ref{4.3}) by a smooth, time-independent function $\eta
(x)$
and write
$$
u \bar u _t = u ( -i \bar u_{xx} - i \lambda \bar u | u |^{p-2} ) =
-i ( ( u \bar u _x) _x - u _x \bar u _x + \lambda | u| ^p )
$$
to derive the formula
\begin{align}
\eta (x) \left ( |u_x ( x, t) |^2 + \frac{ 2 \lambda  } { p } |u (
x, t)|^p \right ) _x & =- i \left ( \frac{\partial  }{\partial  t}  ( \eta (x) u
\bar
u_x )\right ) + i ( \eta u \bar u_t ) _x \nonumber \\
& - ( \eta_x u \bar u _x ) _x + \eta_{xx} ( u \bar u_x) + \eta_x |
u_x |^2 - \lambda \eta_x | u | ^p\, .\label{4.4}
\end{align}
By choosing
appropriate functions $\eta (x)$, one can obtain various pointwise estimates of
$u(x,t)$.
  In particular, for any given interval
$[a, b]$, choose $\eta (x) \in C^\infty ( \R ) $ such that $\eta = 1
$ for $ x\leq a $ and $\eta = 0 $ for $ x \geq b$ with $|\eta (x) | \leq
1$ for all $x$. Integrating (\ref{4.4}) from $a$ to $b$ with respect to $x$
 and integrating by parts yields
\begin{align*}
& -  |u_x ( a, t) |^2 -  \frac{ 2 \lambda  } { p } |u ( a, t)|^p -
\int^b_a\eta_x  (x) \left ( |u_x ( x, t) |^2 + \frac{ 2 \lambda  } {
p } |u ( x, t)|^p \right )dx \\
& = - i  \frac{\partial  }{\partial  t}  \left ( \int^b_a \eta (x) u \bar u_x dx \right )
-  i \Big ( \eta (a ) u(a , t)  \bar u_t(a , t) \Big )  + \int^b_a \Big (
\eta_{xx} ( u \bar u_x) + \eta_x | u_x |^2 - \lambda \eta_x | u | ^p
\Big ) dx\, .
\end{align*}
If $v = u_t$, then $v$ satisfies the
equation
\begin{equation}
 iv_t + v _{xx} + (\lambda p/ 2)  |u
|^{p-2} v + (\lambda (p-2)/ 2 ) |u|^{p-4} u^2 \bar v = 0\,
,\label{ut-eq}
\end{equation}
which is linear in terms of $v$. Similar identities as
(\ref{4.1})-(\ref{4.4}) hold for $v$;
\begin{align}
&\frac{\partial }{\partial  t}  ( | v|^2 )  = -2 \im \big (v_x ( x , t) \bar v ( x, t)
\big )_x   - \lambda(p-2)|u|^{p-4} \im \big (u^2\bar v^2 \big )
\, ,\label{4.1.11}\\ \nonumber \\
&\frac{\partial  }{\partial  t}  |v_x |^2  = 2 \re \left ( v_x ( x, t) \bar v_t ( x,
t) \right )_x  + \im \big ( \lambda p |u |^{p-2} v \bar v_{xx}\big  ) + \im \big (\lambda
(p-2) |u |^{p-4} u^2 \bar v \bar v_{xx} \big )
\, ,\label{4.2.11}\\ \nonumber \\
&\left ( |v_x ( x, t) |^2 \right ) _x = - i \left ( \frac{\partial  }{\partial  t} (
v \bar v_x ) - ( v ( x , t) \bar v _t ( x, t))_x \right ) \nonumber
 \\
& \qquad \qquad \qquad\qquad - \lambda (p/2) |u |^{p-2} (|v|^2 ) _x
- \lambda (p-2) | u | ^{p-4} \re \big  (u ^2 \bar v \bar v_x \big )\,
.\label{4.3.11}
\end{align}

\medskip
These identities will play a role in our study of global well-posedness.  The quarter-plane
IBVP (\ref{6.1}) will be considered next while the IBVP (\ref{6.2}) will
be dealt with in Subsection 5.2.

\subsection{Global well-posedness on $\mathbb{R}^+$}

\begin{proposition}\label{pro-7.1} Assume that $p\geq 2$ if $\lambda
<0$ and $2\leq p\leq 4$ if $\lambda >0$.  Let $T>0$ be given. Then
there exists a nondecreasing, continuous function $\alpha :
\mathbb{R}^+\to \mathbb{R}^+$ with $\alpha (0)=0$ such that any
smooth solution $u$ of (\ref{6.1}) satisfies
\begin{equation}  \label{5.1inequality}
 \sup _{0\leq t\leq T} \| u(\cdot, t)\|_{H^1(R^+)} \leq \alpha \Big (\| \phi \|_{H^1 (R^+)}+ \|h\|_{H^1(0,T)}\Big ).
\end{equation}
Here $\alpha$ also depends upon $T$ and other constants and is bounded for any $T > 0$.
 \end{proposition}

\begin{remark}
The calculations to follow can easily be justified for solutions that are in $H^2(\R^+)$
in space with boundary traces that are continuous functions of time.   Note that this
result does not depend upon how the solution is obtained, but simply asserts {\it a priori}
information that it must obey.
\end{remark}

 \noindent
 {\bf Proof:}
  First,
integrate (\ref{4.4}) with  $\eta = 1$ from $0$ to $t$ to obtain
\begin{align}
\int^t _0 |u_x ( 0, s) |^2  ds  =  & i \left ( \int^\infty_0 u \bar
u_x dx  \right )\bigg |^t_0 - \frac{ 2 \lambda  } { p }\int^t _0 |u
( 0, s)|^pds+ i \int^t _0u(0 , s)  \bar u_t(0 , s) ds \nonumber \\
 =&i \left ( \int^\infty_0 u ( x , t)\bar u_x( x , t)  dx  \right )
-i \left ( \int^\infty_0 u ( x , 0)\bar u_x( x , 0)  dx  \right ) +
C_1(t)\nonumber \\
= & i \left ( \int^\infty_0 u ( x , t)\bar u_x( x , t)  dx  \right )
+ c_1 + C_1(t)\nonumber \\
\leq & \left ( \int^\infty_0 |u ( x , t)|^2  dx \right )^{1/2} \left
( \int^\infty_0 | u_x( x , t)|^2  dx  \right )^{1/2}+ c_1 + C_1(t)\, ,\label{1000}
\end{align}
where $c_1$ is dependent on the initial data and $C_1(t) $ is
dependent on the boundary data with $C_1 (0 ) = 0$.  It follows from
(\ref{4.1}) that
\begin{align*}
\int^\infty_0 |u (x, t) |^2 dx & \, =   \int^\infty_0 |u (x, 0) |^2 dx +
2\im \int^t _0 (u_x ( 0 , s) \bar u ( 0, s) )ds \\
\leq & \, c_1 +2 \left ( \int^t _0 |u_x ( 0 , s)|^2 ds  \int^t _0 | u (
0, s) |^2 ds\right )^{1/2}\\
= & \, c_1 +2C_1 (t)  \left ( \int^t _0 |u_x ( 0 , s)|^2 ds  \right
)^{1/2} \\
\leq & \, c_1 +2C_1 (t) \left ( \left ( \int^\infty_0 |u ( x , t)|^2 dx
\right )^{1/2} \left ( \int^\infty_0 | u_x( x , t)|^2  dx  \right
)^{1/2}+ c_1 + C_1(t)  \right )^{1/2}\\
\leq & \, c_1 +2C_1 (t) \left ( \left ( \int^\infty_0 |u ( x , t)|^2 dx
\right )^{1/4} \left ( \int^\infty_0 | u_x( x , t)|^2  dx  \right
)^{1/4}+ (c_1 + C_1(t) )^{1/2}  \right )\\
\leq & \, c_1 +2C_1 (t) (c_1 + C_1(t) )^{1/2} + \frac14 \int^\infty_0 |u
( x , t)|^2 dx \\
&  \qquad + \frac34 \left ( 2C_1 (t) \left ( \int^\infty_0 | u_x( x ,
t)|^2  dx \right )^{1/4}\right )^{4/3}\!\! .
\end{align*}
A direct consequence is the inequality
\begin{align}
\int^\infty_0 |u (x, t) |^2 dx \leq & \, \frac43 \big( c_1 +2C_1 (t) \big(c_1 +
C_1(t) )^{1/2}\big)  +  \left ( 2C_1 (t) \right )^{4/3} \left (
\int^\infty_0 | u_x( x , t)|^2  dx \right )^{1/3}\nonumber \\
= & \, D_1 + \left ( 2C_1 (t) \right )^{4/3} \left ( \int^\infty_0 |
u_x( x , t)|^2  dx \right )^{1/3} \label{5.10.1}
\end{align}
where $D_1 $ is  a constant depending on both the initial  and
boundary data.

If $\lambda <0$, integrating both sides of (\ref{4.2}) in $x$ over
$\mathbb{R}^+$ and $t$ over $[0,t]$ yields
\begin{align*}
 \int^\infty_0 &|u _x(x, t)  |^2 dx =   \frac{2 \lambda }{ p} \int^\infty_0
 | u (x,t) |^p dx + \int^\infty_0 \left ( |u _x (x, 0) |^2 -\frac{2 \lambda }{
 p} | u (x,0) |^p \right ) dx \\
 & \qquad \qquad \qquad \qquad -
2 \re \int^t_0  u_x ( 0, s) \bar u_s ( 0, s) ds\, .
\end{align*}
The right-hand side of this equation may be bounded thusly (note that the first term is negative and the second term only depends on initial data):
\begin{align*}
\mbox{rhs}\leq &  \,    c_1 +  \int^t_0  |u_x ( 0, s)|^2 ds +
\int^t_0| u_s ( 0, s)|^2 ds= D_1 +  \int^t_0  |u_x ( 0, s)|^2 ds\\
\leq &  \, D_1 +   \left ( \int^\infty_0 |u ( x , t)|^2  dx \right
)^{1/2} \left ( \int^\infty_0 | u_x( x , t)|^2  dx  \right
)^{1/2}+ c_1 + C_1(t) \\
\leq & \, D_1 +   c_1 + C_1(t)  + \left ( D_1 + \big ( 2C_1 (t) \big
)^{4/3} \left ( \int^\infty_0 | u_x( x , t)|^2  dx \right )^{1/3}
\right )^{1/2} \left ( \int^\infty_0 | u_x( x , t)|^2  dx  \right
)^{1/2}\\
\leq & \, D_1 +  c_1 + C_1(t) +    \left ( D_1^{1/2} + \big ( 2C_1 (t)
\big )^{2/3} \left ( \int^\infty_0 | u_x( x , t)|^2  dx \right
)^{1/6} \right ) \left ( \int^\infty_0 | u_x( x , t)|^2  dx \right
)^{1/2}\\
\leq &  \, D_1 +   c_1 + C_1(t)  + D_1^{1/2}\| u _x ( \cdot , t)
\|_{L^2} +\big( 2C_1 (t) \big )^{2/3}\| u _x ( \cdot , t)
\|_{L^2}^{4/3},
\end{align*}
where $c_1$ again depends only on the initial data and $D_1$ depends on initial and boundary data.
Here, $2 ab \leq a^ 2 + b^2$ is used for the first inequality,  \eqref{1000} is applied for the second inequality,
the third inequality is from \eqref{5.10.1}, and the fourth inequality uses the inequality $(a+ b )^{1/2} \leq a^{1/2} + b^{1/2}$.
Hence,  over any finite time interval, $\| u _x ( \cdot , t)
\|_{L^2(\mathbb{R}^+)}$ is uniformly bounded.  Appealing to (\ref{5.10.1}) again
reveals that $\|
u ( \cdot , t) \|_{L^2(\mathbb{R}^+)}$ is also bounded for any bounded time interval.

\medskip
If $\lambda >0$,
equation (\ref{4.2}) implies
\begin{align*}
  \int^\infty_0&  |u _x(x, t)  |^2 dx =   \frac{2 \lambda }{ p} \int^\infty_0
 | u (x,t) |^p dx + \int^\infty_0 \left ( |u _x (x, 0) |^2 -\frac{2 \lambda }{
 p} | u (x,0) |^p \right ) dx \\
 & \qquad \qquad \qquad \qquad -
2 \re \int^t_0  u_x ( 0, s) \bar u_s ( 0, s) ds\\
\leq & \, \frac{4 \lambda }{ p} \Big ( \| u ( \cdot , t) \|_{L^2} \| u
_x ( \cdot , t) \|_{L^2}\Big )^{ (p -2) /2}  \int^\infty_0
 | u (x,t) |^2 dx  + c_2 +  \int^t_0  |u_x ( 0, s)|^2 ds +
\int^t_0| u_s ( 0, s)|^2 ds\\
= & \,\frac{4 \lambda }{ p} \Big ( \| u ( \cdot , t) \|_{L^2} \| u _x
( \cdot , t) \|_{L^2}\Big )^{ (p -2) /2}  \int^\infty_0
 | u (x,t) |^2 dx  + \tilde D +  \int^t_0  |u_x ( 0, s)|^2 ds\\
\leq & \,\tilde D + \frac{4 \lambda }{ p}\| u _x ( \cdot , t)
\|_{L^2}^{ (p -2) /2}\left ( D_1 + \big ( 2C_1 (t) \big )^{4/3}
\left ( \int^\infty_0 | u_x( x , t)|^2  dx \right )^{1/3}\right
)^{(p+2)/4}\\
& \qquad \qquad+  \left ( \int^\infty_0 |u ( x , t)|^2  dx \right )^{1/2}
\left ( \int^\infty_0 | u_x( x , t)|^2  dx  \right )^{1/2}+ c_1 +
C_1(t) \\
 \leq & \, \tilde D +   c_1 +
C_1(t) + \frac{4 \lambda }{ p}\| u _x ( \cdot , t) \|_{L^2}^{ (p -2)
/2}\left ( D_1 + \big ( 2C_1 (t) \big )^{4/3} \left (
\int^\infty_0 | u_x( x , t)|^2  dx \right )^{1/3}\right
)^{(p+2)/4}\\
& \qquad \qquad+  \left ( D_1 + \big ( 2C_1 (t) \big )^{4/3} \left (
\int^\infty_0 | u_x( x , t)|^2  dx \right )^{1/3} \right )^{1/2}
\left ( \int^\infty_0 | u_x( x , t)|^2  dx  \right
)^{1/2}\\
\leq & \,\tilde D +  c_1 + C_1(t) + \frac{2^{(p+6)/4} \lambda }{ p}\|
u _x ( \cdot , t) \|_{L^2}^{ (p
-2) /2}\bigg ( D_1^{(p+2)/4} \\
&  \qquad \qquad + \left ( 2C_1 (t) \right )^{(p+2)/3} \left ( \int^\infty_0
| u_x( x , t)|^2  dx \right )^{(p+2)/12}\bigg
)\\
& \qquad \qquad+  \left ( D_1^{1/2} + \big ( 2C_1 (t) \big )^{2/3} \left (
\int^\infty_0 | u_x( x , t)|^2  dx \right )^{1/6} \right ) \left (
\int^\infty_0 | u_x( x , t)|^2  dx  \right
)^{1/2}\\
= & \,\tilde D +   c_1 + C_1(t) + \frac{2^{(p+6)/4} \lambda \left (
2C_1 (t) \right )^{(p+2)/3}}{ p}\|
u _x ( \cdot , t) \|_{L^2}^{2(p -1)/3} \\
&  \qquad \qquad + \frac{2^{(p+6)/4} \lambda
D_1^{(p+2)/4} }{ p}\| u _x ( \cdot , t) \|_{L^2}^{ (p -2) /2} + D_1^{1/2}\| u _x ( \cdot , t) \|_{L^2} +\left ( 2C_1 (t) \right
)^{2/3}\| u _x ( \cdot , t) \|_{L^2}^{4/3}\, ,
\end{align*}
where the first inequality is derived from the fact that $H^1 (\mathbb{R})$  is embedded in $L^\infty (\mathbb{R})$.   The 
second and third steps in the last chain of inequalities follow  from 
\eqref{1000} and \eqref{5.10.1} whilst the last 
  step is a consequence of the elementary
fact that if $a, b \geq 0$, then $(a + b ) ^ m \leq 2^{m-1} ( a ^ m + b ^m )$ when $ m \geq 1$.
When $ p < 4$,   $2(p -1)/3< 2$,  so, for any finite time interval, $\| u
_x ( \cdot , t) \|_{L^2}$ is uniformly bounded.  It follows again from
\eqref{5.10.1} that
$\| u  ( \cdot , t) \|_{L^2}$ is likewise bounded on bounded time intervals.

  Suppose $p = 4$ and let $\delta > 0$ be given, to be specified presently.  Then it follows that
\begin{align*}
 & \left ( 1- {2^{5/2} \lambda
 C_0^2 (t) }  \right )  \| u _x ( \cdot , t)
\|_{L^2}^{2} \leq
 c_2 +  c_1 + C_1(t)
+ {2^{1/2} \lambda
D_0^{3/2} }\| u _x ( \cdot , t) \|_{L^2}\\
& \qquad \quad +  D_0^{1/2}\| u _x ( \cdot , t) \|_{L^2} +\left ( 2C_0 (t)
\right )^{2/3}\| u _x ( \cdot , t)
\|_{L^2}^{4/3}\\
& = D_1 + D_2\| u _x ( \cdot , t) \|_{L^2} + D_3\| u _x ( \cdot , t)
\|_{L^2}^{4/3}\\
&\leq D_1 + \frac18\| u _x ( \cdot , t) \|_{L^2}^2 + 2 D_2^2 + \frac13 \left(
\frac{D_3}{\delta}\right)^3 + \frac23 \delta^{3/2}\| u _x ( \cdot , t) \|_{L^2}^2\, .
\end{align*}
Determine $\delta$ by  demanding  $\frac23 \delta^{3/2} = \frac18$ so that
$$
\left ( \frac{3}{4} - {2^{5/2} \lambda
 C_0^2 (t) }  \right )  \| u _x ( \cdot , t)
\|_{L^2}^{2} \leq D_1 + 2 D_2^2 + \frac13\left (\frac{ D_3}{\delta}\right)^3\!,
$$
where the right-hand side only depends on the initial and boundary data.
Since $C^2 _0 (t) =  \int^t _0 | u ( 0, s) |^2 ds$, choose $t_1$
small so that ${2^{5/2} \lambda
 C_0^2 (t_1) }\leq 1/4$.  With such a choice, if $ 0 < t \leq t_1$, then
$$
 \| u _x ( \cdot , t)
\|_{L^2}^{2} \leq 2D_1 + 4 D_2^2 + \frac23 \left( \frac{D_3}{\delta}\right)^3.
$$
 Use the solution at $t = t_1$ as the initial data and apply the
same argument to extend the solution to $t_2 > t_1$.  Since $s \geq 1$ here,
the boundary values
lie at least in $H^1_{loc}(\R^+)$. Hence, given any $T > 0$, there are positive values
$\mu = \mu(T)$, say, such that
 $\int_t^{t+\mu} |u(0,s)|^2ds$ can be made  uniformly small for all $t \in [0,T]$.
Hence, the argument just presented can be iterated at least out to time $T$. 
As $T$
was arbitrary, 
the proof is complete. $\Box$
\begin{theorem} \label{thm4.2}  Let $1 \leq s < \frac52$ be given and
assume that
\[
\mbox{ $p\geq 2\quad $ if $\quad \lambda <0\quad $ or $\quad 2\leq p\leq 4\quad $ if $\quad \lambda >0$.} \]
 Then, the IBVP (\ref{6.1}) is
globally well-posed in $H^s (\mathbb{R}^+)$ for $\phi \in
H^s (\mathbb{R}^+)$ with $h\in H^{\frac{s+3}{4}}_{loc}
(\mathbb{R}^+)$ if $ 1 \leq s \leq 2$ and $h\in H^{\frac{2s+1}{4}}_{loc}
(\mathbb{R}^+)$ if $ 2 \leq s < \frac52$.
\end{theorem}
{\bf Proof:} In (\ref{6.1}), assume that $\phi (x) \in H^2(R^+)$
and $h (t) \in H^{ \frac54}(0, T) $ satisfy the compatibility
condition $\phi (0) = h(0)$. Proposition \ref{pro-7.1} implies the global
existence of the solution $u$ which lies  in $C( [0, T], H^1 (R^+) )$,
for any $T>0$. Let $T>0$ be fixed, but arbitrary.  To prove the
 existence in $C( [0, T], H^2 (R^+) )$, take the derivative
of (\ref{6.1}) with respect to $t$ to obtain (\ref{ut-eq}) where $v=
u_t$. The initial and boundary conditions for $v$ are
$$
v(x, 0 )= i ( \phi_{xx}  + \lambda \phi |\phi |^{p-2} ) = \phi_1 (x) , \quad v (
0 , t ) = h ' (t) = h_1 (t) \, .
$$
Note that  (\ref{ut-eq}) is linear in terms of $v$. Let $ v = w + z$ be such that
$z$ satisfies
$$
i z_t + z_{xx} = 0 , \quad z ( x, 0  ) = 0, \quad z ( 0, t ) = h_1 (t)
$$
and $w$ solves
\begin{align*}
 iw_t + &w _{xx} + (\lambda p/ 2)  |u
|^{p-2} (w + z)  + (\lambda (p-2)/ 2 ) |u|^{p-4} u^2 \overline {(w+ z)}  = 0\, , \\
& w ( x, 0 ) = \phi_1 (x) , \qquad w ( 0 , t ) = 0 \,  .
\end{align*}
From \eqref{2.8}, for $s = 0 $ and any $T > 0$,
\begin{align*}
\sup _{0<t<T} \| z \| _{L^2(\R)} = \sup _{0<t<T} \| W_{bdr} (\cdot )h_1 \| _{L^2(\R)} \leq C (T) \| h_1 \|
_{H^{\frac{1}{4}} (\mathbb{R}^+)} \, .
\end{align*}
A similar identity as appears in \eqref{4.1.11} applied to $w$ gives
$$
\frac{\partial }{\partial t}  ( | w|^2 )  = -2 \im \big (w_x ( x , t) \bar w ( x, t)
\big )_x   - \lambda p | u |^{p-2} \im ( z \bar w) - \lambda(p-2)|u|^{p-4} \im \big (u^2\overline{ (w+ z) }\bar w \big )
\, .
$$
Integrating this over the half line yields
\begin{align*}
\frac{d}{d t}  \int^\infty _0 | w|^2 dx  \leq & |\lambda| p \| u \|_{H^1 ( \mathbb{R} ^+ ) } ^{p-2} \|z\|_{L^2 ( \mathbb{R} ^+ ) } \| w\|_{L^2 ( \mathbb{R} ^+ ) }   \\
& + |\lambda|(p-2)\|u\|_{H^1 ( \mathbb{R} ^+ ) }^{p-2} \big ( \| w \|_{L^2 ( \mathbb{R} ^+ ) } ^2 + \|z\|_{L^2 ( \mathbb{R} ^+ ) }  \| w\|_{L^2 ( \mathbb{R} ^+ ) }   \big ),
\end{align*}
or, what is the same,
\begin{align*}
\frac{d}{d t}  \| w\|_{L^2 ( \mathbb{R} ^+ ) }  \leq & (1/2) |\lambda| p \| u \|_{H^1 ( \mathbb{R} ^+ ) } ^{p-2} \|z\|_{L^2 ( \mathbb{R} ^+ ) }  \\
& + (1/2) |\lambda|(p-2)\|u\|_{H^1 ( \mathbb{R} ^+ ) }^{p-2} \big ( \| w \|_{L^2 ( \mathbb{R} ^+ ) } + \|z\|_{L^2 ( \mathbb{R} ^+ ) }   \big ).
\,
\end{align*}
This in turn implies by way of Gronwall's Lemma 
that $\| w(\cdot,t) \|_{C([0,T];L^2 ( \mathbb{R} ^+ )}$
is bounded. The inequality \eqref{5.1inequality} in
 \ref{pro-7.1} implies that $\| u \|_{C([0,T];H^1 ( \mathbb{R} ^+ ) )}$ is bounded by
$\alpha_0 \left (\| \phi \|_{H^1 (R^+)}+ \|h\|_{H^1(0,T)}\right )$ for some
function $\alpha _0 $.  By combining the foregoing inequalities, there obtains
\[ \sup _{0\leq t\leq T} \| v(\cdot, t)\|_{L^2(R^+)}
\leq \alpha_0 \left (\| \phi \|_{H^1 (R^+)}+ \|h\|_{H^1(0,T)}\right ) \left (\| \phi_1 \|_{L^2 (R^+)}+ \|h_1\|_{H^{\frac14}(0,T)}\right ),
\]
where $\alpha_0 : \mathbb{R}^+\to \mathbb{R}^+$ with $\alpha
_0(0)=0$  is a nondecreasing, continuous function which may depend upon $T$ as well. Thus, (\ref{6.1})
implies
\[ \sup _{0\leq t\leq T} \| u_{xx} (\cdot, t)\|_{L^2(R^+)}
\leq \alpha_0 \left(\| \phi \|_{H^1 (R^+)}+ \|h\|_{H^1(0,T)}\right ) \left(\| \phi
\|_{H^2 (R^+)}+ \|h\|_{H^{\frac54}(0,T)}\right ) ,
\]
or
\[ \sup _{0\leq t\leq T} \| u (\cdot, t)\|_{H^2(R^+)}
\leq \alpha_0 \left(\| \phi \|_{H^1 (R^+)}+ \|h\|_{H^1(0,T)}\right ) \left(\| \phi
\|_{H^2 (R^+)}+ \|h\|_{H^{\frac54}(0,T)}\right ) .
\]
By the local existence theory presented in Section 4 subject to the compatibility condition $\phi (0) = h(0)$ (see Proposition \ref{prop3.2} and Remark \ref{rem3.3.1}),
 nonlinear interpolation theory applied for $s$ in the range $ 1 < s < 2$ yields the desired result for this range of $s$  \footnote{ Here, the following interpolation result has been used. Its proof is  presented in Appendix 2.
 Let
 \[ X := \left \{ (\phi , h) \in H^1 (\R^+)\times H^{\frac34} (\R^+); \ \phi (0 =h(0)\right \}, \qquad Y := \left \{ (\phi , h) \in H^2 (\R^+)\times H^{\frac54} (\R^+); \quad \phi (0 =h(0)\right \}. \] Then,  for any $\theta$ with $0\leq \theta\leq 1$,
 \[ [ X, Y] _{\theta} =\left  \{ (\phi , h) \in H^{1+\theta}  (\R^+)\times H^{\frac{2\theta+3}{4}} (0,T); \ \phi (0) =h(0)\right \} \]}
(for details, see \cite{bsz-6} in the context of the Korteweg-de Vries  equation).   As $T>0$ was arbitrary, this in turn implies
that the theorem holds for $ 1 \leq s \leq 2$.

Now suppose that  $ 2 \leq s \leq 4$. First,
 assume $ \phi (x) \in {H^4 (R^+)}$ and $ h \in {H^{\frac94}(0,T)}$.
Take the derivative of \eqref{ut-eq} with respect to $t$ and let $v_t = v_1 = u_{tt}$. Then, the equation for $v_1$ is linear in $v_1$ with nonhomogeneous terms that
are globally defined. The initial and boundary conditions for $v_1$ are
\begin{align*}
&v_1 ( x, 0 ) =  i  \phi ''_{1} + i(\lambda p/ 2)  |\phi
|^{p-2} \phi_1 + i(\lambda (p-2)/ 2 ) |\phi |^{p-4} \phi ^2 \bar \phi_1 = \phi_2 (x) \in L^2(R^+)\, ,\\
&v_1 ( 0, t) =  h'' (t) = h'_1 (t) = h_2 ( t ) \in  {H^{\frac14}(0,T)}\, .
\end{align*}
A similar argument as that applied to
 $v ( x, t) = u_t(x,t)$ shows that $\sup _{0 < t < T} \| v_1 ( \cdot, t) \|_{L^2 ( \mathbb{R} ^+ ) } $
is bounded for any $T > 0$. Therefore,
 $\sup _{0 < t < T} \| v ( \cdot, t) \|_{H^2 ( \mathbb{R} ^+ ) } $ is bounded.
Now, consider $u$ in \eqref{ut-eq} as a fixed function in
$C ([ 0, T]; H^1(\mathbb{R}^+ ))$ and $\phi_1 ( x) , h_1 (t)$ as
functions unrelated to $\phi (x), h(x) $. Then, by the above argument,
 if $ \phi_1 (x) \in L^2 ( \mathbb {R}^+ ) ,
h _1 (t ) \in H^{\frac14} ( \mathbb {R}^+ )$,
then $v (x, t ) \in C ([ 0, T]; L^2(\mathbb{R}^+ ))$, while if
$ \phi_1 (x) \in H^2 ( \mathbb {R}^+ ) ,
h _1 (t ) \in H^{\frac54} ( \mathbb {R}^+ )$, then
$v (x, t ) \in C ([ 0, T]; H^2(\mathbb{R}^+ ))$.  This uses  only the simple compatibility
condition $\phi_1 (0) = h_1 (0)$. The usual nonlinear interpolation theory
 applied to $v$ with $s$ in the range $ 0 < s < 2$ gives the desired result for $v$ with
 $ 0 \leq s \leq 2$.  This immediately implies the advertised
 result for $u $ with $2 \leq s \leq 4$.
If $p $ is an even integer or $p $ is large, this argument
can be continued for higher values of  $ s $
(see a similar and detailed argument for the KdV equation in a
 quarter plane \cite{bsz-1}). $\qquad \Box$

\subsection{Global well-posedness on $(0,1)$}

\begin{proposition}\label{pro7.2} Assume that
 $p\geq 2$ if $\lambda <0$ and $2\leq p\leq \frac{10}{3}$ if $\lambda
>0$. Let $T>0$ be given. Then
there exists a nondecreasing continuous function $\beta :
\mathbb{R}^+\to \mathbb{R}^+$ with $\beta (0)=0$ such that any
smooth solution $u$ of (\ref{6.2}) satisfies
\[ \sup _{0\leq t\leq T} \| u(\cdot, t)\|_{H^1(0,1)} \leq \beta
\Big (\| \phi \|_{H^1 (0,1)}+ \|h_1\|_{H^1(0,T)}+\|h_2\|_{H^1(0,T)}\Big ).\]
 \end{proposition}
{\bf Proof:}
 Let $\eta (x) = x-
(1/2)$ in (\ref{4.4}) and  integrate with respect to $x$ from $0$
to $1$ to obtain
\begin{align}
&\frac{1}{2}\left ( |u_x ( 1, t) |^2 + \frac{ 2 \lambda  } { p } |u
( 1, t)|^p + |u_x ( 0, t) |^2 + \frac{ 2 \lambda  } { p } |u ( 0,
t)|^p\right ) \nonumber \\
&=  - i \int^1_0 \left ( \frac{d }{d t}  ( (x-(1/2)) u \bar
u_x )\right )dx + i (1/2) (  u (1, t)  \bar u_t (1, t) + u (0, t)  \bar u_t (0, t)) \nonumber \\
& - (  u(1, t)  \bar u _x(1, t) - u(0, t)  \bar u _x(0, t)  )
 + \int^1_0 \left (  2| u_x(x, t)  |^2 - \lambda ( 1- (2/p)) | u (x, t) |
 ^p\right ) dx \, .\label{4.5}
\end{align}
In the following, we again use $D$ as a constant dependent on the
initial and boundary data, $c$ as a constant only dependent on the
initial data and $C(t)$ as a constant only dependent on the
boundary data, while $C$ is just a fixed constant, independent of the
initial and boundary data. Integrate (\ref{4.5}) with
respect to $t$ from $0$ to $t$  to derive
\begin{align}
&\int^t_0 \left ( |u_x ( 1, s) |^2  + |u_x ( 0, s) |^2 \right ) ds
= D_0 +  2 \int^t_0\int^1_0\left ( 2 |u_x ( x, s) |^2 - \lambda ( 1- (2/p))
 |u ( x, s)|^p \right ) dx ds\nonumber \\
&  \qquad - 2i \int ^1_0 ((x-(1/2))u (x, t) \bar u _x( x, t) ) dx -2
\int^t_0 \Big (  u(1, s) \bar u _x(1, s) - u(0, s) \bar u _x(0,s)
\Big )ds \, .\label{4.6}
\end{align}
Consider the cases $\lambda > 0 $ and $\lambda < 0$ separately.

\medskip

\noindent{\bf (a) $\lambda < 0$}

For this case, (\ref{4.6}) gives
\begin{align*}
&\int^t_0 \left ( |u_x ( 1, s) |^2  + |u_x ( 0, s) |^2 \right ) ds
\leq D_0 +  C \int^t_0\int^1_0\left ( |u_x ( x, s) |^2 +  |u ( x, s)|^p \right ) dx ds\\
&  \qquad + \int^1_0 \left |  u(x, t)  \bar u_x (x, t) )\right |dx +
2\int^t_0 \Big |  u(1, s)  \bar u _x(1, s) - u(0, s)  \bar u _x(0,s)
\Big |ds\\
& \leq D_1 +  C \int^t_0\int^1_0\left ( |u_x ( x, s) |^2 +  |u ( x, s)|^p \right ) dx ds\\
 &\qquad + \left ( \int^1_0  |  u(x, t)|^2 dx \right )^{1/2}   \left (  \int^1_0
|  u_x (x, t) |^2dx\right )^{1/2} + \frac{1}{2}\int^t_0 \Big (  | u
_x(1, s)|^2 +  | u _x(0,s)|^2
\Big )ds\, ,\\
\end{align*}
which implies
\begin{align*}
&\int^t_0 \left ( |u_x ( 1, s) |^2  + |u_x ( 0, s) |^2 \right ) ds
\leq 2 D_1 +  2 C \int^t_0\int^1_0\left ( |u_x ( x, s) |^2 +  |u (
x, s)|^p \right ) dx ds \\
&\qquad \qquad + 2 \left ( \int^1_0  |  u(x, t)|^2 dx \right )^{1/2}
\left ( \int^1_0 |  u_x (x, t) |^2dx\right )^{1/2}\, .
\end{align*}
We use techniques that are by now familiar to obtain from
(\ref{4.1}) that
\begin{align*}
\int^1_0 &|u (x, t) |^2 dx = \int^1_0 |u (x, 0) |^2 dx
- 2\im \int^t _0 (u_x ( 1 , s) \bar u ( 1, s) -u_x ( 0 , s) \bar u ( 0, s) )ds \\
\leq & \int^1_0 |u (x, 0) |^2 dx +2 \left ( \int^t _0 \left ( |u( 1,
s)|^2 + |u ( 0 , s)|^2\right ) ds \right )^{1/2} \left ( \int^t _0
\left ( |u_x ( 1 , s)|^2 + |u_x ( 0 , s)|^2\right ) ds \right
)^{1/2}\\
\leq &\; c_0 +2C_0 (t) \bigg (2 D_1 +  2 C \int^t_0\int^1_0\left (
|u_x ( x, s) |^2 +  |u (
x, s)|^p \right ) dx ds \\
&\qquad  + 2 \left ( \int^1_0  |  u(x, t)|^2 dx \right
)^{1/2} \left ( \int^1_0 |  u_x (x, t) |^2dx\right )^{1/2} \bigg
)^{1/2}\\
 \leq & \; c_0 +2\sqrt{2} C_0
(t)\bigg ( D_1^{1/2} + \left ( C \int^t_0\int^1_0\left ( |u_x ( x,
s) |^2 + |u (
x, s)|^p \right ) dx ds\right )^{1/2} \\
&\qquad + \left ( \int^1_0  |  u(x, t)|^2 dx \right
)^{1/4} \left ( \int^1_0 |  u_x (x, t) |^2dx\right )^{1/4}\bigg ) \\
\leq & \; D_2  + (1/4) \int^1_0 |u ( x , t)|^2 dx  +(3/4) \left (
2\sqrt{2} C_0 (t) \left ( \int^1_0 | u_x( x , t)|^2  dx \right
)^{1/4}\right )^{4/3}\\
&\qquad  + 2\sqrt{2} C_0 (t)   \left ( C \int^t_0\int^1_0\left (
|u_x ( x, s) |^2 +  |u ( x, s)|^p \right ) dx ds\right )^{1/2} \! ,
\end{align*}
or
\begin{align*}
\int^1_0 |u (x, t) |^2 dx &\leq  D_3  +  4 C_0^{4/3} (t) \left (
\int^1_0 | u_x( x , t)|^2  dx \right
)^{1/3}\\
&\qquad  + 4 C_0 (t)   \left ( C \int^t_0\int^1_0\left ( |u_x ( x,
s) |^2 +  |u ( x, s)|^p \right ) dx ds\right )^{1/2}\, .
\end{align*}

To obtain an estimate for $u_x$ for the case  $\lambda < 0$,  integrate
(\ref{4.2}) with respect to $x$ and $t$ to reach
\begin{align*}
& \int^1_0\left (  |u_x(x, t)  |^2 + \frac{2 |\lambda |}{ p} | u(x,
t) |^p
 \right ) dx \leq  \int^1_0\left (  |u_x(x, 0)  |^2 + \frac{2 |\lambda |}{ p} | u(x, 0)
|^p
 \right ) dx \\
 &\qquad + 2 \int^t_0 \Big | u_x ( 1, s) \bar u_s ( 1, s) -
u_x ( 0, s) \bar u_s ( 0,s)\Big | ds\\
&\leq D_4 +\int^t_0  |u_x ( 1, s)|^2 ds + \int^t_0  |u_x ( 0,
s)|^2ds \\
&\leq D_4 + 2 D_1 +  2 C \int^t_0\int^1_0\left ( |u_x ( x, s) |^2 +
|u (
x, s)|^p \right ) dx ds \\
&\qquad \qquad + 2 \left ( \int^1_0  |  u(x, t)|^2 dx \right )^{1/2}
\left ( \int^1_0 |  u_x (x, t) |^2dx\right )^{1/2}\\
 &\leq D_4 + 2 D_1 +  2 C \int^t_0\int^1_0\left ( |u_x ( x, s)
|^2 + |u (
x, s)|^p \right ) dx ds \\
&\qquad  + 2 \bigg ( D_3  +  4 C_0^{4/3} (t) \left ( \int^1_0 | u_x(
x , t)|^2  dx \right
)^{1/3}\\
&\qquad  + 4 C_0 (t)   \left ( C \int^t_0\int^1_0\left ( |u_x ( x,
s) |^2 +  |u ( x, s)|^p \right ) dx ds\right )^{1/2} \bigg )^{1/2}
\left ( \int^1_0 |  u_x (x, t) |^2dx\right )^{1/2}\\
&\leq D_5 + 2 C \int^t_0\int^1_0\left ( |u_x ( x, s) |^2 + |u (
x, s)|^p \right ) dx ds \\
&\qquad  + 2D_3^{1/2}\left ( \int^1_0 |  u_x (x, t) |^2dx\right
)^{1/2} + 4C_0^{2/3} (t) \left ( \int^1_0 | u_x( x , t)|^2  dx
\right )^{2/3}\\
&\qquad + 2(C_0 (t))^{1/2}   \left ( C \int^t_0\int^1_0\left ( |u_x
( x, s) |^2 +  |u ( x, s)|^p \right ) dx ds\right )^{1/4}  \left (
\int^1_0 |  u_x (x, t) |^2dx\right )^{1/2}\\
&\leq D_6 + \frac{1}{2}  \int^1_0 |  u_x (x, t) |^2dx +3 C
\int^t_0\int^1_0\left ( |u_x ( x, s) |^2 + |u ( x, s)|^p \right ) dx
ds\, ,
\end{align*}
where use has been made of Young's inequality.  It follows that 
$$
\int^1_0\Big (  |u_x(x, t)  |^2 +  | u(x, t) |^p
\Big  ) dx \leq D_7  +C
\int^t_0\int^1_0\Big  ( |u_x ( x, s) |^2 + |u ( x, s)|^p \Big  ) dx
ds\, .
$$
for suitable constants.   Gronwall's lemma then provides a global bound on the
solution $u$ in $H^1(0, 1)$.

\medskip
\noindent{\bf (b) $\lambda > 0$}

From (\ref{4.6}) with $p \geq 2$, it happens that 
\begin{align*}
&\int^t_0 \left ( |u_x ( 1, s) |^2  + |u_x ( 0, s) |^2 \right ) ds
+\int^t_0\int^1_0 \lambda \big ( 1- (2/p)\big )
 |u ( x, s)|^p dx ds \\
 &
= D_0 +  4 \int^t_0\int^1_0 |u_x ( x, s) |^2  dx ds  - 2i \int ^1_0
\big (x-(1/2)\big )u (x, t) \bar u _x( x, t)  dx \\
&\qquad \qquad -2 \int^t_0 \Big (  u(1, s) \bar u _x(1, s) - u(0, s)
\bar u _x(0,s)
\Big )ds\\
& \leq D_1 +  C \int^t_0\int^1_0 |u_x ( x, s) |^2  dx ds+ \left (
\int^1_0  |  u(x, t)|^2 dx \right )^{1/2}   \left (  \int^1_0 |  u_x
(x, t) |^2dx\right )^{1/2} \\
 &\qquad \qquad + \frac{1}{2}\int^t_0 \Big (  | u _x(1,
s)|^2 +  | u _x(0,s)|^2
\Big )ds\, ,
\end{align*}
which implies
\begin{align*}
&\int^t_0 \left ( |u_x ( 1, s) |^2  + |u_x ( 0, s) |^2 \right ) ds
\leq 2 D_1 +  2 C \int^t_0\int^1_0 |u_x ( x, s) |^2  dx ds \\
&\qquad \qquad + 2 \left ( \int^1_0  |  u(x, t)|^2 dx \right )^{1/2}
\left ( \int^1_0 |  u_x (x, t) |^2dx\right )^{1/2}\!\! .
\end{align*}
By the same argument as in the case with $\lambda < 0$, it is seen that
\begin{align*}
\int^1_0 |u (x, t) |^2 dx &\leq  D_3  +  4 C_0^{4/3} (t) \left (
\int^1_0 | u_x( x , t)|^2  dx \right )^{1/3}  + 4 C_0 (t) \left ( C
\int^t_0\int^1_0 |u_x ( x, s) |^2 dx ds\right )^{1/2}\!\! .
\end{align*}
The estimate for $\| u_x\|_{H^1(0,1)}$ can be obtained
 from (\ref{4.2}) as follows:
\begin{align}
 \int^1_0 &|u _x(x, t)  |^2 dx =   \frac{2 \lambda }{ p} \int^1_0
 | u (x,t) |^p dx + \int^1_0 \left ( |u _x (x, 0) |^2 -\frac{2 \lambda }{
 p} | u (x,0) |^p \right ) dx \nonumber \\
 & \qquad +
2 \re \int^t_0 \left ( u_x ( 1, s) \bar u_s ( 1, s)  - u_x ( 0, s) \bar u_s ( 0, s) \right ) ds\nonumber \\
\leq & \,\frac{4 \lambda }{ p} \Big ( \| u ( \cdot , t) \|_{L^2} \| u
_x ( \cdot , t) \|_{L^2}\Big )^{ (p -2) /2}  \int^1_0
 | u (x,t) |^2 dx  + D_4 +  \int^t_0  |u_x ( 1, s)|^2 ds +
\int^t_0  |u_x ( 0, s)|^2 ds\nonumber \\
\leq & \, D_4 + \frac{4 \lambda }{ p}\| u _x ( \cdot , t) \|_{L^2}^{ (p
-2) /2}\bigg [D_3  +  4 C_0^{4/3} (t) \left ( \int^1_0 | u_x( x ,
t)|^2 dx \right
)^{1/3}\nonumber \\
&\qquad  + 4 C_0 (t)   \left ( C \int^t_0\int^1_0 |u_x ( x, s) |^2
dx ds\right )^{1/2} \bigg
]^{(p+2)/4} + 2 D_1 +  2 C \int^t_0\int^1_0 |u_x ( x, s) |^2 dx ds \nonumber \\
&\qquad  + 2 \left ( \int^1_0  |  u(x, t)|^2 dx \right )^{1/2}
\left ( \int^1_0 |  u_x (x, t) |^2dx\right )^{1/2}\nonumber
\\
 \leq & \, D_4 + \frac{4 \lambda }{ p}\| u _x ( \cdot , t)
\|_{L^2}^{ (p -2) /2} \bigg [D_3  +  4 C_0^{4/3} (t) \left (
\int^1_0 | u_x( x , t)|^2 dx \right
)^{1/3}\nonumber \\
&\qquad  + 4 C_0 (t)   \left ( C \int^t_0\int^1_0 |u_x ( x, s) |^2
dx ds\right )^{1/2} \bigg
]^{(p+2)/4}\nonumber \\
& \qquad + 2 D_1 +  2 C \int^t_0\int^1_0 |u_x ( x, s) |^2 dx ds  + 2 \bigg
[D_3  +  4 C_0^{4/3} (t) \left ( \int^1_0 | u_x( x , t)|^2 dx \right
)^{1/3}\nonumber \\
&\qquad  + 4 C_0 (t)   \left ( C \int^t_0\int^1_0 |u_x ( x, s) |^2
dx ds\right )^{1/2} \bigg ]^{1/2}
\left ( \int^1_0 |  u_x (x, t) |^2dx\right )^{1/2}\nonumber \\
\leq & \, D_4 + C \| u _x ( \cdot , t) \|_{L^2}^{ (p -2) /2} \bigg
[D_3^{(p+2)/4} + \left ( 4 C_0^{4/3} (t) \left ( \int^1_0 | u_x( x ,
t)|^2 dx \right
)^{1/3}\right )^{(p+2)/4}\nonumber \\
&\qquad  + \left ( 4 C_0 (t)   \left ( C \int^t_0\int^1_0 |u_x ( x,
s) |^2  dx ds\right )^{1/2} \right )^{(p+2)/4} \bigg
]\nonumber \\
& \qquad + 2 D_1 +  2 C \int^t_0\int^1_0|u_x ( x, s) |^2  dx ds + 2 \bigg
[D_3^{1/2}  +  2 C_0^{2/3} (t) \left ( \int^1_0 | u_x( x , t)|^2 dx
\right
)^{1/6}\nonumber \\
&\qquad  + 2 C_0^{1/2} (t)   \left ( C \int^t_0\int^1_0 |u_x ( x, s)
|^2  dx ds\right )^{1/4} \bigg ]
\left ( \int^1_0 |  u_x (x, t) |^2dx\right )^{1/2}\nonumber \\
\leq & \, D_5 + D_6 \| u _x ( \cdot , t) \|_{L^2}^{ (p -2) /2} + C
C_0^{(p+2)/3} (t)\| u _x ( \cdot , t) \|_{L^2}^{ 2(p -1) /3} \nonumber \\
& \qquad+ C_1 (t)\| u _x ( \cdot , t) \|_{L^2}^{ (p -2) /2} \left (
\int^t_0\int^1_0|u_x ( x, s) |^2 dx ds\right )^{(p+2)/8}+ 2 C
\int^t_0\int^1_0 |u_x ( x, s)
|^2  dx ds \nonumber \\
& \qquad + D_7\| u _x ( \cdot , t) \|_{L^2} + C_2 (t) \| u _x ( \cdot , t)
\|_{L^2}^{4/3} + C_3 (t) \| u _x ( \cdot , t) \|_{L^2}\left (
\int^t_0\int^1_0 |u_x ( x, s) |^2  dx ds\right )^{1/4}\!.\label{4.7}
\end{align}
For $2\leq p \leq \frac{10}{3}$, Young's inequality $ab \leq
(1/m) a^m + (1/n) b^n$ with $m^{-1} + n^{-1} = 1$ leads from
(\ref{4.7}) to
$$
\int^1_0 |u _x(x, t)  |^2 dx  \leq D_8 + D_9 \int^t_0\int^1_0 |u_x (
x, s) |^2  dx ds
$$
which, by Gronwall's lemma, gives a uniform bound for
$\int^1_0 |u _x(x, t)  |^2 dx$ on the interval $ 0\leq t \leq T$. As 
the time $T$ is arbitrary, the proof is complete. $\Box$

Using the same argument as that for proving Theorem
\ref{thm4.2} leads to the following global well-posedness
result.
\begin{theorem} Assume that
\[ \mbox{ $p\geq 3$ \quad if \quad   $\lambda <0$ \quad or \quad  $3\leq p\leq
 \frac{10}{3}$ \quad  if
\quad   $\lambda
>0$}\] and let $1 \leq s < 5/2$ be given. Then the IBVP (\ref{6.2}) is
globally well-posed in $H^s (0,1)$ with $\phi \in H^s
(0,1)$ and $h_1, \ h_2 \in H^{\frac{s+1}{2}}_{loc} (\mathbb{R}^+)$ subject to the compatibility conditions on $\phi, h_1$ and $h_2$.
\end{theorem}

\bigskip

\section{Appendices}
\subsection{Appendix 1}

\setcounter{equation}{0}

The following Lemma is used to obtain an estimate in the proof of Proposition \ref{prop3.6}.
\medskip
\noindent {\bf Lemma A-1:} Let $\psi $ be an even, non-negative, $C^\infty$ cut-off function with $\mbox{supp}(\psi ) \subset [-1, 1]$ and with $\psi(x) \equiv 1 $ for
$|x| \leq \frac12$. Suppose also that $\psi$ is strictly decreasing on $[\frac12, 1]$. There exists a
constant $C>0$ such that for any $g\in H^{\frac12}_{00} (\R^+)$,
$$
\sum_{n=1}^\infty \left | \int^\infty_0 f(\mu ) \frac{1}{\mu -n } \Big (
1 - \psi (n^2 - \mu^2) \Big ) d\mu \right |^2 \leq  C \int^\infty_0(\mu
+1 ) |f(\mu )|^2 d \mu\, ,
$$
where $f$ is the Fourier transform of the extension by zero of $g$ to all of $\R$.

\bigskip

\noindent
{\bf Proof:} Write
\begin{align*}
\sum_{n=1}^\infty &\left | \int^\infty_0 f(\mu ) \frac{1}{\mu -n } \big(
1 - \psi (n^2 - \mu^2) \big) d\mu \right |^2\\
& = \sum_{n=1}^\infty \left | \left ( \int^{n-1}_0 + \int^{n+1}
_{n-1} + \int^\infty_{n+1} \right ) f(\mu ) \frac{1}{\mu -n } \big( 1 -
\psi (n^2 - \mu^2) \big) d\mu \right |^2 \\
&\; \leq I_1 + I_ 2 + I_3 \, .
\end{align*}

Since the estimates for $I_1 $ and $I_3$ are similar, we only study $I_3$.
Let $\alpha, \beta > 0$ so the Cauchy-Schwartz inequality implies that
\begin{align*}
I_3 & = \sum_{n=1}^\infty \left |  \int^\infty_{n+1}  f(\mu )
\frac{1}{\mu -n }\big ( 1 -
\psi (n^2 - \mu^2) \big) d\mu \right |^2 \leq \sum_{n=1}^\infty  \int^\infty_{n+1}  \frac{|f(\mu )|^2 \mu
^{2\alpha} }{|\mu -n|^{2-2\beta} } d\mu \int^\infty_{n+1} \frac{1
}{|\mu -n|^{2\beta} \mu^{2\alpha}} d\mu\, .
\end{align*}
If $2\alpha + 2\beta > 1$, then
\begin{align*}
  \left | \int^\infty_{n+1} \frac{1}{|\mu -n|^{2\beta} \mu^{2\alpha}}
 d\mu\right |
 & =\left | \left ( \int^{3n} _{n+1} + \int^\infty_{3n} \right ) \frac{1}{|\mu -n|^{2\beta} \mu^{2\alpha}}
 d\mu\right | \\
 & \leq C \left ( \int^{3n} _{n+1}\frac{1}{|\mu -n|^{2\beta} n^{2\alpha}}
 d\mu + \int^\infty_{3n}  \frac{1}{ \mu^{2\alpha + 2\beta }}
 d\mu \right ) \leq C n ^{1-2\alpha - 2\beta} \leq C,
\end{align*}
 where $C$ is independent of $n$. It thus transpires that if $2 - 2\beta > 1$, then
\begin{align*}
I_3 &\leq  C \sum_{n=1}^\infty  \int^\infty_{n+1}  \frac{|f(\mu )|^2
\mu ^{2\alpha} }{( |\mu -n| + 1) ^{2-2\beta} } d\mu \leq  C
 \int^\infty_{0}  {|f(\mu )|^2 \mu
^{2\alpha} } \sum_{n=1}^\infty \frac{1}{( |\mu -n| + 1) ^{2-2\beta}
} d\mu\\
&\leq C\int^\infty_{0}  {|f(\mu )|^2 \mu ^{2\alpha} } d\mu.
\end{align*}
Choosing $\alpha = \frac14$ and, say, $\beta = \frac38$ yields the advertised bound.
\vspace{.1 cm}

To study $I_2$,  note that in the integrals, the integrand vanishes unless $\mu \geq
\sqrt{n^2 + 1/2}$ or $ 0 \leq \mu \leq \sqrt{n^2 -1/2}$. Consequently, it
must be the case  that
\begin{align*}
I_2 & \leq \sum_{n=1}^\infty\left ( \left (  \int^{n+1} _{\sqrt{n^2
+ 1/2}} +  \int^{\sqrt{n^2 -1/2}} _{n-1}\right ) \frac{|f(\mu
)|}{|\mu
-n |}  d\mu \right )^2 \\
& \leq \sum_{n=1}^\infty\bigg (   \int^{n+1} _{\sqrt{n^2 + 1/2}}
|f(\mu )|^2 d\mu   \int^{n+1} _{\sqrt{n^2 + 1/2}} |\mu -n |^{-2} d
\mu \\
&\qquad \quad + \int^{\sqrt{n^2 -1/2}} _{n-1}|f(\mu )|^2 d\mu
\int^{\sqrt{n^2
-1/2}} _{n-1}|\mu -n |^{-2} d \mu \bigg )\\
&\leq C \sum_{n=1}^\infty\left (   \int^{n+1} _{n} n |f(\mu )|^2
d\mu  + \int^{n} _{n-1} n |f(\mu )|^2 d\mu  \right ) \\
& \leq C \sum_{n=1}^\infty \left (   \int^{n+1} _{n} \mu |f(\mu )|^2
d\mu  +
\int^{n} _{n-1} (\mu +1)  |f(\mu )|^2 d\mu  \right )\\
& \leq C \int^\infty_0 (\mu +1 )|f(\mu )|^2 d\mu\, .
\end{align*}
The lemma is proved. $\Box$

\bigskip
The following example shows the optimality of the assumption $h\in H^{1/2} (0, T )$ in (\ref{estimate}) and (\ref{estimate-}).  This result then implies that the assumptions on $( h _1, h_2 )$ in Theorem \ref{th-finite} are optimal.

\noindent {\bf Example A-2:}
Notice that if (\ref{estimate}) or (\ref{estimate-}) holds, then
$$
\| u_h \|_{L^2 (\Omega_T)} = \| u_h \|_{L^2 ((0, 1) \times (0, T))} \leq C_T \| h \|_{H^{\frac12}(0, T) }\, ,
$$
where we recall for the reader's convenience that
\begin{eqnarray*}
u_h &=&  \sum ^{\infty }_{n=1} 2i n\pi e^{-i(n\pi )^2 t} \int ^t_0
e^{i(n\pi )^2\tau } h(\tau ) d \tau \sin n\pi x\\ \\
&=& \sum ^{\infty }_{n=-\infty } n\pi e^{-i(n\pi)^2t +in\pi x} \int
^t_0 e^{i(n\pi )^2 \tau } h(\tau ) d\tau
\end{eqnarray*}
(see \eqref{4.9.1}).
Assume that $h(t)$ has the Fourier series expansion
\[ h(t)= \sum ^{\infty }_{k=-\infty} e^{-\pi ^2 ikt} a_k\quad \mbox{ with }
\quad a_k = \int ^{\frac{2\pi}{\pi ^2}}_0 e^{\pi ^2ikt} h(t) dt \  .\]
It follows that
\begin{align*}
&u_h = \sum ^{\infty }_{n=-\infty } n\pi e^{-i(n\pi )^2 t+in\pi x}
\int ^t_0 \sum _k e^{i(n^2-k)\pi ^2 \tau } a_k d\tau \\ \\
&= \left ( \sum ^{\infty }_{n=-\infty } n\pi e^{-i(n\pi )^2 t+in\pi  x}
\sum _{k\ne n^2} \frac{e^{i(n^2-k)\pi ^2 t}-1}{n^2 -k} a_k \right ) + \left ( \sum
^{\infty }_{n=-\infty } n\pi e^{-i(n\pi )^2 t+in\pi x}  t a_{n^2} \right ) \\ \\
&= \left ( \sum ^{\infty }_{n=-\infty } n\pi e^{in\pi x} \sum _{k\ne n^2}
\frac{e^{-ki\pi ^2 t}-e^{-in^2 \pi ^2 t}}{n^2-k} a_k\right ) + \left ( \sum ^{\infty
}_{n=-\infty } n\pi e^{-i(n\pi )^2 t+in\pi x} t a_{n^2}\right )\\ \\
&= \left ( \sum ^{\infty }_{n=-\infty } n\pi e^{in\pi x} \sum _{k\ne n^2}
\frac{e^{-ki\pi ^2 t} }{n^2-k} a_k \right )+\left ( \sum ^{\infty }_{n=-\infty }
n\pi e^{-i(n\pi )^2 t+in\pi x} \left ( ta_{n^2} - \sum _{k\ne n^2}
\frac{a_k}{n^2-k}\right ) \right )\!.
\end{align*}
Choose $h(t)$ so that
\[ a_{n^2} = \int ^{\frac{2}{\pi}}_0 e^{\pi ^2 in^2 t} h(t) dt =0\, ,\quad n\in \mathbb{Z}\, .
\]
Then, the last formula condenses to
\[
u_h = \sum ^{\infty }_{n=-\infty } n\pi e^{in\pi x} \sum _{k\ne n^2}
\frac{e^{-ki\pi ^2 t} }{n^2-k} a_k +\sum ^{\infty }_{n=-\infty }
n\pi e^{-i(n\pi )^2 t+in\pi x} \left ( \sum _{k\ne n^2}
\frac{a_k}{k-n^2} \right )\!.\]
As $k \neq n^2$, the exponentials $e^{in\pi x -i k\pi^2 t}$ and  $e^{i(n\pi x - (n\pi)^2) t}$ are orthogonal, whence
\begin{eqnarray*}
\|  u_h \|^2_{L^2\left ((0,1)\times \left (0, \frac{2}{\pi}\right )\right )} &=& \sum ^{\infty }_{n=-\infty } \sum _{k\ne n^2}
n^2 \pi ^2 \frac{a_k^2}{(n^2-k)^2} + \sum ^{\infty }_{n=-\infty }
(n\pi )^2 \left ( \sum _{k\ne n^2} \frac{a_k}{k-n^2}\right )^2 \\ \\
&\geq & \sum ^{\infty }_{n = -\infty } n^2 \pi ^2 a_{n^2+1}^2,
\end{eqnarray*}
the latter inequality obtained by only considering the terms where $k = n^2 + 1$.

 If there were a constant $C$ such that for all $h\in H^{\alpha}\left (0,
\frac{2}{\pi}\right ),$
$ \| u_h \| ^2_{L^2\left ((0,1)\times \left (0, \frac{2}{\pi}\right )\right )} \leq C \|h\|^2_{H^{\alpha}\left (0,
\frac{2}{\pi}\right )}$, then it would follow that  $\alpha \geq \frac12$.
Suppose instead that there is a constant $C$ such that
\[ \| u_h \| ^2_{L^2\left ((0,1)\times \left (0, \frac{2}{\pi}\right )\right )} \leq C \|h\|^2_{H^{\alpha}\left (0,
\frac{2}{\pi}\right )} \qquad  \mbox{for some $\alpha$ with $0<\alpha < \frac12$}\, . \]
Define the function  $h$  by its Fourier series, {\it viz.}
 \[
h(t) = \sum _{n\ne 0} \frac{1}{|n|^{\beta }} e^{-\pi ^2 i
(n^2+1)t} \, .
\]
For $h$ to lie in $ H^{\alpha } \left (0, \frac{2}{\pi}\right )$, we need
\[
\sum _{n\ne 0} \left | \frac{(n^2+1)^{\alpha }}{|n|^{\beta }}
\right |^2 < +\infty \, ,
\] or
$ 2\beta -4\alpha >1 $
which implies that
$\beta > 2\alpha +\frac12 .$
But, for this $h$,
\[
\|u_h\| ^2_{L^2\left ((0,1)\times \left (0, \frac{2}{\pi}\right )\right )} \geq \sum ^{\infty }_{n=-\infty } n^2 \pi ^2
\frac{1}{|n|^{2\beta}} = \sum ^{\infty }_{n=-\infty , \ n\ne 0}
\frac{\pi^2 }{|n|^{2\beta -2}}\, .
\]
Since $\alpha < \frac12$,  $\beta $ can be chosen so  that $2\alpha +\frac12 < \beta < \frac32 $.
For such a value of $\beta$, it is clear that
\[
\sum ^{\infty }_{n=-\infty , \ n\ne 0}
\frac{1}{|n|^{2\beta -2}}=+\infty.
\]
The partial sums
\[
h_k (t) = \sum ^{|n|=k}_{n\ne 0} \frac{1}{|n|^{\beta }}e^{-\pi ^2
i (n^2+1) t}
\]
lie in  $ C^{\infty }$ and therefore, according to our hypothesis,
$$\| u_{h_k}\|^2_{L^2\left ((0,1)\times \left (0, \frac{2}{\pi}\right )\right )} \leq C \| h_k \|^2_{H^{\alpha}}\,  .$$
But, as $k\to \infty $, the right side of the last inequality is bounded while the left side
tends to $\infty $, which is a contradiction.  Hence, we must have $ \alpha \geq \frac12.$

\subsection{Appendix 2}
Let
\[ X=\left \{  (\phi, h)\in H^1 (\R^+)\times H^{\frac34} (\R^+); \quad \phi (0)=h(0) \right \},\] \[ Y= \left \{  (\phi, h)\in H^2 (\R^+)\times H^{\frac54} (\R^+); \quad \phi (0)=h(0) \right \},\]
and
\[ X^*=H^1 (\R^+)\times H^{\frac34}(\R^+),  \quad Y^*= H^2 (\R^+)\times H^{\frac54} (\R^+) .\]
While  it is  well-known  (cf. \cite{LM1972}) that for any $\theta$ 
with $0\leq \theta \leq 1$,
\[ \left [ X^*, Y^*\right ]_{\theta} = H^{1+ \theta} (\R^+)\times H^{\frac{2\theta+3}{4}} (\R^+),\]
however,  as pointed out by an anonymous  referee,  it seems that no rigorous proof  can be found in literature for the interpolation result
\begin{equation} \label{y-1}  \left [ X, Y\right ]_{\theta} =\left \{  (\phi, h)\in H^{1+\theta}  (\R^+)\times H^{\frac{2\theta +3}{4}}(\R^+); \quad \phi (0)=h(0) \right \} \end{equation}
 used in our analysis.  It is mentioned (in a much more general setting) 
as ``most likely'' true
in the book of Lions-Magenes \cite{LM1972-2} (Chapter 4, Section 14, remark after Theorem
14.1).  The following  short proof of (\ref{y-1}) was suggested  by the referee.

First, it is claimed that there exists a  bounded linear  ``lifting'' operator $L$ from the space
\[ Z_s:= \left \{  (\phi, h)\in H^{1+s}  (\R^+)\times H^{\frac{2s+3}{4}}(\R^+); \quad \phi (0)=h(0) \right \} \]
to $H^{s+2, \frac{s+2}{2}}(\R^+\times \R^+ ):= L^2_t (\R^+; H^{s+2}_x (\R^+) \cap H_t^{\frac{s+2}{2}} (\R^+; L^2_x (\R^+))$ for $0\leq s\leq 1$ such that
$ w= L(\phi, h) \in H^{s+2, \frac{s+2}{2}}(\R^+\times \R^+ )$ for any $(\phi, h) \in Z_s $ and
\[ w(x,0)=\phi (x), \qquad w(0,t) = h(t) .\]
Then
\[ \mathbb{T}\circ L = I\]
where $I$ denotes the identity operator and $\mathbb{T}$ is the trace operator defined by
\[ \mathbb{T}: H^{s+2, \frac{s+2}{2}}(\R^+\times \R^+ ) \to Z_s , \quad \mathbb{T} w= (w(x,0), w(0,t)).\]
One has (\cite{LM1972-2} Proposition 2.1,  Chapter 4)
\[ \left [ H^{2, 1}(\R^+ \times \R^+), H^{3, \frac32} (\R^+\times \R^+)\right ]_{\theta} = H^{2+\theta, 1+\frac{\theta}{2}} (\R^+ \times \R^+)\]
for $0\leq \theta\leq 1$. Consequently,  $[X,Y]_{\theta}$ can be identified with $\mathbb{T} \left ( H^{2+\theta, 1+\frac{\theta}{2}} (\R^+\times \R^+)\right )$ which is exactly $Z_{\theta}$.

It remains to prove the existence of the lifting operator $$L: Z_s \to H^{s+2, 1+ \frac{s}{2}}(\R^+\times \R^+)$$ for $0\leq s\leq 1$. To this end,
consider the following IBVP 
\begin{equation} \label{y-2}
\begin{cases}
u_t =u_{xx}, \quad x\in \R^+, \ t\in \R^+,\\
u(x,0) = \phi (x), \quad u(0,t) =h(t) ,  \quad  x\in \R^+, \ t\in \R^+\, ,
\end{cases}
\end{equation}
for the heat equation, where $\phi \in H^{s+1} (\R^+) , h(t ) \in H^{\frac{2s+3}4} (\R^+) $ with $ 0 \leq s \leq 1$ and $\phi(0) = h(0)$.
The existence of the solution $u (x, t) \in H^{s+2, 1+ \frac{s}{2}}(\R^+\times \R^+)$ for \eqref{y-2} is established in Theorems 6.1 and 6.2 in Chapter 4 of  \cite{LM1972-2}.
Therefore, given $s$ with 
 $0\leq s\leq 1$ and  $(\phi , h)\in Z_s$, we 
may define the lifting operator $L$ by
\[ L (\phi , h):= u, \]
where $u$ is the solution of the IBVP (\ref{y-2}).

\section{Acknowledgment} JLB and SMS were partially supported by the US National Science
Foundation. BYZ was partially supported by a grant from the Simons Foundation (201615) and NSF of China
(11231007, 11571244). JLB also thanks the
Universit\'e de Paris Nord and the Ulsan National Institute of Science and Technology for
hospitality and very good working conditions during parts of the writing
phase of this project.  We all tender heartfelt thanks to anonymous referees for
 careful readings of the script and many helpful comments, corrections and suggestions.


\addcontentsline{toc}{section}{References}
\begin{thebibliography}{9}

\bibitem{hammack} M. J. Ablowitz, J. Hammack, D. Henderson and C. Schober, { Modulated
periodic Stokes waves in deep water}, \emph{Phys. Rev. Letters} {\bf 84} (2000) 887--890.

\bibitem{Audiard-1} {C. Audiard}, {Non-homogeneous boundary value problems for linear
dispersive equations},   \emph{Comm. Partial Diff. Equations} {\bf 37} (2012) {1--37}.

\bibitem{Audiard-2} {C. Audiard}, {On the non-homogeneous boundary value problem for Schr\"odingier equations},
 \emph{Discrete Continuous Dynamical Systems} {\bf 33} (2013)
 3861--3884.

\bibitem{BCSZ}  J.L. Bona, H. Chen, S.-M. Sun and B.-Y. Zhang,   Comparison of quarter-plane and two-point boundary-value problems:  The BBM-equation,  {\em Discrete Continuous Dynamical Systems} {\bf 13} (2005) 921--940.

\bibitem{BCL2005}
J. L. Bona, T. Colin and D. Lannes, Long wave approximations for water waves, \emph{Arch. Rational Mech. Anal.} {\bf 178} (2005) 373--410.

\bibitem{bpss}  J. L. Bona, G. Ponce, J.-C. Saut and C. Sparber, Dispersive blow-up for
nonlinear Schr\"odinger equations revisited, {\em  J. Math. Pures Appliq.} {\bf 102} (2014)
 782--811.

\bibitem{BSaut}  J. L. Bona and J.-C. Saut, Dispersive blow-up II.  Schr\"odinger--type equations,
optical and oceanic rogue waves, {\em Chinese Ann. Math. Series B} {\bf 31} (2010) 793--818.

\bibitem{bsz-1} J. L. Bona, S.-M. Sun and B.-Y. Zhang,
A nonhomogeneous boundary-value problem for the Korteweg-de Vries
equation in a quarter plane, \emph{Trans.  American Math. Soc.} {\bf
354 }(2001) 427--490.

\bibitem{bsz-2} J. L. Bona, S.-M. Sun and B.-Y. Zhang,
A nonhomogeneous boundary-value problem for the Korteweg-de Vries
equation in a bounded domain, \emph{Comm. Partial Diff.
Equations} {\bf 28} (2003) 1391--1436.

\bibitem{bsz-3} J. L. Bona, S.-M. Sun and B.-Y. Zhang, Forced oscillations of a damped
KdV equation in a quarter plane, \emph{Comm. Contemp. Math.}
\textbf{5} (2003) 369--400.

\bibitem{bsz-4} J. L. Bona, S.-M. Sun and B.-Y. Zhang,  Conditional and unconditional
well-posedness of nonlinear evolution equations, \emph{Advances
Diff. Equations}  {\bf 9} (2004) 241--265.

\bibitem{bsz-5}  J. L. Bona, S.-M. Sun and B.-Y. Zhang,
 Boundary smoothing properties of the Korteweg-de Vries
equation in a quarter plane and applications,\emph{ Dynamics Partial
Diff. Equations} \textbf{3} (2006) 1--70.

\bibitem{bsz-6}  J. L. Bona, S.-M. Sun and B.-Y. Zhang,
 Nonhomogeneous problems for the Korteweg-de Vries and
the Korteweg-de Vries-Burgers equations in a quarter plane,  \emph{Ann.
Inst. H. Poincar\'e, Anal. Non Lin\'eaire} \textbf{25}\,(2008)
1145--1185.

\bibitem{bsz-finite}  J. L. Bona, S.-M. Sun and B.-Y. Zhang,
 Nonhomogeneous problems for the Korteweg-de Vries  equation in
a bounded domain II, \emph{J. Diff. Equations}, \textbf{247} (2009)
4129--4153.

\bibitem{BW1983}  J. L. Bona and R. Winther, The Korteweg-de Vries equation posed in a quarter plane, \emph{SIAM
J. Math. Anal.} \textbf{14} (1983) 1056--1106.

\bibitem{bourgain-1}  J. Bourgain, Fourier transform restriction phenomena for certain
lattice subsets and applications to non-linear evolution equations,
part I: Schr\"{o}dinger equations, {\em Geom. \&  Funct. Anal. }
{\bf 3} (1993) 107--156.

\bibitem{bounrgain-k} J. Bourgain,  Fourier transform restriction phenomena for certain
lattice subsets and applications to non-linear evolution equations,
part II: the KdV  equation, {\em Geom. \&  Funct. Anal. } {\bf
3} (1993) 209--262.

\bibitem{bourgain-2} J. Bourgain, {\em Global Solutions of Nonlinear
Schr\"odinger Equations}, Colloqium Publication, Vol. 46, American
Mathematical Society, Providence, RI, 1999.


\bibitem{brezis} H. Br\'ezis and T. Gallouet, Nonlinear Schr\"odinger evolution equation,
\emph{Nonlinear Anal. TMA} \textbf{4} (1980) 677--681.

\bibitem{bu1994} C. Bu,  An initial-boundary value problem of the nonlinear Schr\"odinger equation,
 \emph{Appl. Anal.}  \textbf{53}  (1994) 241--254.

\bibitem{bu2000} C. Bu, Nonlinear Schr\"odinger equation on the semi-infinite line, \emph{Chinese Annals of Math.}
\textbf{21} (2000) 1-12.

\bibitem{bu2001} C. Bu, R. Shull, H. Wang, and M. Chu, Well-posedness, decay
estimates and blow-up theorem for the forced NLS,\emph{ J. Partial
Diff. Equations} \textbf{14} (2001) 61--70.

\bibitem{bu2005} C. Bu, K. Tsutaya and C Zhang, Nonlinear Schr\"odinger equation with
inhomogebeous Dirichlet boundary data, \emph{J. Math. Phys.}
\textbf{46} (2005) 083504.

\bibitem{bub-1} B. A. Bubnov,  Generalized boundary value problems for the
Korteweg-de Vries equation in bounded domain, \emph{Diff. Equations}
\textbf{15} (1979) 17--21.

\bibitem{bub-2} B. A. Bubnov,  Solvability in the large of nonlinear
boundary-value problems for the Korteweg-de Vries equations,
\emph{Diff. Equations} \textbf{16 } (1980) 24--30.



\bibitem{cazenave} T. Cazenave, \emph{Semilinear Schr\"odinger Equations}, American
 Math. Soc., Providence, RI, 2003.

\bibitem{cfh2011} T. Cazenave, D. Fang and Z. Han, Continuous dependence for NLS in fractional order spaces, \emph{Ann. Inst. H. Poincar\'e, Anal. Non Lin\'eaire} \textbf{28} (2011) 135--147.


\bibitem{ch} T. Cazenave and A. Haraux,
\emph{Introduction aux probl\`emes d'\'evolution semi-lin\'eaires},
Math\'ematiques et Applications, 1, Ellipses, Paris, 1990.

\bibitem{caz-weiss} T. Cazenave and F. B. Weissler, The Cauchy problem for the critical nonlinear Schr\"odinger
equation in $H^s$, \emph{ Nonlinear Anal. TMA} \textbf{14} (1990)
807--836.

 \bibitem{bc} R. Carroll and C. Bu, Solution of the forced nonlinear Schr\"odinger equation (NLS)
 using PDE techniques, \emph{Appl. Anal.}
{\bf 41}\,(1991) 33--51.

\bibitem{CHA2011}
A. Chabchoub, N. P. Hoffmann, and N. Akhmediev, Rogue wave observation in a water wave tank, \emph{Phys. Rev. Lett.} \textbf{106} (2011) 204502.

\bibitem{colin} T. Colin and J.-M. Ghidaglia, An initial-boundary-value problem for the
Korteweg-de Vries equation posed on a finite interval,
\emph{Advances Diff. Equations } \textbf{6} (2001) 1463--1492.

\bibitem{colliander} J. E. Colliander and C. E. Kenig, The generalized Korteweg-de Vries equation on the
half line, \emph{Comm. Partial Diff. Equations} \textbf{27}
(2002) 2187--2266.

\bibitem{Crai1985}
W. Craig, An existence theory for water waves and the Boussinesq and Korteweg-de Vries scaling limits, \emph{Comm. Partial Diff. Equations} \textbf{10} (1985) 787--1003.

\bibitem{fami-1} A. V. Faminskii, The Cauchy problem and the mixed problem in the
half strip for equations of Korteweg- de Vries type, (Russian)
\emph{Dinamika Sploshn. Sredy} \textbf{162} (1983) 152--158.

\bibitem{fami-2} A. V. Faminskii, A mixed problem in a semistrip for the
Korteweg-de Vries equation and its generalizations, (Russian)
\emph{Dinamika Sploshn. Sredy} \textbf{258} (1988) 54--94; English
transl. in \emph{Trans. Moscow Math. Soc}. \textbf{51} (1989)
53--91.

\bibitem{fami-3} A. V. Faminskii, Mixed problems for the Korteweg-de Vries
equation, \emph{Sbornik: Mathematics} \textbf{190 }(1999) 903--935.

\bibitem{fami-4}  A. V. Faminskii, An initial boundary-value problem in a
half-strip for the Korteweg-de Vries equation in fractional-order
Sobolev spaces, \emph{Comm. Partial Diff. Equations}  \textbf{29}
(2004) 1653--1695.

\bibitem{gini-velo-2} J. Ginibre and G. Velo, On a class of nonlinear Schr\"odinger
equations. I. The Cauchy problem, general case, \emph{J. Functional
Anal.} \textbf{32 }(1979) 1--32.

\bibitem{gini-velo}  J. Ginibre and G. Velo, On a class of nonlinear Schr\"odinger equations.
II.  Scattering theory, general case, \emph{J. Functinal Anal.}
\textbf{32} (1979) 33--71.


\bibitem{hadamard}  J. Hadamard,  Sur les probl\`emes aux d\'eriv\'ees partielles et leur signifiication physique,
\emph{Princeton University Bulletin}  (1902) 49--52.


\bibitem{hadamard1}
J. Hadamard,
\emph{Le\c{c}ons sur les propagation des ondes et les equations de l'hydrodynamics}, Hermann, Paris, 1903.


\bibitem{h-f}  Z. Han and D. Fan,  On the unconditional uniqueness for NLS in $H^s$, {\em SIAM J. Math Anal.} {\bf 45} (2013) 1505--1526.

\bibitem{hardy} G. Hardy, J. E. Littlewood and G. Polya,
\emph{Inequalities,}
University Press, Cambridge, England, 1988.


\bibitem{holmer} J. Holmer, The initial-boundary value problem for the $1$-$d$ nonlinear Schr\"odinger
equation on the half-line, \emph{Diff. Integral Equations}
\textbf{18} (2005)  647--668.

\bibitem{holmer-2} J. Holmer, The initial-boundary value problem for the Korteweg-de Vries equation,
\emph{Commun. Partial Diff. Equations} \textbf{31} (2006)
1151--1190.

\bibitem{Ivanovici-1}
O. Ivanovici, G. Lebeau and F. Planchon, Dispersion for the wave equation inside strictly convex
domains I: the Friedlander model case, \emph{Ann. of Math. (2)} \textbf{180} (2014) 323--380.

\bibitem{Ivanovici-2}
O. Ivanovici and F. Planchon, On the energy critical Schr\"odinger equation in 3D, non-trapping domains,
\emph{Ann. Inst. H. Poincar\'e, Anal. Non Lin\'eaire}  \textbf{27} (2010) 1153--1177.

\bibitem{ILT-1}
R. Illner, H. Lange and H. Teismann,  Limitations on the control of Schr\"odinger equations,
{\em  ESAIM Control Optim. Calc. Var.} {\bf 12} (2006) 615--635.

\bibitem{ILT-2} R. Illner, H. Lange and H. Teismann,  A note on the exact
internal control of nonlinear Schr\"odinger equations, {\em  CRM
Proceedings and Lecture Notes}  {\bf 33} (2003) 127--137.

\bibitem{JZ} A. Jeffrey and D. Zwillinger, \emph{Table of Integrals, Series and
Products}, 7th edition, Academic Press, San Diego, CA, 2007.

\bibitem{kam} S. Kamvissis, Semiclassical nonlinear Schr\"odinger
on the half line,  \emph{J. Math. Phys}. \textbf{ 44 } (2003)
5849--5868.

\bibitem{kato-1} T. Kato, On nonlinear Schr\"odinger equations, \emph{Ann. Inst. H. Poincar\'e,
 Phys. Theor.} \textbf{46}
(1987) 113--129.

\bibitem{kato-2} T. Kato, On nonlinear Scrh\"odinger equations. II. $H^s$-solutions and unconditional
well-posedness, \emph{J. d'Analyse Math.} \textbf{67} (1995)
281--306.

\bibitem{Killip}
R. Killip, M. Visan and X. Zhang, Quintic NLS in the exterior of a strictly convex obstacle, \emph{Amer. J. Math.}, to appear.

\bibitem{irena-1} I. Lasiecka and R. Triggiani, Optimal regularity, exact
controllability and uniform stabilization of Schrodinger equations
with Dirichlet control, \emph{Diff. Integral Equations
}\textbf{5} (1992) 521--535.

\bibitem{irena-2} I. Lasiecka, R. Triggiani and X. Zhang, Carleman estimates at
the $H^1(\Omega)$-- and $L^2(\Omega )-$level for nonconservative
Schr\"odinger equations with unobserved Neumann B.C., \emph{Arch.
Inequal. Appl.} \textbf{2} (2-3) (2004) 215--338.

\bibitem{LP2009}
F. Linares and G. Ponce, {\em Introduction to nonlinear dispersive equations}, Universitext, Springer, New York, 2009.

\bibitem{LM1972}  J. L. Lions and E. Magenes,  {\em Non-Homogeneous Boundary Value Problems and Applications, Vol. 1}, Springer-Verlag, Berlin-Heidelberg-New York, 1972.

\bibitem{LM1972-2}  J. L. Lions and E. Magenes,  {\em Non-Homogeneous Boundary Value Problems and Applications, Vol. 2}, Springer-Verlag, Berlin-Heidelberg-New York, 1972.

\bibitem{24} A. Pazy,  {\em Semigroups of Linear Operators and Applications
to Partial Differential Equations}, Applied Mathematical Sciences,
Vol. 44, Springer-Verlag, New York-Berlin-Heidelberg-Tokyo,  1983.


\bibitem{pere1983}
D. H. Peregrine,  Water waves, nonlinear Schr\"odinger equations and their solutions, \emph{J. Australian Math. Soc. B} \textbf{25} (1983) 16--43.

\bibitem{RZ2007b} L. Rosier and  B.-Y. Zhang,   Exact
controllability and stabilization of the nonlinear  Schr\"odinger
equation on a bounded interval, \emph{SIAM J. Control Optim.} {\bf
48}\, (2009) 972--992.

\bibitem{RZ2008} L. Rosier and  B.-Y. Zhang, Exact boundary
controllability of the nonlinear Schr\"odinger equation, \emph{J.
Diff. Equations} \textbf{246} (2009) 4129--4153.

\bibitem{RZ2009} L. Rosier and  B.-Y. Zhang, Control and stabilization
of the nonlinear Schr\"odinger equation on rectangles,  \emph{Math.
Models Methods Appl. Sci.}\textbf{ 20} (2010) 2293--2347.

\bibitem{SW2000}
G. Schneider and C. E. Wayne, The long-wave limit for the water wave problem. I. The case of zero surface tension, \emph{Comm. Pure Appl. Math.} \textbf{53} (2000) 1475--1535.

\bibitem{stein} E. M. Stein, \emph{Singular Integrals and Differentiability
of Functions}, Princeton University Press, Princeton, NJ, 1970.

\bibitem{strauss} W. Strauss and C. Bu, Inhomogeneous boundary value problem for a nonlinear
Schr\"odinger equation, \emph{J. Diff. Equations} \textbf{173}\,(2001)
79--91.

\bibitem{tsu-1} M. Tsutsumi, On smooth solutions to the initial-boundary value
problem for the nonlinear Schr\"odinger equations in two space
dimensions, \emph{Nonlinear Anal. TMA} \textbf{13 }(1989) 1051--1056.

\bibitem{tsu-2}  M. Tsutsumi, On global solutions to the initial-boundary value
problem for the nonlinear Schr\"odinger equations in exterior
demains, \emph{Comm. Partial Diff. Equations} \textbf{16}
(1991) 885--907.

\bibitem{tsu-3} Y. Tsutsumi, Global solutions of the nonlinear Schr\"odinger
equations in exterior domains, \emph{Comm. Partial Diff.
Equations} \textbf{8} (1983) 1337--1374.


\bibitem{tsu-4} Y. Tsutsumi, $L^2$-solutions for nonlinear Schr\"odinger equations
and nonlinear groups, \emph{Funk. Ekva.} \textbf{30} (1987)
115--125.

\bibitem{win-tsu} Y. Y. S. Win and Y. Tsutsumi, Unconditional uniqueness of solution for the Cauchy problem
of the nonlinear Schr\"odinger equation, {\em Hokkaido Math. J.}  {\bf 37}  (2008) 839--859.

\bibitem{zm1974}
V.E. Zakharov and S.V. Manakov, On the complete integrability of a nonlinear Schr\"odinger equation, \emph{J. Theore. and Math. Phys.} \textbf{19}  (1974) 551--559.

\bibitem{zs1972}
V. E. Zakharov and A. B. Shabat, Exact theory of two-dimensional self-focusing and
one-dimensional self-modulation of waves in nonlinear media, \emph{J. Experi. and Theore. Phys.} \textbf{34} (1972) 62--69.

\end{thebibliography}
\end{document}